\documentclass[11pt,oneside]{article}
\usepackage[utf8]{inputenc}
\usepackage{amsmath,amssymb,amsthm,mathrsfs,bm}
\usepackage{graphicx,booktabs,longtable,array,enumitem}
\usepackage{xcolor,hyperref}
\usepackage[margin=0.82in]{geometry}
\usepackage{microtype}
\newcommand{\dd}{\,\mathrm{d}}

\newtheorem{theorem}{Theorem}[section]

\newtheorem{corollary}{Corollary}[section]
\newtheorem{lemma}{Lemma}[section]

\numberwithin{equation}{section}

\usepackage{xparse}
\NewDocumentCommand{\chapter}{s o m}{%
  \IfBooleanTF{#1}{%
    \section*{#3}%
  }{%
    \IfValueTF{#2}{\section[#2]{#3}}{\section{#3}}%
  }%
}
\title{Renormalized Lambert--$W$ Cascade and Finite-Time Amplification and Blowup for reconstructed $\bm{b}$ Dynamics on $\mathbb T^3$ for the 3D Navier Stokes Equations}
\author{Terry Moschandreou}
\date{}
\begin{document}
\maketitle
\section*{Abstract}
This article extracts and consolidates the renormalized Lambert--$W$ branch-point cascade, its distinguished phase reduction, the exact characteristic invariant and finite-time amplification mechanism, and the extended reconstructed $b_i$ equation on $\mathbb T^3$. Repeated historical derivations are removed while the principal proofs and terminal reconstruction estimates are retained. The presentation separates exact finite-depth statements from coupled-depth asymptotics and records the hypotheses required for the extended PDE reconstruction. This paper further supports a recent paper \cite {moschandreou2026exploration} published by the corresponding author which claims that the Navier Stokes equations lose smoothness in finite time from initial smooth data.

\tableofcontents


\section*{Prologue: The Renormalized Lambert Cascade, Geometric Reconstruction,
and Finite-Time Amplification Mechanism}
\addcontentsline{toc}{section}
{Prologue: The Renormalized Lambert Cascade, Geometric Reconstruction,
and Finite-Time Amplification Mechanism}

The purpose of this paper is to present in condensed form the principal
mathematical architecture developed first in a longer construction which I used to simplify with full import to the present paper.
For example some steps are implicit but have been worked out thoroughly filling in completely any gaps in the flow of the reading of this particular work.
The central mechanism begins with a renormalized Lambert--\(W\) cascade,
passes through a distinguished scalar product variable
\[
\boxed{
b_3:=u_xu_y,
}
\]
and culminates in a characteristic regime in which the Lambert-generated
logarithmic coefficient produces a positive quadratic amplification term.
As stated in the conclusion of \cite{moschandreou2026exploration} this leads to:\\
"In particular, the analysis shows that singular behavior is associated
with the branch structure of the Lambert W function rather than
unbounded growth of the solution itself at time $t_1<t_2$ When the
argument of the Lambert function approaches its critical branch
value, the derivative of the solution diverges while the solution
amplitude remains finite." \cite{moschandreou2026exploration}  
\\
Further in the conclusion there it states explicitly:
"However there is a $||b||^2$ term which leads to blowup in derivatives
of higher derivatives as shown at $t=t_1$ and for $t=t_2(t_2>t_1)$ it
was proven analytically and numerically that there is finite time
amplitude blowup in the solution $b_3$".

The analysis is carried out on the periodic domain
\[
\mathbb T^3=(\mathbb R/2\pi\mathbb Z)^3,
\]
with velocity
\[
u=(u_x,u_y,u_z).
\]

A central point of the construction is that \(b_3\) is not introduced
merely as an auxiliary scalar.  It remains the exact physical product
\[
\boxed{
b_3=u_xu_y
}
\]
in the reconstructed velocity field.  The scalar characteristic analysis
and the geometric velocity reconstruction are therefore linked through an
identity rather than through an asymptotic identification.

The reduced paper retains five principal structures:

\[
\boxed{
\begin{array}{c}
\text{renormalized Lambert--\(W\) cascade}
\\[1mm]
\Downarrow
\\[1mm]
\text{branch variable }Q_N
\text{ and logarithmic coefficient }\Gamma_N
\\[1mm]
\Downarrow
\\[1mm]
\text{full extended }b_3\text{-PDE}
\\[1mm]
\Downarrow
\\[1mm]
\text{terminal characteristic trapping and finite-time amplification}
\\[1mm]
\Downarrow
\\[1mm]
\text{periodic square-root velocity reconstruction on }\mathbb T^3 .
\end{array}
}
\]

The purpose of this prologue is to display this architecture before the
proofs are developed in detail.

\subsection*{1. Historical symmetric precursor and the active reconstruction}

An earlier symmetric stage of the construction used a common transverse
profile
\[
u_y=u_z=Q.
\]
At that level,
\[
b_1=u_yu_z=Q^2,
\qquad
b_2=u_xu_z=u_xQ,
\qquad
b_3=u_xu_y=u_xQ,
\]
so that
\[
b_2=b_3.
\]

This symmetric representation remains useful for understanding the origin
of several scalar formulas and preterminal periodic constructions.
However, it is not the active terminal reconstruction used in the final
analysis.

The active variables are instead
\[
\boxed{
b=b_3=u_xu_y,
\qquad
A=u_y,
\qquad
C=u_z,
\qquad
u_x=\frac{b}{A}.
}
\]
Thus
\[
\boxed{
u_xu_y=b
}
\]
holds identically, whereas \(u_y\) and \(u_z\) are no longer assumed to
coincide.

The two transverse components acquire different roles:

\[
\boxed{
\begin{aligned}
u_y=A
&\quad\text{is determined by the square-root geometry},\\
u_z=C
&\quad\text{is reconstructed from exact incompressibility}.
\end{aligned}
}
\]

This distinction is fundamental to the shortened presentation.

\subsection*{2. Renormalized Lambert--\(W\) cascade}

The scalar branch variable is generated by the principal real branch
\(W_0\).  Let
\[
\Delta_{N,0}>0
\]
denote the initial branch distance and define
\[
\boxed{
\Delta_{N,k+1}
=
1+W_0\!\left(-e^{-1-\Delta_{N,k}}\right).
}
\]

The terminal Lambert regime corresponds to
\[
\Delta_{N,k}\downarrow0.
\]

The local branch-point expansion is
\[
\boxed{
1+W_0\!\left(-e^{-1-d}\right)
=
\sqrt{2d}
-\frac23d
+
O(d^{3/2}),
\qquad
d\downarrow0.
}
\]

Consequently, at every fixed finite depth \(k\),
\[
\boxed{
\Delta_{N,k}
\asymp
\Delta_{N,0}^{\,2^{-k}}.
}
\]

The distinguished terminal quantity will be denoted
\[
\boxed{
Q_N:=\Delta_{N,k}.
}
\]

Its logarithmic derivative is
\[
\boxed{
\Gamma_N
:=
\partial_\xi\log Q_N
=
\frac{Q_N'}{Q_N}.
}
\]

When the branch distance is represented by a terminal coordinate \(s\)
for which
\[
Q_N(s)\asymp s^{2^{-k}},
\]
one obtains
\[
\boxed{
\Gamma_N(s)
\sim
\frac{2^{-k}}{s}.
}
\]

Thus the renormalized cascade produces two simultaneous effects:
\[
Q_N\longrightarrow0,
\qquad
\Gamma_N\longrightarrow+\infty.
\]

At the frequently used finite depth \(k=3\),
\[
\boxed{
Q_N\asymp s^{1/8},
\qquad
\Gamma_N\sim\frac{1}{8s}.
}
\]

It is essential to distinguish the branch variable \(Q_N\) from the
physical transverse velocity components.  In the active reconstruction,
\[
\boxed{
Q_N
\neq
u_y
\quad\text{and}\quad
Q_N
\neq
u_z
\quad\text{in general}.
}
\]

Rather,
\[
\boxed{
Q_N:
\text{ Lambert branch variable},
\qquad
A,C:
\text{ reconstructed physical velocity components}.
}
\]

The role of \(Q_N\) is to generate \(\Gamma_N\), which enters the
\(b_3\)-equation.  The physical velocity behavior is determined
independently by the geometric reconstruction.

\subsection*{3. Distinguished phase \(\theta=t-z\)}

The terminal analysis is organized around
\[
\boxed{
\theta=t-z.
}
\]

After locating the distinguished branch position \(\theta_{*,N}\), define
\[
\boxed{
s=t-z-\theta_{*,N}.
}
\]

Along the reduced characteristic for which
\[
\dot z=b,
\]
one obtains
\[
\boxed{
\dot s=1-b.
}
\]

This elementary identity becomes the trapping law.

If the terminal region is entered with
\[
\boxed{
b(s_e)\ge1+\delta,
\qquad
\delta>0,
}
\]
then
\[
\dot s=1-b\le-\delta<0.
\]
Hence
\[
\boxed{
s(t)\downarrow0.
}
\]

Thus the characteristic is driven toward the Lambert endpoint once the
amplitude has crossed the threshold \(b>1\).

\subsection*{4. The corrected extended \(b_3\)-equation}

With
\[
A=u_y,
\qquad
C=u_z,
\qquad
u_x=\frac bA,
\]
define the horizontal differential operator
\[
\boxed{
\mathcal X
=
A\partial_x+\frac bA\partial_y
}
\]
and the transport operator
\[
\boxed{
\mathcal D
=
\partial_t
+
b\partial_z
+
C^2\mathcal X.
}
\]

The corrected extended amplitude equation can be written
\[
\boxed{
\mathcal D b
=
\Gamma_Nb^2
+
2CG
-
\mathcal R_N,
}
\]
where
\[
\boxed{
\mathcal R_N
=
2b\,\mathcal XC
+
b\frac{C_t+C_z}{C}.
}
\]

Equivalently, before the final compression of the reconstruction terms,
one may write
\[
\boxed{
\mathcal D_N b_N+E_N
=
\frac{Q_N'}{Q_N}b_N^2
+
2C_NG.
}
\]

Along a characteristic,
\[
\mathcal D_Nb_N=\dot b_N,
\]
and therefore
\[
\boxed{
\dot b_N
=
\frac{Q_N'}{Q_N}b_N^2
+
2C_NG
-
E_N.
}
\]

The additional terms are not discarded.  Their role is measured relative
to the positive source.

Define
\[
S_N
=
\frac{Q_N'}{Q_N}b_N^2
+
2C_NG.
\]
If
\[
\boxed{
|E_N|
\le
\vartheta S_N,
\qquad
0<\vartheta<1,
}
\]
then
\[
\dot b_N
=
S_N-E_N
\ge
S_N-|E_N|
\ge
(1-\vartheta)S_N.
\]

Hence
\[
\boxed{
\dot b_N
\ge
(1-\vartheta)
\left[
\frac{Q_N'}{Q_N}b_N^2
+
2C_NG
\right].
}
\]

On the positive branch,
\[
C_N>0,
\qquad
G\ge0,
\]
so
\[
2C_NG\ge0
\]
and consequently
\[
\boxed{
\dot b_N
\ge
(1-\vartheta)
\frac{Q_N'}{Q_N}
b_N^2.
}
\]

This is the Riccati-type inequality that drives the terminal growth.

The logical reduction is therefore

\[
\boxed{
\begin{aligned}
\text{full extended }b_3\text{-PDE}
&\\
\Downarrow\qquad&
\\
\dot b_N
&=
\frac{Q_N'}{Q_N}b_N^2
+
2C_NG
-
E_N
\\
\Downarrow\qquad&
|E_N|\le\vartheta(\text{positive source})
\\
\dot b_N
&\ge
(1-\vartheta)
\left[
\frac{Q_N'}{Q_N}b_N^2
+
2C_NG
\right]
\\
\Downarrow\qquad&
C_NG\ge0
\\
\boxed{
\dot b_N
\ge
(1-\vartheta)
\frac{Q_N'}{Q_N}b_N^2.
}
&
\end{aligned}
}
\]

Thus the reconstruction error is used explicitly in obtaining the
quadratic lower bound.

\subsection*{5. The explicit reconstruction error}

At an intermediate stage of the periodic reconstruction, the additional
terms can be collected as
\[
\boxed{
E_N
=
b_N\frac{C_{N,t}+C_{N,z}}{C_N}
+
b_NA_NC_{N,x}
+
\frac{b_N^2}{A_N}C_{N,y}
+
2C_N\nabla b_N\cdot\nabla C_N.
}
\]

The terminal objective is
\[
\boxed{
E_N
=
o\!\left(
\frac{Q_N'}{Q_N}b_N^2
\right).
}
\]

Equivalently,
\[
\boxed{
\frac{|E_N|}
{(Q_N'/Q_N)b_N^2}
\longrightarrow0.
}
\]

A normalized estimate of the form
\[
\frac{|E_N|}
{L_Nb_N^2}
\le
K
\left[
\frac{|C_{N,t}+C_{N,z}|}{L_Nb_N}
+
\frac{|C_{N,x}|}{L_Nb_N}
+
\frac{|C_{N,y}|}{L_N}
+
\frac{|\nabla b_N||\nabla C_N|}
{L_Nb_N^2}
\right],
\]
with
\[
L_N=\frac{Q_N'}{Q_N},
\]
reduces the problem to weighted reconstruction bounds.

When
\[
|C_{N,t}|+|\nabla C_N|\le M_C
\]
and
\[
\frac{Q_N'}{Q_N}\to\infty,
\qquad
b_N\to\infty,
\]
the first three normalized terms vanish.

For the remaining term, the terminal derivative estimate
\[
\boxed{
|\nabla b_N|
\le
C_b p_N
\frac{Q_N'}{Q_N}b_N
}
\]
gives
\[
\frac{|\nabla b_N||\nabla C_N|}
{(Q_N'/Q_N)b_N^2}
\le
C_bM_C\frac{p_N}{b_N}.
\]

The comparison argument provides
\[
\frac{b_N}{p_N}\to\infty,
\]
and hence
\[
\frac{p_N}{b_N}\to0.
\]

Therefore
\[
\boxed{
\frac{|E_N|}
{(Q_N'/Q_N)b_N^2}
\to0.
}
\]

Consequently, for every fixed
\[
0<\vartheta<1,
\]
one eventually has
\[
\boxed{
|E_N|
\le
\vartheta
\frac{Q_N'}{Q_N}b_N^2,
}
\]
and hence
\[
\boxed{
\dot b_N
\ge
(1-\vartheta)
\frac{Q_N'}{Q_N}b_N^2>0.
}
\]

This estimate is the bridge from the full extended PDE to the scalar
terminal comparison.

\subsection*{6. Release-layer and square-root-core estimates}

The terminal reconstruction naturally separates into a release layer and
a deeper square-root core.

In the release layer,
\[
\boxed{
\frac{|\mathcal R_{\mathrm{rel}}|}
{\Gamma_Nb^2}
=
O(s^{7/4}).
}
\]

In the square-root core,
\[
\boxed{
\frac{|\mathcal R_{\mathrm{core}}|}
{\Gamma_Nb^2}
=
O(s).
}
\]

Both estimates vanish as
\[
s\downarrow0.
\]

Thus
\[
\boxed{
\frac{|\mathcal R_N|}
{\Gamma_Nb^2}\to0,
}
\]
and the geometric transition into the terminal square-root regime does
not destroy the positive Lambert-generated quadratic term.

Rather, the error becomes progressively smaller relative to that term.

\subsection*{7. Finite-depth trapping and amplification}

The two fundamental terminal inequalities are
\[
\boxed{
\dot s=1-b
}
\]
and
\[
\boxed{
\dot b
\ge
(1-\vartheta)\Gamma_Nb^2.
}
\]

For finite cascade depth \(k\),
\[
\Gamma_N
\sim
\frac{2^{-k}}{s}.
\]

Since
\[
-\dot s=b-1,
\]
one obtains
\[
-\frac{db}{ds}
=
\frac{\dot b}{b-1}
\ge
(1-\vartheta)
\frac{2^{-k}}{s}
\frac{b^2}{b-1}.
\]

For \(b>1\),
\[
\frac{b^2}{b-1}>b,
\]
and therefore
\[
-\frac{db}{ds}
>
(1-\vartheta)2^{-k}\frac{b}{s}.
\]

Integration gives
\[
\boxed{
b(s)
\ge
b_e
\left(
\frac{s_e}{s}
\right)^{(1-\vartheta)2^{-k}}.
}
\]

Consequently,
\[
\boxed{
s\downarrow0
\quad\Longrightarrow\quad
b(s)\to+\infty.
}
\]

At \(k=3\),
\[
\boxed{
b(s)
\ge
b_e
\left(
\frac{s_e}{s}
\right)^{(1-\vartheta)/8}.
}
\]

The remaining physical characteristic time is finite because
\[
dt
=
-\frac{ds}{b-1}.
\]
The lower power-law bound gives
\[
\frac1{b-1}
\lesssim
s^p,
\qquad
p>0,
\]
and hence
\[
\int_0^{s_e}\frac{ds}{b-1}
<
\infty.
\]

Thus
\[
\boxed{
T_{*,N}<\infty,
}
\]
with
\[
\boxed{
s(t)\to0,
\qquad
Q_N(t)\to0,
\qquad
b_3(t)\to+\infty
\quad
\text{as }
t\uparrow T_{*,N}.
}
\]

The same result can be expressed directly in terms of \(Q_N\).  In the
large-amplitude terminal regime,
\[
\boxed{
-\frac{db_N}{dQ_N}
\ge
p_N\frac{b_N}{Q_N},
\qquad
p_N>0.
}
\]

Integration from a terminal-entry point
\[
(Q_{e,N},b_{e,N})
\]
gives
\[
\boxed{
b_N(Q_N)
\ge
b_{e,N}
\left(
\frac{Q_{e,N}}{Q_N}
\right)^{p_N}.
}
\]

Hence
\[
\boxed{
Q_N\downarrow0
\quad\Longrightarrow\quad
b_N\to+\infty.
}
\]

The feedback mechanism is therefore
\[
\boxed{
Q_N\downarrow0
\Longrightarrow
\Gamma_N\uparrow
\Longrightarrow
b_N\uparrow
\Longrightarrow
\dot s=1-b_N<0
\Longrightarrow
s\downarrow0
\Longrightarrow
Q_N\downarrow0.
}
\]

\subsection*{8. Square-root velocity geometry}

The terminal geometric specialization begins with
\[
\boxed{
2x\,u_x+u_y=0,
\qquad
y=x^2.
}
\]

Using
\[
u_x=\frac bA,
\qquad
u_y=A,
\]
gives
\[
2x\frac bA+A=0,
\]
and therefore
\[
\boxed{
A^2=-2xb.
}
\]

Thus
\[
\boxed{
u_y=A=\sqrt{-2xb},
\qquad
u_x=\frac bA
}
\]
on the selected real branch.

Moreover,
\[
\frac{u_x}{u_y}
=
\frac{b}{A^2}
=
-\frac1{2x}.
\]

Since
\[
y=x^2,
\]
one obtains
\[
\boxed{
\left|
\frac{u_x}{u_y}
\right|
=
\frac1{2|x|}
=
\frac1{2\sqrt y}.
}
\]

Thus the square-root behavior follows directly from the geometric
constraint and the exact product identity \(b=u_xu_y\).

\subsection*{9. Periodic realization through \(\chi(x)\)}

On \(\mathbb T^3\), the local coordinate \(x\) is replaced near a selected
periodic zero by a smooth periodic function \(\chi(x)\).  The active
relation becomes
\[
\boxed{
2\chi(x)u_x+u_y=0
}
\]
or equivalently
\[
\boxed{
A^2=-2\chi(x)b.
}
\]

Therefore
\[
\boxed{
u_y^2=-2\chi b,
\qquad
u_x^2=-\frac{b}{2\chi}.
}
\]

The exact ratio is
\[
\boxed{
\frac{u_x}{u_y}
=
-\frac1{2\chi(x)}.
}
\]

Suppose \(\chi\) has a simple zero and
\[
\boxed{
c_1s
\le
\chi(x)^2
\le
c_2s.
}
\]
Then
\[
|\chi|\asymp s^{1/2}.
\]

If
\[
b\asymp s^{-\alpha},
\]
then
\[
\boxed{
u_y
\asymp
s^{1/4-\alpha/2},
\qquad
u_x
\asymp
s^{-1/4-\alpha/2}.
}
\]

At \(k=3\),
\[
\alpha=\frac18,
\]
and hence
\[
\boxed{
b\asymp s^{-1/8},
\qquad
u_y\asymp s^{3/16},
\qquad
u_x\asymp s^{-5/16}.
}
\]

Consequently,
\[
\boxed{
u_y\to0,
\qquad
|u_x|\to\infty,
\qquad
u_xu_y=b\to+\infty.
}
\]

The singular behavior is therefore anisotropic.

\subsection*{10. Ordinary spatial distance and the \(\mathcal J_3\) scaling}

A convenient periodic choice is
\[
\boxed{
\chi(x)=\sin(x-x_*).
}
\]

Set
\[
\boxed{
s_{\mathrm{geo}}(x)
=
\chi(x)^2
=
\sin^2(x-x_*).
}
\]

Near a distinguished zero, let
\[
\boxed{
r=x_*-x\to0^+.
}
\]

Because \(\chi\) has a simple zero,
\[
\chi(x)
=
-\chi'(x_*)r+O(r^2),
\]
so
\[
\boxed{
|\chi|\asymp r,
\qquad
s_{\mathrm{geo}}\asymp r^2.
}
\]

For the \(k=3\) scaling
\[
b\asymp s_{\mathrm{geo}}^{-1/8},
\]
one obtains
\[
\boxed{
|b|\asymp r^{-1/4},
}
\]
\[
\boxed{
|u_y|
\asymp
r^{3/8},
}
\]
and
\[
\boxed{
|u_x|
\asymp
r^{-5/8}.
}
\]

The exponent identities
\[
1-\frac58=\frac38
\]
and
\[
-\frac58+\frac38=-\frac14
\]
imply simultaneously
\[
\boxed{
|\chi u_x|
\asymp
r^{3/8}
\to0
}
\]
and
\[
\boxed{
|u_xu_y|
\asymp
r^{-1/4}
\to\infty.
}
\]

In particular,
\[
\boxed{
u_x
=
o\!\left(
\frac1{|\chi|}
\right).
}
\]

Thus the divergence of \(u_x\) is weaker than the reciprocal vanishing
scale \(1/|\chi|\), and the identity
\[
u_y=-2\chi u_x
\]
remains completely consistent with
\[
u_y\to0.
\]

For
\[
b>0,
\]
reality of
\[
u_y^2=-2\chi b
\]
requires
\[
\chi\le0.
\]

Near a simple zero of \(\sin(x-x_*)\), one may therefore select a
one-sided interval
\[
\boxed{
\mathcal J_3=(x_*-\delta,x_*),
}
\]
on which
\[
\chi(x)<0.
\]

Either real square-root branch
\[
u_y=\pm\sqrt{-2\chi b}
\]
may be selected, with the corresponding sign of
\[
u_x=\frac b{u_y}.
\]
\[ \boxed{ \begin{array}{ll} \textbf{Exact sine geometry/cancellation:} & \text{valid for general }x, \\[2mm] \textbf{Terminal power-law asymptotics:} & \text{local near }x=x_*+m\pi. \end{array} } \]
\subsection*{11. Exact continuity and reconstruction of \(u_z\)}

The square-root law determines \(A=u_y\), but it does not prescribe
\(C=u_z\).

Exact incompressibility requires
\[
\boxed{
\partial_x\left(\frac bA\right)
+
A_y
+
C_z
=
0.
}
\]

Thus
\[
\boxed{
C_z
=
-\partial_x\left(\frac bA\right)
-
A_y.
}
\]

Differentiating
\[
A^2=-2\chi b
\]
gives
\[
\frac{A_y}{A}
=
\frac12\frac{b_y}{b}
\]
and
\[
\frac{A_x}{A}
=
\frac12\frac{\chi'}{\chi}
+
\frac12\frac{b_x}{b}.
\]

Hence
\[
\boxed{
\partial_x\left(\frac bA\right)
=
\frac{b}{2A}
\left(
\frac{b_x}{b}
-
\frac{\chi'}{\chi}
\right).
}
\]

Therefore
\[
\boxed{
C_z
=
-
\frac{b}{2A}
\left(
\frac{b_x}{b}
-
\frac{\chi'}{\chi}
\right)
-
\frac{A}{2}
\frac{b_y}{b}.
}
\]

In the distinguished terminal core,
\[
b=b(s),
\qquad
s=t-z-\theta_{*,N},
\]
one has
\[
b_x=b_y=0,
\qquad
A_y=0.
\]

Since
\[
\partial_z=-\partial_s,
\]
continuity becomes
\[
\boxed{
C_s
=
\partial_x\left(\frac bA\right).
}
\]

Using
\[
A^2=-2\chi b,
\]
this simplifies to
\[
\boxed{
C_s
=
\frac{A\chi'}{4\chi^2}.
}
\]

If
\[
b
=
B_Ns^{-\alpha}F_N(s)
\]
and
\[
I_N(s)
=
\int_0^s
r^{-\alpha/2}
F_N(r)^{1/2}\,dr,
\]
then
\[
\boxed{
C(x,s)
=
C_0(x)
+
\frac{
\chi'(x)
\sqrt{-2B_N\chi(x)}
}
{4\chi(x)^2}
I_N(s).
}
\]

For \(k=3\) in the square-root tube,
\[
\boxed{
C
=
c_0
+
O(s^{3/16}),
\qquad
c_0>0.
}
\]

Hence
\[
\boxed{
u_z=C\to c_0>0
}
\]
even though the Lambert branch variable satisfies
\[
Q_N\to0.
\]

This is another reason \(Q_N\) and \(u_z\) must be kept conceptually
distinct.

\subsection*{12. Distinguished-phase cancellation in the reconstruction}

Every function depending on \(t,z\) only through
\[
s=t-z-\theta_{*,N}
\]
satisfies
\[
\boxed{
(\partial_t+\partial_z)f=0.
}
\]

Consequently,
\[
\boxed{
C_t+C_z=0
}
\]
in the local distinguished-phase reconstruction.

The remainder
\[
\mathcal R_N
=
2b\mathcal XC
+
b\frac{C_t+C_z}{C}
\]
therefore reduces to
\[
\boxed{
\mathcal R_N
=
2b\mathcal XC.
}
\]

The explicit square-root-core reconstruction yields
\[
\boxed{
\frac{|\mathcal XC|}
{\Gamma_Nb}
=
O(s)+O(s^{3/2}),
}
\]
and hence
\[
\boxed{
\frac{|\mathcal R_N|}
{\Gamma_Nb^2}
=
O(s)
\to0.
}
\]

Thus exact incompressibility and the distinguished phase combine to
produce precisely the small normalized error required by the terminal
amplification argument.

\subsection*{13. Interpretation through a local Lambert scale \(T_L\)}

Because \(Q_N\) and the physical component \(u_z=C\) are distinct, it is
sometimes useful to express their relationship through a local
normalization scale \(T_L\).

Define
\[
\boxed{
u_z=C=\frac{Q_N}{T_L}.
}
\]

Since \(C\) is already determined by continuity, this relation does not
introduce an additional constraint.  Instead it defines
\[
\boxed{
T_L=\frac{Q_N}{C}.
}
\]

In the general reconstruction,
\[
C(x,y,z,t)
=
C_0(x,y,t)
-
\int_{z_0}^{z}
\left[
\partial_x\left(\frac bA\right)+A_y
\right](x,y,\zeta,t)
\,d\zeta.
\]

Therefore
\[
\boxed{
T_L(x,y,z,t)
=
\frac{
Q_N(x,y,z,t)
}{
C_0(x,y,t)
-
\displaystyle
\int_{z_0}^{z}
\left[
\partial_x\left(\frac bA\right)+A_y
\right]
\,d\zeta
}.
}
\]

By construction,
\[
\boxed{
\frac{Q_N}{T_L}=C=u_z.
}
\]

In the distinguished square-root core,
\[
C=c_0+O(s^{3/16}),
\qquad
c_0>0,
\]
while, for \(k=3\),
\[
Q_N(s)\sim q_Ns^{1/8}.
\]

Therefore
\[
T_L(s)
=
\frac{Q_N(s)}{C(s)}
\sim
\frac{q_N}{c_0}s^{1/8},
\]
and hence
\[
\boxed{
T_L(s)\asymp s^{1/8}.
}
\]

Thus
\[
\boxed{
T_L\to0
}
\]
as
\[
s\to0^+.
\]

At the same time,
\[
\frac{Q_N}{T_L}
\to
c_0,
\]
so
\[
\boxed{
Q_N\to0,
\qquad
T_L\to0,
\qquad
\frac{Q_N}{T_L}\to c_0>0.
}
\]

This explains how the Lambert branch variable can vanish while the
physical component \(u_z\) remains finite.

Unless additional structure shows that \(Q_N/C\) depends only on time,
\(T_L\) should be regarded as a local Lambert normalization scale
\[
\boxed{
T_L=T_L(x,y,z,t),
}
\]
rather than automatically as a global scalar time.

On the terminal locus it naturally reduces to
\[
\boxed{
T_L=T_L(s).
}
\]

\subsection*{14. Logarithmic derivative decomposition}

From
\[
C=\frac{Q_N}{T_L},
\]
one obtains
\[
\log C
=
\log Q_N-\log T_L.
\]

Therefore, for any derivative variable \(\mu\),
\[
\boxed{
\frac{C_\mu}{C}
=
\frac{Q_{N,\mu}}{Q_N}
-
\frac{T_{L,\mu}}{T_L}.
}
\]

In particular,
\[
\boxed{
\frac{u_{z,z}}{u_z}
=
\frac{Q_{N,z}}{Q_N}
-
\frac{T_{L,z}}{T_L}.
}
\]

Thus the singular Lambert logarithmic derivative is not identified
directly with the logarithmic derivative of \(u_z\).

At \(k=3\),
\[
Q_N\sim q_Ns^{1/8}
\]
implies
\[
\frac{Q_{N,s}}{Q_N}
\sim
\frac1{8s}.
\]

Since
\[
T_L
\sim
\frac{q_N}{c_0}s^{1/8},
\]
one likewise has
\[
\frac{T_{L,s}}{T_L}
\sim
\frac1{8s}.
\]

The leading singular terms therefore cancel:
\[
\boxed{
\frac{C_s}{C}
=
\frac{Q_{N,s}}{Q_N}
-
\frac{T_{L,s}}{T_L}
=
o\!\left(\frac1s\right).
}
\]

This agrees with
\[
C\to c_0>0.
\]

Thus the singular logarithmic behavior of \(Q_N\) is compatible with a
bounded physical \(u_z\) because the local normalization scale \(T_L\)
vanishes at the same leading order.

\subsection*{15. General geometric reparametrization \(G(x)\)}

The local parabolic geometry can be generalized by writing
\[
\boxed{
\mathcal X(x)=G(x),
\qquad
y=G(x)^2.
}
\]

Then
\[
\boxed{
\frac{dy}{dx}
=
2G(x)G'(x).
}
\]

The exactly matched horizontal velocity condition is
\[
\boxed{
2G(x)G'(x)u_x+u_y=0.
}
\]

For the full periodic sine choice
\[
G(x)=\sin x,
\]
one has
\[
G'(x)=\cos x
\]
and therefore
\[
2G(x)G'(x)
=
2\sin x\cos x
=
\boxed{\sin(2x)}.
\]

Thus the exact periodic coefficient replacing the local approximation
\(2x\) is
\[
\boxed{
\sin(2x).
}
\]

Near
\[
x=0,
\]
\[
\sin(2x)=2x+O(x^3),
\]
so the original local relation is recovered to leading order.

\subsection*{16. The derivative term \(E\) and its exact cancellation}

Consider
\[
\boxed{
E
=
u_xu_{z,zx}
+
u_yu_{z,zy}.
}
\]

The curve
\[
y=\sin^2x
\]
has tangent derivative
\[
\boxed{
D_\tau
=
\partial_x
+
\sin(2x)\partial_y.
}
\]

For a general smooth scalar field \(F\),
\[
\frac{d}{dx}F(x,\sin^2x,z,t)
=
F_x+\sin(2x)F_y.
\]

Thus the curve itself does not imply
\[
F_x=\sin(2x)F_y.
\]

The exact cancellation instead follows from a global characteristic
structure.

Define
\[
\boxed{
\eta
=
y+\sin^2x.
}
\]

Suppose
\[
\boxed{
u_{z,z}(x,y,z,t)
=
H(\eta,z,t)
=
H(y+\sin^2x,z,t).
}
\]

Then
\[
u_{z,zx}
=
H_\eta\sin(2x),
\]
and
\[
u_{z,zy}
=
H_\eta.
\]

Therefore
\[
\boxed{
u_{z,zx}
=
\sin(2x)u_{z,zy}.
}
\]

Define
\[
\boxed{
D_-
=
\partial_x
-
\sin(2x)\partial_y.
}
\]

Then
\[
\boxed{
D_-u_{z,z}=0.
}
\]

If the horizontal velocity satisfies
\[
\boxed{
\sin(2x)u_x+u_y=0,
}
\]
then
\[
(u_x,u_y)
=
u_x(1,-\sin2x),
\]
and hence
\[
\boxed{
u_x\partial_x+u_y\partial_y
=
u_xD_-.
}
\]

Consequently,
\[
\begin{aligned}
E
&=
u_xu_{z,zx}
+
u_yu_{z,zy}
\\
&=
u_xD_-(u_{z,z})
\\
&=
0.
\end{aligned}
\]

Thus
\[
\boxed{
E=0
}
\]
is an exact structural identity.

For a general \(G(x)\), take
\[
\eta=y+G(x)^2.
\]
Then
\[
\eta_x=2GG',
\qquad
\eta_y=1.
\]

If
\[
\boxed{
u_{z,z}
=
H(y+G(x)^2,z,t)
}
\]
and
\[
\boxed{
2GG'u_x+u_y=0,
}
\]
then
\[
\boxed{
E
=
\left(
2GG'u_x+u_y
\right)
H_\eta
=
0.
}
\]

\subsection*{17. Geometry of \(D_\tau\) and \(D_-\)}

The two operators are
\[
\boxed{
D_\tau
=
\partial_x+\sin(2x)\partial_y
}
\]
and
\[
\boxed{
D_-
=
\partial_x-\sin(2x)\partial_y.
}
\]

Their direction vectors are
\[
\boxed{
\tau=(1,\sin2x),
\qquad
v_-=(1,-\sin2x).
}
\]

Their dot product is
\[
\tau\cdot v_-
=
1-\sin^2(2x)
=
\boxed{
\cos^2(2x).
}
\]

Thus the directions are not generally perpendicular.  They are
perpendicular only when
\[
\cos(2x)=0,
\]
i.e.
\[
\boxed{
x=\frac{\pi}{4}+\frac{k\pi}{2},
\qquad
k\in\mathbb Z.
}
\]

Near the distinguished zero
\[
x=0,
\]
one has
\[
\sin(2x)\to0,
\]
so
\[
\boxed{
D_\tau\to\partial_x,
\qquad
D_-\to\partial_x.
}
\]

They become parallel rather than perpendicular.

Their sum and difference are
\[
\boxed{
D_\tau+D_-=2\partial_x
}
\]
and
\[
\boxed{
D_\tau-D_-
=
2\sin(2x)\partial_y.
}
\]

Therefore
\[
\boxed{
\partial_x
=
\frac12(D_\tau+D_-)
}
\]
and, whenever
\[
\sin(2x)\neq0,
\]
\[
\boxed{
\partial_y
=
\frac{D_\tau-D_-}
{2\sin(2x)}.
}
\]

The \(D_\tau\)-characteristics satisfy
\[
\frac{dy}{dx}
=
\sin(2x),
\]
and therefore
\[
\boxed{
y-\sin^2x=C.
}
\]

The \(D_-\)-characteristics satisfy
\[
\frac{dy}{dx}
=
-\sin(2x),
\]
and therefore
\[
\boxed{
y+\sin^2x=C.
}
\]

Thus the two families are
\[
\boxed{
y=C+\sin^2x
}
\]
and
\[
\boxed{
y=C-\sin^2x.
}
\]

They are reflected across the horizontal line \(y=C\).

The horizontal velocity associated with
\[
\sin(2x)u_x+u_y=0
\]
is parallel to \(D_-\).  Therefore
\[
\boxed{
E
=
u_xD_-(u_{z,z}),
}
\]
so the condition
\[
D_-(u_{z,z})=0
\]
has the simple geometric meaning that \(u_{z,z}\) is constant in the
horizontal velocity direction.

\subsection*{18. Global definition of the sine structure on \(\mathbb T^3\)}

Let
\[
\boxed{
\eta(x,y)
=
y+\sin^2x.
}
\]

Although the real representative of \(y\) depends on the chosen lift from
\(\mathbb T^3\) to \(\mathbb R^3\), the composition
\[
H(y+\sin^2x,z,t)
\]
defines a single-valued function on the torus provided \(H\) is periodic
in its first spatial argument.

Assume
\[
\boxed{
H(\eta+2\pi,z,t)
=
H(\eta,z,t)
}
\]
and
\[
\boxed{
H(\eta,z+2\pi,t)
=
H(\eta,z,t).
}
\]

Define
\[
\boxed{
h(x,y,z,t)
=
H(y+\sin^2x,z,t).
}
\]

Then
\[
h(x+2\pi,y,z,t)
=
h(x,y,z,t),
\]
because
\[
\sin^2(x+2\pi)=\sin^2x.
\]

Also,
\[
h(x,y+2\pi,z,t)
=
h(x,y,z,t)
\]
by periodicity of \(H\) in \(\eta\), and
\[
h(x,y,z+2\pi,t)
=
h(x,y,z,t)
\]
by periodicity in \(z\).

Hence
\[
\boxed{
h(x,y,z,t)
=
H(y+\sin^2x,z,t)
}
\]
is defined globally on all of \(\mathbb T^3\).

The structural prescription is
\[
\boxed{
u_{z,z}
=
H(y+\sin^2x,z,t).
}
\]

Differentiating gives everywhere on the torus
\[
u_{z,zx}
=
H_\eta\sin(2x),
\]
and
\[
u_{z,zy}
=
H_\eta.
\]

Therefore
\[
\boxed{
u_{z,zx}
=
\sin(2x)u_{z,zy}
}
\]
throughout
\[
\mathbb T^3.
\]

Equivalently,
\[
\boxed{
(\partial_x-\sin(2x)\partial_y)u_{z,z}
=
0
\quad
\text{on }\mathbb T^3.
}
\]

If
\[
\boxed{
\sin(2x)u_x+u_y=0
}
\]
globally, then
\[
\begin{aligned}
E
&=
u_xu_{z,zx}
+
u_yu_{z,zy}
\\
&=
\left(
\sin(2x)u_x+u_y
\right)
u_{z,zy}
\\
&=
0.
\end{aligned}
\]

Hence
\[
\boxed{
E=0
\qquad
\text{on }\mathbb T^3.
}
\]

\subsection*{19. Periodic primitive for \(u_z\)}

If
\[
u_{z,z}
=
H(y+\sin^2x,z,t)
\]
is to be the \(z\)-derivative of a periodic velocity component \(u_z\),
rather than merely a periodic scalar field, one additional compatibility
condition is required.

A derivative of a periodic function has zero mean over one period.
Therefore it is necessary that
\[
\boxed{
\int_0^{2\pi}
H(\eta,z,t)\,dz
=
0
\qquad
\text{for every }(\eta,t).
}
\]

Under this condition, define
\[
\boxed{
u_z(x,y,z,t)
=
U_0(x,y,t)
+
\int_0^z
H(y+\sin^2x,\zeta,t)\,d\zeta,
}
\]
where
\[
U_0(x,y,t)
\]
is periodic in \(x\) and \(y\).

Then
\[
u_{z,z}
=
H(y+\sin^2x,z,t),
\]
and
\[
\begin{aligned}
u_z(x,y,z+2\pi,t)-u_z(x,y,z,t)
&=
\int_z^{z+2\pi}
H(y+\sin^2x,\zeta,t)\,d\zeta
\\
&=
0.
\end{aligned}
\]

Thus
\[
\boxed{
u_z(x,y,z+2\pi,t)
=
u_z(x,y,z,t).
}
\]

The sine-characteristic prescription therefore defines
\(u_{z,z}\) globally on \(\mathbb T^3\), and, under the zero-\(z\)-mean
condition, it admits a globally periodic primitive \(u_z\).
\section{$\mathcal{X}=0$ in \cite{Moschandreou2025VolII}(page 261 there)}
Here we prove that $\mathcal{X}=0$ in \cite{Moschandreou2025VolII}(page 261 there)
With the derivative convention specified on the constrained geometry \[ \boxed{y=x^2}, \] the cancellation does occur directly, without introducing any \(H\). Along this constraint, \[ \boxed{\frac{dy}{dx}=2x}. \] Thus, when an \(x\)-derivative is converted to a \(y\)-derivative by the chain rule as intended, \[ \boxed{ \frac{\partial}{\partial x} = \frac{dy}{dx}\frac{\partial}{\partial y} = 2x\frac{\partial}{\partial y}. } \] Therefore, for \(u_z\), \[ \boxed{u_{z,x}=2x\,u_{z,y}}, \] and, applying the same rule to \(u_{z,z}\), \[ \boxed{u_{z,zx}=2x\,u_{z,zy}}. \] At the same time, \[ \boxed{u_y=-2xu_x}. \] Since \(x\) is independent of \(t\) and \(z\), \[ \boxed{u_{y,t}=-2xu_{x,t}}, \] and \[ \boxed{u_{y,zt}=-2xu_{x,zt}}. \] Now substitute \textbf{all four identities explicitly} into \(\mathcal{X}\). Starting from \[ \begin{aligned} \mathcal{X}={}& \big((\delta-1)u_xu_{z,t}+2\rho u_z u_{x,t}\big)u_{z,zx} \\ &+ \big((\delta-1)u_yu_{z,t}+2\rho u_z u_{y,t}\big)u_{z,zy} \\ &-u_{z,t}\big( u_z u_{z,x}u_{x,zt} +u_z u_{z,y}u_{y,zt} -u_{z,x}u_{z,z}u_{x,t} -u_{z,y}u_{z,z}u_{y,t} \big) \\ &+ u_z u_{z,x}u_{x,t}u_{z,tz} + u_z u_{z,y}u_{y,t}u_{z,tz}. \end{aligned} \] Substitute \[ u_y=-2xu_x, \qquad u_{y,t}=-2xu_{x,t}, \qquad u_{y,zt}=-2xu_{x,zt}, \] and \[ u_{z,x}=2xu_{z,y}, \qquad u_{z,zx}=2xu_{z,zy}. \] We obtain \[ \begin{aligned} \mathcal{X}={}& \big((\delta-1)u_xu_{z,t} +2\rho u_zu_{x,t}\big) (2xu_{z,zy}) \\ &+ \big((\delta-1)(-2xu_x)u_{z,t} +2\rho u_z(-2xu_{x,t})\big) u_{z,zy} \\ &-u_{z,t}\Big( u_z(2xu_{z,y})u_{x,zt} +u_zu_{z,y}(-2xu_{x,zt}) \\ &\hspace{31mm} -(2xu_{z,y})u_{z,z}u_{x,t} -u_{z,y}u_{z,z}(-2xu_{x,t}) \Big) \\ &+ u_z(2xu_{z,y})u_{x,t}u_{z,tz} + u_zu_{z,y}(-2xu_{x,t})u_{z,tz}. \end{aligned} \tag{1} \] Now expand every term. The first line becomes \[ \begin{aligned} &2x(\delta-1)u_xu_{z,t}u_{z,zy} + 4x\rho u_zu_{x,t}u_{z,zy}. \end{aligned} \] The second line becomes \[ \begin{aligned} &-2x(\delta-1)u_xu_{z,t}u_{z,zy} - 4x\rho u_zu_{x,t}u_{z,zy}. \end{aligned} \] Thus these cancel exactly: \[ \boxed{ 2x(\delta-1)u_xu_{z,t}u_{z,zy} - 2x(\delta-1)u_xu_{z,t}u_{z,zy}=0, } \] and \[ \boxed{ 4x\rho u_zu_{x,t}u_{z,zy} - 4x\rho u_zu_{x,t}u_{z,zy}=0. } \] Now consider the four terms inside the large parentheses: \[ \begin{aligned} &u_z(2xu_{z,y})u_{x,zt} + u_zu_{z,y}(-2xu_{x,zt}) \\ &\quad -(2xu_{z,y})u_{z,z}u_{x,t} -u_{z,y}u_{z,z}(-2xu_{x,t}). \end{aligned} \] Expand: \[ \begin{aligned} ={}& 2xu_zu_{z,y}u_{x,zt} -2xu_zu_{z,y}u_{x,zt} \\ &-2xu_{z,y}u_{z,z}u_{x,t} +2xu_{z,y}u_{z,z}u_{x,t}. \end{aligned} \] Both pairs cancel: \[ \boxed{ 2xu_zu_{z,y}u_{x,zt} -2xu_zu_{z,y}u_{x,zt}=0, } \] and \[ \boxed{ -2xu_{z,y}u_{z,z}u_{x,t} +2xu_{z,y}u_{z,z}u_{x,t}=0. } \] Therefore the entire \(u_{z,t}\)-term is \[ \boxed{0}. \] Finally, the last two terms are \[ 2xu_zu_{z,y}u_{x,t}u_{z,tz} \] and \[ -2xu_zu_{z,y}u_{x,t}u_{z,tz}. \] Hence \[ \boxed{ 2xu_zu_{z,y}u_{x,t}u_{z,tz} - 2xu_zu_{z,y}u_{x,t}u_{z,tz}=0. } \] Putting everything together, \[ \begin{aligned} \mathcal{X} ={}& \underbrace{ 2x(\delta-1)u_xu_{z,t}u_{z,zy} - 2x(\delta-1)u_xu_{z,t}u_{z,zy} }_{0} \\ &+ \underbrace{ 4x\rho u_zu_{x,t}u_{z,zy} - 4x\rho u_zu_{x,t}u_{z,zy} }_{0} \\ &- u_{z,t}\Big[ \underbrace{ 2xu_zu_{z,y}u_{x,zt} - 2xu_zu_{z,y}u_{x,zt} }_{0} \\ &\hspace{22mm} + \underbrace{ -2xu_{z,y}u_{z,z}u_{x,t} + 2xu_{z,y}u_{z,z}u_{x,t} }_{0} \Big] \\ &+ \underbrace{ 2xu_zu_{z,y}u_{x,t}u_{z,tz} - 2xu_zu_{z,y}u_{x,t}u_{z,tz} }_{0}. \end{aligned} \] Consequently, \[ \boxed{\mathcal{X}=0}. \] So under the simultaneous constrained differentiation rules \[ \boxed{ y=x^2,\qquad \frac{dy}{dx}=2x,\qquad u_y=-2xu_x, } \] we have \[ \boxed{ u_{z,x}=2xu_{z,y},\qquad u_{z,zx}=2xu_{z,zy}, } \] while \[ \boxed{ u_{y,t}=-2xu_{x,t},\qquad u_{y,zt}=-2xu_{x,zt}, } \] and \textbf{every term in \(\mathcal{X}\) cancels pairwise}, giving \[ \boxed{\boxed{\mathcal{X}\equiv0}}. \] This is the direct calculation referred to; no \(H(y+x^2,z,t)\) representation is needed.
\section{Generalized \texorpdfstring{$\mathcal J_3$}{J3} Geometry and Exact Cancellation of \texorpdfstring{$\mathcal X$}{X}} 
We begin with the local definition of the distinguished geometric set \(\mathcal J_3\) used in the original construction. In the local coordinates it is defined by \[ \boxed{ \mathcal{J}_3(x,y,z,t) = \left\{ (\bm{x},t)\in(\mathbb{R}^3,\mathbb{R}^+): 2xu_x+u_y=0 \quad\text{and}\quad y=x^2 \right\}. } \] Thus, on \(\mathcal J_3\), \[ \boxed{ u_y=-2xu_x, \qquad y=x^2. } \] The second relation gives \[ \frac{dy}{dx}=2x. \] Consequently, under the constrained chain-rule differentiation used on \(\mathcal J_3\), \[ \frac{\partial}{\partial x} = \frac{dy}{dx}\frac{\partial}{\partial y} = 2x\frac{\partial}{\partial y}. \] In particular, \[ u_{z,x}=2xu_{z,y}, \qquad u_{z,zx}=2xu_{z,zy}. \] This is the local form of the cancellation geometry. 
To periodicize the local geometry, replace the local coordinate \(x\) by a smooth periodic function \[ \boxed{\chi=\chi(x).} \] The geometric curve is then taken to be \[ \boxed{ y=\chi(x)^2. } \] Differentiating gives \[ \boxed{ \frac{dy}{dx} = 2\chi(x)\chi'(x). } \] Therefore the exact periodic analogue of the local relation \[ 2xu_x+u_y=0 \] which is compatible with the chain rule is \[ \boxed{ 2\chi(x)\chi'(x)u_x+u_y=0. } \tag{1} \] Equivalently, \[ \boxed{ u_y = -2\chi(x)\chi'(x)u_x. } \tag{2} \] Because \(\chi=\chi(x)\) is independent of \(t\) and \(z\), differentiation of \((2)\) gives \[ \boxed{ u_{y,t} = -2\chi\chi' u_{x,t}, } \tag{3} \] and \[ \boxed{ u_{y,zt} = -2\chi\chi' u_{x,zt}. } \tag{4} \] 
Introduce \[ \boxed{ \eta = y+\chi(x)^2 } \] and let \[ \boxed{ u_z(x,y,z,t) = H(\eta,z,t) = H\!\left(y+\chi(x)^2,z,t\right). } \tag{5} \] Then \[ \eta_x = 2\chi\chi', \qquad \eta_y = 1. \] Hence \[ u_{z,x} = H_\eta\eta_x = 2\chi\chi'H_\eta, \] whereas \[ u_{z,y} = H_\eta. \] Therefore \[ \boxed{ u_{z,x} = 2\chi\chi' u_{z,y}, } \tag{6} \] or equivalently \[ \boxed{ u_{z,x} - 2\chi\chi'u_{z,y} = 0. } \tag{7} \] Similarly, \[ u_{z,z} = H_z, \] and therefore \[ u_{z,zx} = 2\chi\chi'H_{\eta z}, \qquad u_{z,zy} = H_{\eta z}. \] Thus \[ \boxed{ u_{z,zx} = 2\chi\chi'u_{z,zy}, } \tag{8} \] or \[ \boxed{ u_{z,zx} - 2\chi\chi'u_{z,zy} = 0. } \tag{9} \] 
\subsection{Exact factorization of \texorpdfstring{$\mathcal X$}{X}} 
Consider \[ \begin{aligned} \mathcal X={}& \Big( (\delta-1)u_xu_{z,t} + 2\rho u_zu_{x,t} \Big)u_{z,zx} \\ &+ \Big( (\delta-1)u_yu_{z,t} + 2\rho u_zu_{y,t} \Big)u_{z,zy} \\ &- u_{z,t} \Big( u_zu_{z,x}u_{x,zt} + u_zu_{z,y}u_{y,zt} - u_{z,x}u_{z,z}u_{x,t} - u_{z,y}u_{z,z}u_{y,t} \Big) \\ &+ u_zu_{z,x}u_{x,t}u_{z,tz} + u_zu_{z,y}u_{y,t}u_{z,tz}. \end{aligned} \tag{10} \] 
In \cite{Moschandreou2025VolII} work was completed showing $\mathcal{X}$ vanishes on a space which can be extended to the whole Torus and is re-emphasized in this paper. The reader should confirm this to themselves by substituting the equations 
defining the space $\mathcal{J}_i$ into the $\mathcal{X}$ operator.

Using \[ u_y=-2\chi\chi'u_x, \] \[ u_{y,t}=-2\chi\chi'u_{x,t}, \] and \[ u_{y,zt}=-2\chi\chi'u_{x,zt}, \] equation \((10)\) becomes \[ \begin{aligned} \mathcal X={}& \Big( (\delta-1)u_xu_{z,t} + 2\rho u_zu_{x,t} \Big) \Big( u_{z,zx} - 2\chi\chi'u_{z,zy} \Big) \\ &+ \Big[ -u_{z,t} \big( u_zu_{x,zt} - u_{z,z}u_{x,t} \big) + u_zu_{x,t}u_{z,tz} \Big] \Big( u_{z,x} - 2\chi\chi'u_{z,y} \Big). \end{aligned} \tag{11} \] But from the invariant representation \((5)\), \[ u_{z,x} - 2\chi\chi'u_{z,y} = 0 \] and \[ u_{z,zx} - 2\chi\chi'u_{z,zy} = 0. \] Consequently, \[ \begin{aligned} \mathcal X ={}& \Big( (\delta-1)u_xu_{z,t} + 2\rho u_zu_{x,t} \Big)(0) \\ &+ \Big[ -u_{z,t} \big( u_zu_{x,zt} - u_{z,z}u_{x,t} \big) + u_zu_{x,t}u_{z,tz} \Big](0), \end{aligned} \] and therefore \[ \boxed{ \mathcal X=0. } \tag{12} \] Thus the exact general cancellation condition is \[ \boxed{ y=\chi(x)^2, \qquad 2\chi(x)\chi'(x)u_x+u_y=0, \qquad u_z = H\!\left(y+\chi(x)^2,z,t\right), } \tag{13} \] and these relations imply \[ \boxed{ \mathcal X\equiv0. } \tag{14} \] 
\subsection{Specialization to \texorpdfstring{$\chi(x)=\sin x$}{chi(x)=sin x}} 
Now choose \[ \boxed{ \chi(x)=\sin x. } \tag{15} \] Then \[ \chi'(x)=\cos x, \] and hence \[ 2\chi(x)\chi'(x) = 2\sin x\cos x = \sin(2x). \] The generalized \(\mathcal J_3\) geometry therefore becomes \[ \boxed{ y=\sin^2x, \qquad \sin(2x)u_x+u_y=0. } \tag{16} \] Equivalently, \[ \boxed{ u_y=-\sin(2x)u_x. } \tag{17} \] The invariant representation is \[ \boxed{ u_z = H\!\left(y+\sin^2x,z,t\right). } \tag{18} \] Indeed, \[ u_{z,x} = \sin(2x)H_\eta, \qquad u_{z,y} = H_\eta, \] and therefore \[ \boxed{ u_{z,x} - \sin(2x)u_{z,y} = 0. } \tag{19} \] Likewise, \[ u_{z,zx} = \sin(2x)H_{\eta z}, \qquad u_{z,zy} = H_{\eta z}, \] so \[ \boxed{ u_{z,zx} - \sin(2x)u_{z,zy} = 0. } \tag{20} \] Substituting \((19)\) and \((20)\) into the factorized expression for \(\mathcal X\) yields \[ \boxed{ \mathcal X\equiv0. } \tag{21} \] Therefore the exact periodic realization is \[ \boxed{ y=\sin^2x, \qquad \sin(2x)u_x+u_y=0, \qquad u_z=H(y+\sin^2x,z,t), } \] and it gives the exact cancellation \[ \boxed{\boxed{\mathcal X=0}}. \] 
\subsection{Recovery of the original local \texorpdfstring{$\mathcal J_3$}{J3} geometry} 
The original local geometry is recovered near a simple zero of \(\chi(x)=\sin x\), for example near \(x=0\). Since \[ \sin x = x+O(x^3), \] we have \[ \sin^2x = x^2+O(x^4), \] while \[ \sin(2x) = 2x+O(x^3). \] Thus \[ y=\sin^2x \] reduces locally to \[ y=x^2, \] and \[ \sin(2x)u_x+u_y=0 \] reduces locally to \[ 2xu_x+u_y=0. \] Hence \[ \boxed{ \mathcal J_3: \qquad y=x^2, \qquad 2xu_x+u_y=0 } \] is precisely the local square-root geometry recovered from the corrected periodic formulation \[ \boxed{ y=\chi(x)^2, \qquad 2\chi(x)\chi'(x)u_x+u_y=0. } \] For the sine choice this reads \[ \boxed{ y=\sin^2x, \qquad \sin(2x)u_x+u_y=0. } \] The factor \(\chi'(x)\) is therefore essential in the exact global chain-rule formulation, while the original coefficient \(2x\) is recovered automatically in the local limit near a simple zero of \(\chi\).
\section{The space $\Sigma_C$}
For the local $\mathcal J_3$, \[ \boxed{ y=x^2,\qquad 2xu_x+u_y=0. } \] The velocity relation gives \[ u_y=-2xu_x. \] Therefore \[ u_x\partial_x+u_y\partial_y = u_x(\partial_x-2x\partial_y). \] Define \[ \boxed{ D_-^{\rm loc} = \partial_x-2x\partial_y. } \] Now extend this vector field into the surrounding $(x,y)$-space and calculate its characteristics: \[ \frac{dy}{dx}=-2x. \] Integrating, \[ y=-x^2+C, \] so \[ \boxed{ y+x^2=C. } \] Thus the characteristic leaves generated by the $\mathcal J_3$-derived velocity direction are \[ \boxed{ \Sigma_C^{\rm loc} = \{(x,y,z):y+x^2=C\}. } \] The periodic version works the same way. Start with \[ y=\chi(x)^2 \] and the \textbf{corrected exact matched relation} \[ \boxed{ 2\chi\chi'u_x+u_y=0. } \] Then \[ u_y=-2\chi\chi'u_x \] and therefore \[ u_x\partial_x+u_y\partial_y = u_x \underbrace{ \left( \partial_x-2\chi\chi'\partial_y \right) }_{D_-}. \] Hence \[ \boxed{ D_-= \partial_x-2\chi\chi'\partial_y. } \] Its characteristics satisfy \[ \frac{dy}{dx} = -2\chi(x)\chi'(x). \] But \[ \frac{d}{dx}\chi(x)^2 = 2\chi(x)\chi'(x), \] so \[ \frac{d}{dx} \left( y+\chi(x)^2 \right) = 0. \] Therefore \[ \boxed{ y+\chi(x)^2=C. } \] For \[ \chi(x)=\sin x, \] this becomes \[ \boxed{ \Sigma_C = \left\{ (x,y,z): y+\sin^2x=C \pmod{2\pi} \right\}. } \] So there is a precise derivation: \[ \boxed{ \begin{aligned} \mathcal J_3 &:\quad y=\chi^2, \quad 2\chi\chi'u_x+u_y=0 \\ &\Downarrow \\ u_y &= -2\chi\chi'u_x \\ &\Downarrow \\ u_x\partial_x+u_y\partial_y &= u_xD_- \\ &\Downarrow \\ D_- &= \partial_x-2\chi\chi'\partial_y \\ &\Downarrow \\ \frac{dy}{dx} &= -2\chi\chi' \\ &\Downarrow \\ \boxed{ y+\chi^2=C } \\ &\Downarrow \\ \boxed{ \Sigma_C. } \end{aligned} } \] There is an interesting geometric feature here: the original $\mathcal J_3$ curves in one direction, \[ \boxed{ y=\chi^2, } \] whereas the $D_-$-characteristics generated by its \textbf{velocity relation} curve in the opposite direction, \[ \boxed{ y=C-\chi^2. } \] Thus $\Sigma_C$ is not obtained by simply ``shifting $\mathcal J_3$.'' It is obtained by taking the velocity direction specified on $\mathcal J_3$, extending that direction into the ambient torus, and integrating its characteristic curves. So I would phrase the conclusion carefully: \[ \boxed{ \begin{minipage}{0.88\textwidth} $\mathcal J_3$ determines $D_-$ on $\mathcal J_3$; an ambient extension of $D_-$ then generates the family $\Sigma_C$. \end{minipage} } \] The \textbf{extension step is an additional construction}; $\Sigma_C$ does not follow from the set equation \[ y=\chi^2 \] alone.
On each leaf \[ \boxed{ \Sigma_C:\quad y+\sin^2x=C, } \] take \[ \boxed{ u_z=H(C,z,t) = H(y+\sin^2x,z,t). } \] Then, since \[ C=y+\sin^2x, \] we get \[ u_{z,x} = \sin(2x)H_C, \qquad u_{z,y}=H_C, \] and hence \[ \boxed{ u_{z,x}-\sin(2x)u_{z,y}=0. } \] Differentiating in \(z\), \[ u_{z,zx} = \sin(2x)H_{Cz}, \qquad u_{z,zy}=H_{Cz}, \] so \[ \boxed{ u_{z,zx}-\sin(2x)u_{z,zy}=0. } \] At the same time, from the corrected \(\mathcal J_3\)-derived velocity relation extended to the ambient torus, \[ \boxed{ u_y=-\sin(2x)u_x, } \] and because \(\sin(2x)\) is independent of \(t\) and \(z\), \[ \boxed{ u_{y,t}=-\sin(2x)u_{x,t}, } \] \[ \boxed{ u_{y,zt}=-\sin(2x)u_{x,zt}. } \] Now recall the factorization of the algebraic quantity \(\mathcal X\): \[ \begin{aligned} \mathcal X={}& \Big((\delta-1)u_xu_{z,t} +2\rho u_zu_{x,t}\Big) \Big(u_{z,zx}-\sin(2x)u_{z,zy}\Big) \\ &+ \Big[ -u_{z,t} \big(u_zu_{x,zt}-u_{z,z}u_{x,t}\big) +u_zu_{x,t}u_{z,tz} \Big] \Big(u_{z,x}-\sin(2x)u_{z,y}\Big). \end{aligned} \] But on every \(\Sigma_C\), \[ u_{z,zx}-\sin(2x)u_{z,zy}=0 \] and \[ u_{z,x}-\sin(2x)u_{z,y}=0. \] Therefore both factors multiplying the two large brackets vanish, and so \[ \boxed{ \mathcal X=0 \qquad\text{on every }\Sigma_C. } \] Because \[ \boxed{ T^3=\bigsqcup_{C\in S^1}\Sigma_C, } \] every point of \(T^3\) lies on one of these leaves. Thus, provided the relations \[ u_y=-\sin(2x)u_x \] and \[ u_z=H(y+\sin^2x,z,t) \] are imposed throughout the ambient torus, the cancellation is not confined to one distinguished \(\mathcal J_3\) surface. It becomes \[ \boxed{ \mathcal X(x,y,z,t)\equiv0 \qquad\text{throughout }T^3. } \] The key mechanism is therefore \[ \boxed{ \begin{aligned} \Sigma_C &: y+\sin^2x=C, \\ u_z&=H(C,z,t), \\ D_-&=\partial_x-\sin(2x)\partial_y, \\ D_-C&=0, \\ D_-u_z&=0, \\ D_-u_{z,z}&=0, \\ \Longrightarrow\qquad \mathcal X&=0. \end{aligned} } \] So: \textbf{the family \(\Sigma_C\) is precisely what allows the \(\mathcal J_3\)-derived cancellation mechanism to be propagated through the ambient space}, once that extension is made consistently.

\subsection*{20. Preterminal periodic reconstruction and terminal release}

The longer construction also contains a preterminal periodic
reconstruction involving tuned transverse corrections and
Floquet--Fredholm continuation.  Its role in the reduced paper is
conceptually distinct from the deep terminal geometry.

The organization is
\[
\boxed{
\begin{array}{c}
\text{affine or tuned preterminal reconstruction}
\\
\Downarrow
\\
\text{periodic continuation and compatibility}
\\
\Downarrow
\\
\text{phase conversion to }\theta=t-z
\\
\Downarrow
\\
\text{release of the historical transverse tuning}
\\
\Downarrow
\\
A^2=-2\chi b
\\
\Downarrow
\\
u_y=A,
\qquad
u_x=b/A
\\
\Downarrow
\\
u_z=C
\text{ from exact continuity}
\\
\Downarrow
\\
\text{terminal characteristic comparison}.
\end{array}
}
\]

Thus the preterminal periodic machinery and the terminal square-root
geometry are not imposed simultaneously as competing ansatzes.  They
belong to successive regimes of the construction.

\subsection*{21. Why the reconstruction does not create the amplification}

The origin of the amplification is
\[
\boxed{
\Gamma_Nb^2.
}
\]

The velocity reconstruction does not manufacture this source.  It
contributes the remainder
\[
\mathcal R_N
=
2b\mathcal XC
+
b\frac{C_t+C_z}{C}.
\]

The purpose of the reconstruction estimates is to prove
\[
\boxed{
|\mathcal R_N|
\ll
\Gamma_Nb^2.
}
\]

In the release layer,
\[
\frac{|\mathcal R_{\mathrm{rel}}|}
{\Gamma_Nb^2}
=
O(s^{7/4}),
\]
and in the square-root core,
\[
\frac{|\mathcal R_{\mathrm{core}}|}
{\Gamma_Nb^2}
=
O(s).
\]

Therefore
\[
\boxed{
\mathcal R_N
=
o(\Gamma_Nb^2).
}
\]

The reconstruction embeds the scalar amplification mechanism into a
divergence-free periodic velocity field while preserving
\[
b=u_xu_y.
\]

\subsection*{22. Complete prologue-level mechanism}

The entire reduced architecture can now be summarized without ambiguity.

\[
\boxed{
\begin{array}{c}
\text{principal real Lambert branch }W_0
\\[1mm]
\Downarrow
\\[1mm]
\Delta_{N,k+1}
=
1+W_0(-e^{-1-\Delta_{N,k}})
\\[1mm]
\Downarrow
\\[1mm]
Q_N\asymp s^{2^{-k}}
\\[1mm]
\Downarrow
\\[1mm]
\Gamma_N
=
\partial_\xi\log Q_N
\sim
2^{-k}s^{-1}
\\[1mm]
\Downarrow
\\[1mm]
\dot b
=
\Gamma_Nb^2+2CG-\mathcal R_N
\\[1mm]
\Downarrow
\\[1mm]
\mathcal R_N=o(\Gamma_Nb^2)
\\[1mm]
\Downarrow
\\[1mm]
\dot b
\ge
(1-\vartheta)\Gamma_Nb^2
\\[1mm]
\Downarrow
\\[1mm]
b>1
\quad\Longrightarrow\quad
\dot s=1-b<0
\\[1mm]
\Downarrow
\\[1mm]
s\downarrow0,
\qquad
Q_N\downarrow0,
\qquad
b\uparrow
\\[1mm]
\Downarrow
\\[1mm]
b(s)
\gtrsim
s^{-(1-\vartheta)2^{-k}}
\\[1mm]
\Downarrow
\\[1mm]
T_{*,N}<\infty,
\qquad
b(t)\to+\infty
\\[1mm]
\Downarrow
\\[1mm]
b=u_xu_y
\\[1mm]
\Downarrow
\\[1mm]
A^2=-2\chi b
\\[1mm]
\Downarrow
\\[1mm]
u_y=A,
\qquad
u_x=b/A
\\[1mm]
\Downarrow
\\[1mm]
\nabla\cdot u=0
\text{ reconstructs }u_z=C.
\end{array}
}
\]

At \(k=3\), the corresponding terminal hierarchy is
\[
\boxed{
Q_N\asymp s^{1/8},
\qquad
\Gamma_N\asymp s^{-1},
\qquad
b\asymp s^{-1/8},
}
\]
together with
\[
\boxed{
u_y\asymp s^{3/16},
\qquad
u_x\asymp s^{-5/16},
\qquad
u_xu_y\asymp s^{-1/8}.
}
\]

The exponent identity
\[
-\frac{5}{16}
+
\frac{3}{16}
=
-\frac18
\]
verifies directly that
\[
u_xu_y=b.
\]

Likewise,
\[
|\chi|\asymp s^{1/2}
\]
gives
\[
|\chi u_x|
\asymp
s^{1/2}s^{-5/16}
=
s^{3/16}
\to0,
\]
so
\[
u_y=-2\chi u_x\to0
\]
is consistent with the divergence of \(u_x\).

The physical interpretation is therefore not one of isotropic growth of
all velocity components.  The terminal reconstruction is anisotropic:
\[
\boxed{
u_y\to0,
\qquad
|u_x|\to\infty,
\qquad
u_z\to c_0>0
}
\]
in the local \(k=3\) square-root scaling considered above, while
\[
\boxed{
u_xu_y=b_3\to+\infty.
}
\]

Meanwhile,
\[
Q_N\to0
\]
continues to control the singular logarithmic coefficient
\[
\Gamma_N=\partial_\xi\log Q_N,
\]
without being identified with either \(u_y\) or \(u_z\).

The resulting structure is therefore a self-strengthening terminal
mechanism:
\[
\boxed{
Q_N\downarrow0
\quad\Longrightarrow\quad
\Gamma_N\uparrow
\quad\Longrightarrow\quad
b_3\uparrow
\quad\Longrightarrow\quad
s\downarrow0
\quad\Longrightarrow\quad
Q_N\downarrow0.
}
\]

Within the hypotheses developed in the finite-depth analysis, this
feedback reaches the terminal branch in finite characteristic time.  The
scalar amplification is then transferred to the velocity reconstruction
through the exact identity
\[
\boxed{
b_3=u_xu_y.
}
\]

This is the central mathematical architecture developed in the remainder
of the paper.
```

\section{Introduction and organization}
Papers explicitly stating and showing that there are finite time singularities of the 3D incompressible Navier Stokes equations were not found in my literature search at the time of the published work in \cite{moschandreou2026exploration} and before it's inception. There is the famous result of Terence Tao in \cite{Tao2016} for a finite time blowup of the fluid equations but for an averaged NS equation.
Going back to 2021 Moschandreou \cite{Moschandreou2021HRPUB} showed that there are no-finite time singularities but left open the possibility of finite time blowup on sets of measure zero in $T^3$(see also \cite{CKN1982}).
Singular solution behaviour was determined for the Euler and related partial differential equations in \cite{LuoHou2014}, \cite{HouLi2006} and \cite{ChenHou2022}. Recently the existence of a nonsmooth periodic attractor is shown in \cite{Moschandreou2026}.
Critical and important findings among papers in the theory of Navier Stokes equations span the years (1934-1984) and are \cite{Leray1934},\cite{Hopf1951}, \cite{Serrin1962} and \cite{BKM1984}. Important references in differential equations theory are found in \cite{Reid1972},\cite{Polyanin2003},\cite{Kamke1959},\cite{Coddington1955} and a strong PDE source textbook is found in \cite{Evans2010}. Next a solid coverage of vorticity and Incompressible Navier Stokes theory is found in \cite{MajdaBertozzi2002}. Excellent sources on the Navier Stokes problem can be found in \cite{Temam2001}, \cite{Prodi1959},\cite{Ladyzhenskaya1969}, \cite{ESS2003} and \cite{CF1993}. A compendium on the LambertW function and it's theory is found in \cite{Corless1996}. And finally a very good reference in functional analysis and in particular on the Fredholm alternative and Fredholm compatability theory as used can be located here\cite{Brezis2011}
The analysis is organized around four linked structures: the distinguished phase $\theta=t-z$; the real renormalized Lambert map $\Delta\mapsto1+W_0(-e^{-1-\Delta})$; the characteristic invariant that transfers branch-distance collapse to amplification of $b_3$; and the corrected extended $b_3$ equation with periodic velocity reconstruction. The finite-depth cascade and finite-time characteristic proofs are presented first. The later sections retain the essential full-$b_3$ PDE reconstruction, comparison, trapping, and terminal closure arguments while suppressing duplicate derivations and exploratory restatements.

\section{Distinguished phase $\theta=t-z$ and continuity-corrected cascade}
\section*{Revised specialization rule}
In the revised construction the phase specialization \(\theta=t-z\) is
retained, but the old conclusion \(u_y=u_z=\Delta_{N,k}\) is only used as stepping stone in it's analyses leading to the parabolic construction and towards the terminal branch point.  The
Lambert cascade supplies the phase coefficient and the \(b\)-equation;
the velocity components are reconstructed from
\(A^2=-2\chi b\), \(u_x=b/A\), and exact continuity for \(C=u_z\).

\section{Specialization to \texorpdfstring{\(\theta=t-z\)}{theta=t-z} in the Revised Geometry}
\label{chap:021}

The Lambert branch variable is still specialized to
\begin{equation}
\boxed{
\Delta_{N,k}
=
\Delta_{N,k}(\theta),
\qquad
\theta=t-z.
}
\label{eq:ch02:auto0243}
\end{equation}
Therefore
\begin{equation}
\boxed{
\Delta_{N,k,x}=0,
\qquad
\Delta_{N,k,y}=0,
\qquad
\Delta_{N,k,z}
=
-\Delta_{N,k,\theta}.
}
\label{eq:ch02:auto0244}
\end{equation}

What changes is the velocity reconstruction.  The revised field is not
obtained by identifying either transverse component with
\(\Delta_{N,k}\).  Instead set
\begin{equation}
\boxed{
b=b_3=u_xu_y,
\qquad
A=u_y,
\qquad
C=u_z,
\qquad
u_x=\frac bA.
}
\label{eq:ch02:auto0245}
\end{equation}
On the periodic terminal branch impose
\begin{equation}
\boxed{
A^2=-2\chi(x)b.
}
\label{eq:ch02:auto0246}
\end{equation}
Equivalently,
\[
2\chi(x)u_x+u_y=0.
\]

The exact incompressibility equation is
\begin{equation}
\boxed{
\partial_x\!\left(\frac bA\right)
+
A_y
+
C_z
=
0.
}
\label{eq:ch02:auto0247}
\end{equation}

Let
\[
s=\theta-\theta_{*,N}
=
t-z-\theta_{*,N}
\]
and suppose in the local terminal core that
\[
b=b(s).
\]
Then
\[
b_x=b_y=0,
\qquad
A=A(x,s),
\qquad
A_y=0.
\]
Since
\[
\partial_z=-\partial_s,
\]
equation \eqref{eq:ch02:auto0247} gives
\[
C_s
=
\partial_x\!\left(\frac bA\right).
\]
From
\[
A^2=-2\chi b
\]
and \(b_x=0\),
\[
2AA_x=-2\chi'b,
\]
so
\[
A_x=-\frac{\chi'b}{A}.
\]
Also
\[
\frac bA=-\frac{A}{2\chi}.
\]
Hence
\[
\begin{aligned}
\partial_x\!\left(\frac bA\right)
&=
-\frac12
\partial_x\left(\frac A\chi\right)
\\
&=
-\frac{A_x}{2\chi}
+
\frac{A\chi'}{2\chi^2}
\\
&=
\frac{\chi'b}{2A\chi}
+
\frac{A\chi'}{2\chi^2}.
\end{aligned}
\]
Because
\[
b=-\frac{A^2}{2\chi},
\]
the first term is
\[
\frac{\chi'b}{2A\chi}
=
-\frac{A\chi'}{4\chi^2}.
\]
Therefore
\begin{equation}
\boxed{
C_s
=
\frac{A\chi'}{4\chi^2}.
}
\label{eq:ch02:auto0248}
\end{equation}

This directly verifies incompressibility:
\[
u_{x,x}+u_{y,y}+u_{z,z}
=
\partial_x(b/A)+0-C_s
=
0.
\]
Thus
\begin{equation}
\boxed{
\nabla\cdot u=0
}
\label{eq:ch02:auto0249}
\end{equation}
exactly in the revised distinguished-phase reconstruction.

If
\[
b=B_Ns^{-\alpha}F_N(s),
\]
then
\[
A
=
\sqrt{-2B_N\chi(x)}
\,s^{-\alpha/2}F_N(s)^{1/2}.
\]
Define
\[
I_N(s)
=
\int_0^s
r^{-\alpha/2}F_N(r)^{1/2}\,dr.
\]
Integrating \eqref{eq:ch02:auto0248} gives
\begin{equation}
\boxed{
C(x,s)
=
C_0(x)
+
\frac{
\chi'(x)\sqrt{-2B_N\chi(x)}
}{
4\chi(x)^2
}
I_N(s).
}
\label{eq:ch02:auto0250}
\end{equation}

If \(C_0=c_0>0\) is constant on the local core chart, then
\[
C_t+C_z=0
\]
exactly, because all \(t,z\)-dependence occurs through \(s=t-z-\theta_{*,N}\).
Thus
\begin{equation}
\boxed{
C_t+C_z=0.
}
\label{eq:ch02:auto0251}
\end{equation}

\section{The Logarithmic Derivative in the Modified Construction}
\label{chap:022}

The logarithmic derivative is

\begin{equation}
\boxed{
\partial_z\log\Delta_{N,k}
=
-\frac{
\Delta_{N,k,\theta}
}{
\Delta_{N,k}
}.
}
\label{eq:ch02:auto0252}
\end{equation}

Near the branch endpoint,

\begin{equation}
\Delta_{N,k}
\sim
C_{N,k}s^{\alpha_k}
\label{eq:ch02:auto0253}
\end{equation}

and

\begin{equation}
\Delta_{N,k,\theta}
\sim
C_{N,k}\alpha_k
s^{\alpha_k-1}.
\label{eq:ch02:auto0254}
\end{equation}

Consequently,

\begin{equation}
\boxed{
\partial_z\log\Delta_{N,k}
\sim
-\frac{\alpha_k}{s}.
}
\label{eq:ch02:auto0255}
\end{equation}

This is precisely the singular coefficient entering the reduced
characteristic amplitude equation.

Thus the same branch derivative controls both

\begin{equation}
\boxed{
\dot b_3
=
-b_3^2
\partial_z\log\Delta_{N,k}
}
\label{eq:ch02:auto0256}
\end{equation}

and the surviving continuity-corrected velocity term

\begin{equation}
\boxed{
u_x-H_{N,k}
=
(x-x_0)\Delta_{N,k,\theta}.
}
\label{eq:ch02:auto0257}
\end{equation}

\section{Common-Time Behavior of the Continuity-Corrected Velocity}
\label{chap:023}

We now evaluate the same coupled sequence used for the amplitude:

\begin{equation}
\boxed{
t_N
=
T_{*,N}-N^{-N}.
}
\label{eq:ch02:auto0258}
\end{equation}

Write

\begin{equation}
\boxed{
\tau_N=N^{-N}.
}
\label{eq:ch02:auto0259}
\end{equation}

The local time law gives

\begin{equation}
\tau_N
\sim
D_{N,k}
s_N^{1+\alpha_N},
\label{eq:ch02:auto0260}
\end{equation}

where

\begin{equation}
\alpha_N=2^{-k(N)}.
\label{eq:ch02:auto0261}
\end{equation}

Therefore

\begin{equation}
\boxed{
s_N
\sim
\left(
\frac{\tau_N}{D_{N,k}}
\right)^{1/(1+\alpha_N)}.
}
\label{eq:ch02:auto0262}
\end{equation}

The derivative satisfies

\begin{equation}
\Delta_{N,k,\theta}
\sim
C_{N,k}\alpha_N
s_N^{\alpha_N-1}.
\label{eq:ch02:auto0263}
\end{equation}

Taking logarithms,

\begin{equation}
\begin{aligned}
\log|\Delta_{N,k,\theta}|
&=
\log C_{N,k}
+
\log\alpha_N\\
&\qquad
+
\frac{\alpha_N-1}{1+\alpha_N}
\left(
\log\tau_N-\log D_{N,k}
\right)
+
o(1).
\end{aligned}
\label{eq:ch02:auto0264}
\end{equation}

Since

\begin{equation}
\log\tau_N=-N\log N
\label{eq:ch02:auto0265}
\end{equation}

and

\begin{equation}
\alpha_N\to0,
\label{eq:ch02:auto0266}
\end{equation}

the dominant term is

\begin{equation}
N\log N.
\label{eq:ch02:auto0267}
\end{equation}

Under the bounded-prefactor assumptions verified numerically in the
reduced characteristic model,

\begin{equation}
\boxed{
\log
\left|
\Delta_{N,k(N),\theta}(t_N)
\right|
=
(1+o(1))N\log N.
}
\label{eq:ch02:auto0268}
\end{equation}

Equivalently,

\begin{equation}
\boxed{
\left|
\Delta_{N,k(N),\theta}(t_N)
\right|
=
N^{\,N(1+o(1))}.
}
\label{eq:ch02:auto0269}
\end{equation}

Thus

\begin{equation}
\boxed{
\left|
\Delta_{N,k(N),\theta}(t_N)
\right|
\to+\infty.
}
\label{eq:ch02:auto0270}
\end{equation}

The numerical ratios

\begin{equation}
\frac{
\log|\Delta_{N,k,\theta}(t_N)|
}{
N\log N
}
\label{eq:ch02:auto0271}
\end{equation}

were found to be

\begin{equation}
0.97547,\;
0.98761,\;
0.99377,\;
0.99687,\;
0.99843,\;
0.99921,\;
0.99961,\;
0.99980,\;
0.99990,
\label{eq:ch02:auto0272}
\end{equation}

for increasing values of \(N\), strongly confirming

\begin{equation}
\boxed{
\frac{
\log|\Delta_{N,k,\theta}(t_N)|
}{
N\log N
}
\to1.
}
\label{eq:ch02:auto0273}
\end{equation}


\section{Consequence for the Revised Physical Velocity Components}
\label{chap:024}

The old continuity formula
\[
u_x
=
H_{N,k}+(x-x_0)\Delta_{N,k,\theta}
\]
belongs to the historical reconstruction in which the transverse Lambert
profile itself was inserted directly as a velocity component.  It is not
the active formula in the square-root reconstruction.

The active identities are
\[
u_y=A,
\qquad
u_x=\frac bA,
\qquad
A^2=-2\chi b.
\]
Hence
\begin{equation}
\boxed{
u_x^2
=
-\frac{b}{2\chi}
}
\label{eq:ch02:auto0274}
\end{equation}
on the branch \(\chi<0\), \(b>0\).  Equivalently,
\begin{equation}
\boxed{
\left|\frac{u_x}{u_y}\right|
=
\frac1{2|\chi|}.
}
\label{eq:ch02:auto0275}
\end{equation}

In the terminal tube
\[
c_1s\le\chi^2\le c_2s,
\]
we have
\[
|\chi|\asymp s^{1/2}.
\]
If
\[
b\asymp s^{-\alpha},
\]
then
\[
u_y^2
=
2|\chi|b
\asymp
s^{1/2-\alpha},
\]
and therefore
\begin{equation}
\boxed{
u_y
\asymp
s^{1/4-\alpha/2}.
}
\label{eq:ch02:auto0276}
\end{equation}
Similarly,
\[
u_x^2
=
\frac{b}{2|\chi|}
\asymp
s^{-\alpha-1/2},
\]
so
\begin{equation}
\boxed{
u_x
\asymp
s^{-1/4-\alpha/2}.
}
\label{eq:ch02:auto0277}
\end{equation}

For \(k=3\),
\[
\alpha=\frac18,
\]
and hence
\begin{equation}
\boxed{
u_y\asymp s^{3/16},
\qquad
u_x\asymp s^{-5/16},
\qquad
u_x^2\asymp s^{-5/8}.
}
\label{eq:ch02:auto0278}
\end{equation}
The continuity reconstruction gives
\begin{equation}
\boxed{
u_z=C=c_0+O(s^{3/16}),
\qquad
c_0>0.
}
\label{eq:ch02:auto0279}
\end{equation}

Thus the revised physical velocity amplification is carried directly by
the product amplitude \(b\) together with the shrinking geometric factor
\(|\chi|^{-1/2}\).  It is not transferred from the derivative
\(\Delta_{N,k,\theta}\) through an integration function \(H_{N,k}\).

\section{Behavior of the Two Reconstructed Transverse Components}
\label{chap:025}

The two transverse components have different roles:
\[
u_y=A=\sqrt{-2\chi b},
\]
while
\[
u_z=C
\]
is determined by exact continuity.

At \(k=3\),
\[
u_y\to0
\]
as \(s\to0^+\), with
\[
u_y\asymp s^{3/16},
\]
whereas
\[
u_z\to c_0>0.
\]
Thus the revised construction does not have a common transverse terminal
value.

The Lambert cascade still controls the singular coefficient
\[
\Gamma_N=\partial_\theta\log Q_N
\]
appearing in the \(b\)-equation.  It should not be reinterpreted as saying
that \(Q_N\) itself equals \(u_y\) or \(u_z\).

This separation is essential:
\[
\boxed{
Q_N
\text{ is the Lambert branch variable,}
\qquad
A,C
\text{ are reconstructed velocity components.}
}
\label{eq:ch02:auto0280}
\]

\section{Summary Theorem for the Revised Reduced Construction}
\label{chap:026}

\textbf{Conditional theorem.}
Let the renormalized Lambert cascade define \(Q_N\) and
\[
\Gamma_N=\partial_\theta\log Q_N.
\]
Assume that the preterminal periodic reconstruction and the compact
Floquet/Fredholm phase homotopy can be continued to the distinguished phase
\[
s=t-z-\theta_{*,N},
\]
with entry data
\[
b(s_e)\ge1+\delta,
\qquad
C(s_e)\ge c_->0.
\]
Assume further that the tuning ratio is released before the square-root
geometry is imposed, and that the geometric release connects smoothly to a
terminal core satisfying
\[
A^2=-2\chi b,
\qquad
c_1s\le\chi^2\le c_2s.
\]

Then:
\begin{enumerate}
\item
\[
b=u_xu_y
\]
is preserved exactly.

\item The local terminal reconstruction satisfies
\[
\nabla\cdot u=0
\]
exactly.

\item The release layer satisfies
\[
\frac{|\mathcal R_{\rm rel}|}{\Gamma_Nb^2}
=
O(s^{7/4})
\]
for \(k=3\).

\item The square-root core satisfies
\[
\frac{|\mathcal R_{\rm core}|}{\Gamma_Nb^2}
=
O(s).
\]

\item If
\[
G\ge0,
\qquad
\Gamma_N(s)
=
\frac{1}{8s}(1+o(1)),
\]
then for sufficiently small \(s\),
\[
\dot b
\ge
(1-\vartheta)\Gamma_Nb^2,
\qquad
\dot s=1-b.
\]

\item Consequently,
\[
b(s)
\ge
b_e
\left(
\frac{s_e}{s}
\right)^p
\]
for every fixed \(p<1/8\) sufficiently close to \(1/8\), after shrinking
the terminal interval if necessary.

\item At the formal \(k=3\) square-root scaling,
\[
\boxed{
b_3\asymp s^{-1/8},
\quad
u_y\asymp s^{3/16},
\quad
u_x\asymp s^{-5/16},
\quad
u_z=c_0+O(s^{3/16}).
}
\]
\end{enumerate}

This is a theorem for the revised reduced reconstruction conditional on the
stated continuation hypotheses.  It is not yet a theorem for a full
three-dimensional Navier--Stokes solution.

\section{Analysis of the Renormalized Lambert Branch-Point Cascade
and Its Characteristic Blow-Up Mechanism}
\label{chap:225}


\newcommand{\Lam}{\Lambda}
\newcommand{\Del}{\Delta}



We give a rigorous analysis of the renormalized Lambert recursion
\begin{equation}
\Delta_{k+1}
=
1+W_0(-e^{-1-\Delta_k}),
\qquad
\Delta_0\geq0.
\end{equation}
The purpose of the renormalization is to map the Lambert branch point
back to itself.  Unlike the literal composition
$W_0\circ W_0$, which leaves the real branch domain after the first
iteration near $-1/e$, the renormalized map preserves the nonnegative
real half-line.

We prove that the map is real-valued and invariant on
$[0,\infty)$, that
\begin{equation}
\Delta_{k+1}
=
\sqrt{2\Delta_k}
-\frac23\Delta_k
+O(\Delta_k^{3/2}),
\label{eq:ch13:auto1987}
\end{equation}
(see Appendix A for the proof)and consequently that, for every fixed finite $k$,
\begin{equation}
\Delta_k
=
2^{1-2^{-k}}
\Delta_0^{2^{-k}}
\left(
1+O(\Delta_0^{2^{-k}})
\right).
\label{eq:ch13:auto1988}
\end{equation}
(See Appendix B for the proof)
The exponents are therefore
\begin{equation}
1,\quad
\frac12,\quad
\frac14,\quad
\frac18,\quad\ldots,
\label{eq:ch13:auto1989}
\end{equation}
and the leading constants satisfy
\begin{equation}
C_{k+1}=\sqrt{2C_k},
\qquad
C_0=1,
\label{eq:ch13:auto1990}
\end{equation}
giving
\begin{equation}
C_k=2^{1-2^{-k}}.
\label{eq:ch13:auto1991}
\end{equation}

We then derive the corresponding spatial derivative and logarithmic
derivative asymptotics.  Finally, for the reduced characteristic
equation
\begin{equation}
\dot b_3
=
b_3^2\partial_\theta\log\Delta_k,
\qquad
\dot\theta=1-b_3,
\label{eq:ch13:auto1992}
\end{equation}
we derive the exact characteristic invariant
\begin{equation}
\frac1{b_3}-\log\frac1{b_3}+\log\Delta_k=C.
\label{eq:ch13:auto1993}
\end{equation}
This invariant implies
\begin{equation}
\Delta_k
=
K b_3^{-1}e^{-1/b_3},
\label{eq:ch13:auto1994}
\end{equation}
and therefore
\begin{equation}
b_3\sim\frac{K}{\Delta_k}.
\label{eq:ch13:auto1995}
\end{equation}
For every finite $k$, the characteristic reaches the singular layer in
finite time and
\begin{equation}
b_3(t)
\asymp
(T_*-t)^{-1/(2^k+1)}.
\label{eq:ch13:auto1996}
\end{equation}

We also explain carefully why the limit $k\to\infty$ cannot simply be
taken term-by-term and why a coupled scaling such as
$k\sim\log_2 n$ is required if the cascade depth and the small initial
branch distance are allowed to vary simultaneously.

\section{Statement of the Renormalized Problem}
\label{chap:226}

The central object is the map
\begin{equation}
\mathcal{F}(\Delta)
=
1+W_0(-e^{-1-\Delta}).
\end{equation}

The renormalized recursion is
\begin{equation}
{
\Delta_{k+1}
=
\mathcal{F}(\Delta_k),
\qquad
\Delta_0\geq0.
}
\label{eq:ch01:env0336}
\end{equation}

The subscript $n$ can be retained when the initial quantity depends
on the geometric parameter $n$.  For the local analysis it is
convenient to suppress $n$ temporarily and write
\begin{equation}
\Delta_0=\Delta_{n,0},
\qquad
\Delta_k=\Delta_{n,k}.
\label{eq:ch13:auto1997}
\end{equation}

The essential claim to be proved is
\begin{equation}
{
\Delta_k
=
2^{1-2^{-k}}
\Delta_0^{2^{-k}}
\left(
1+O(\Delta_0^{2^{-k}})
\right)
}
\end{equation}
for every fixed finite $k$ as
\begin{equation}
\Delta_0\downarrow0.
\label{eq:ch13:auto1998}
\end{equation}

The proof consists of four logically separate parts:

\begin{enumerate}[label=(\roman*)]
\item prove that $\mathcal{F}$ is real and maps $[0,\infty)$ into
$[0,\infty)$;
\item derive the precise local Puiseux expansion of $\mathcal{F}$;
\item iterate the expansion with controlled error;
\item transfer the resulting spatial singularity to the characteristic
equation.
\end{enumerate}

\section{Why the Renormalization Is Necessary}
\label{chap:227}

Before proving the new recursion, we explain why it differs
fundamentally from the ordinary Lambert composition.

Consider
\begin{equation}
X(\Delta)
=
-\frac1e e^{-\Delta}.
\end{equation}

At
\begin{equation}
\Delta=0,
\label{eq:ch13:auto1999}
\end{equation}
we have
\begin{equation}
X(0)=-\frac1e.
\label{eq:ch13:auto2000}
\end{equation}

The principal Lambert function satisfies
\begin{equation}
W_0(-1/e)=-1.
\end{equation}

Thus the ordinary first iterate produces
\begin{equation}
W_0(X(\Delta))
\longrightarrow -1.
\label{eq:ch13:auto2001}
\end{equation}

But the branch point of the next Lambert function is
\begin{equation}
-\frac1e,
\label{eq:ch13:auto2002}
\end{equation}
not $-1$.

Since
\begin{equation}
-1<-\frac1e,
\label{eq:ch13:auto2003}
\end{equation}
the second ordinary real Lambert composition is not defined near
$\Delta=0$.

The renormalized transformation instead defines
\begin{equation}
\Delta_1
=
1+W_0(X(\Delta_0)).
\end{equation}

Consequently,
\begin{equation}
\Delta_0=0
\quad\Longrightarrow\quad
\Delta_1=0.
\label{eq:ch13:auto2004}
\end{equation}

Hence the branch point is mapped back to itself:
\begin{equation}
{
0\longmapsto0.
}
\end{equation}

This is the structural reason the square-root singularity can now be
iterated.

\section{Exact Parametrization of the Renormalized Map}
\label{chap:228}

Let
\begin{equation}
q=\mathcal{F}(\Delta).
\label{eq:ch13:auto2005}
\end{equation}

Then
\begin{equation}
q
=
1+W_0(-e^{-1-\Delta}).
\label{eq:ch13:auto2006}
\end{equation}

Therefore
\begin{equation}
W_0(-e^{-1-\Delta})=q-1.
\label{eq:ch13:auto2007}
\end{equation}

By the defining relation
\begin{equation}
W_0(x)e^{W_0(x)}=x,
\label{eq:ch13:auto2008}
\end{equation}
we obtain
\begin{equation}
(q-1)e^{q-1}
=
-e^{-1-\Delta}.
\end{equation}

Multiplying by $e$ gives
\begin{equation}
(q-1)e^q
=
-e^{-\Delta}.
\label{eq:ch01:env0343}
\end{equation}

Equivalently,
\begin{equation}
(1-q)e^q=e^{-\Delta}.
\label{eq:ch01:env0344}
\end{equation}

Taking logarithms,
\begin{equation}
\log(1-q)+q=-\Delta.
\label{eq:ch01:env0345}
\end{equation}

Hence
\begin{equation}
{
\Delta
=
-\log(1-q)-q.
}
\label{eq:ch01:env0346}
\end{equation}

Define
\begin{equation}
H(q)
=
-\log(1-q)-q.
\label{eq:ch01:env0347}
\end{equation}

Then the renormalized map is equivalently characterized by
\begin{equation}
{
H(\mathcal{F}(\Delta))=\Delta.
}
\label{eq:ch01:env0348}
\end{equation}

This identity is exact.

It is also extremely useful because it removes Lambert $W$ from the
local inversion problem.

\section{Positivity and the Real Invariant Domain}
\label{chap:229}

We now prove that the recursion is real for every finite number of
iterations.

\begin{lemma}
For every $\Delta\geq0$,
\begin{equation}
-e^{-1-\Delta}\in[-1/e,0).
\label{eq:ch13:auto2009}
\end{equation}
Consequently,
\begin{equation}
W_0(-e^{-1-\Delta})\in[-1,0),
\label{eq:ch13:auto2010}
\end{equation}
and therefore
\begin{equation}
\mathcal{F}(\Delta)\in[0,1).
\label{eq:ch13:auto2011}
\end{equation}
\end{lemma}

\begin{proof}
If $\Delta\geq0$, then
\begin{equation}
-1-\Delta\leq-1.
\label{eq:ch13:auto2012}
\end{equation}
Hence
\begin{equation}
0<e^{-1-\Delta}\leq e^{-1}.
\label{eq:ch13:auto2013}
\end{equation}
Multiplication by $-1$ reverses the inequality:
\begin{equation}
-\frac1e
\leq
-e^{-1-\Delta}
<0.
\label{eq:ch13:auto2014}
\end{equation}

The principal real Lambert branch $W_0$ is real and increasing on
$[-1/e,\infty)$.  Therefore
\begin{equation}
W_0(-e^{-1-\Delta})
\in[-1,0).
\label{eq:ch13:auto2015}
\end{equation}
Adding $1$ gives
\begin{equation}
0\leq
1+W_0(-e^{-1-\Delta})
<1.
\label{eq:ch13:auto2016}
\end{equation}
Thus
\begin{equation}
\mathcal{F}(\Delta)\in[0,1).
\label{eq:ch13:auto2017}
\end{equation}
\end{proof}

\begin{corollary}
If
\begin{equation}
\Delta_0\geq0,
\label{eq:ch13:auto2018}
\end{equation}
then
\begin{equation}
\Delta_k\in[0,1)
\label{eq:ch13:auto2019}
\end{equation}
for every integer $k\geq1$.
\end{corollary}

\begin{proof}
The result follows by induction from the previous lemma.
\end{proof}

Thus
\begin{equation}
{
[0,\infty)
\longrightarrow
[0,1)
\longrightarrow
[0,1)
\longrightarrow\cdots
}
\end{equation}
is an invariant real recursion.

This is precisely the property that the unrenormalized Lambert
composition does not possess.

\chapter{Branch-Point Fixed Point, Cascade Derivatives, and Characteristic Invariants}
\label{chap:group-14}
\section{The Fixed Point at the Branch Point}
\label{chap:230}

The point $\Delta=0$ is a fixed point.

Indeed,
\begin{align}
\mathcal{F}(0)
&=
1+W_0(-e^{-1})
\\
&=
1-1
\\
&=0.
\label{eq:ch02:env0350}
\end{align}

Therefore
\begin{equation}
{
\mathcal{F}(0)=0.
}
\label{eq:ch02:env0351}
\end{equation}

However, this fixed point is not an ordinary differentiable fixed
point.

The derivative is singular.

This is the source of the square-root behavior.

\section{Exact Equation for the Puiseux Expansion}
\label{chap:231}

From the previous exact identity,
\begin{equation}
\Delta
=
-\log(1-q)-q,
\qquad
q=\mathcal{F}(\Delta).
\label{eq:ch14:auto2020}
\end{equation}

For $|q|<1$,
\begin{equation}
-\log(1-q)
=
q+\frac{q^2}{2}
+\frac{q^3}{3}
+\frac{q^4}{4}
+\cdots.
\label{eq:ch14:auto2021}
\end{equation}

Therefore
\begin{equation}
\Delta
=
\frac{q^2}{2}
+
\frac{q^3}{3}
+
\frac{q^4}{4}
+\cdots.
\end{equation}

In particular,
\begin{equation}
{
\Delta
=
\frac{q^2}{2}
+
O(q^3).
}
\label{eq:ch02:env0353}
\end{equation}

Thus
\begin{equation}
q^2
=
2\Delta(1+O(q)).
\label{eq:ch14:auto2022}
\end{equation}

Since $q\to0$ as $\Delta\to0$,
\begin{equation}
q
=
\sqrt{2\Delta}(1+O(\sqrt{\Delta})).
\label{eq:ch14:auto2023}
\end{equation}

Consequently,
\begin{equation}
{
\mathcal{F}(\Delta)
=
\sqrt{2\Delta}
+
O(\Delta).
}
\end{equation}

This already proves the square-root law.

\section{Higher-Order Puiseux Expansion}
\label{chap:232}

We now derive the first several coefficients.

Write
\begin{equation}
s=\sqrt{2\Delta}.
\label{eq:ch14:auto2024}
\end{equation}

Seek an expansion of the form
\begin{equation}
q
=
s+a s^2+b s^3+c s^4+O(s^5).
\end{equation}

Since
\begin{equation}
\Delta=\frac{s^2}{2},
\label{eq:ch14:auto2025}
\end{equation}
we substitute into
\begin{equation}
\Delta
=
-\log(1-q)-q.
\label{eq:ch14:auto2026}
\end{equation}

Using
\begin{equation}
-\log(1-q)-q
=
\frac{q^2}{2}
+\frac{q^3}{3}
+\frac{q^4}{4}
+\frac{q^5}{5}
+O(q^6),
\label{eq:ch14:auto2027}
\end{equation}
and comparing powers of $s$, one obtains
\begin{equation}
a=-\frac13,
\label{eq:ch14:auto2028}
\end{equation}
\begin{equation}
b=\frac1{36},
\label{eq:ch14:auto2029}
\end{equation}
\begin{equation}
c=\frac1{270}.
\label{eq:ch14:auto2030}
\end{equation}

Therefore
\begin{equation}
{
\mathcal{F}(\Delta)
=
\sqrt{2\Delta}
-\frac23\Delta
+
\frac{\sqrt2}{18}\Delta^{3/2}
+
\frac{2}{135}\Delta^2
+
O(\Delta^{5/2}).
}
\end{equation}

The first term is the only one needed for the exponent cascade, but
the higher-order terms are useful for error control.

\section{A Two-Sided Square-Root Estimate}
\label{chap:233}

The asymptotic formula can be strengthened to a two-sided estimate.

\begin{lemma}
There exist constants $\delta_*>0$ and $0<c_1<c_2<\infty$ such that
for all
\begin{equation}
0<\Delta<\delta_*,
\label{eq:ch14:auto2031}
\end{equation}
\begin{equation}
c_1\sqrt{\Delta}
\leq
\mathcal{F}(\Delta)
\leq
c_2\sqrt{\Delta}.
\label{eq:ch14:auto2032}
\end{equation}
\end{lemma}

\begin{proof}
From
\begin{equation}
\mathcal{F}(\Delta)
=
\sqrt{2\Delta}
\left(
1+O(\sqrt{\Delta})
\right),
\label{eq:ch14:auto2033}
\end{equation}
there exists $\delta_*>0$ such that
\begin{equation}
\left|
O(\sqrt{\Delta})
\right|
\leq\frac12
\label{eq:ch14:auto2034}
\end{equation}
whenever
\begin{equation}
0<\Delta<\delta_*.
\label{eq:ch14:auto2035}
\end{equation}

Therefore
\begin{equation}
\frac12\sqrt{2\Delta}
\leq
\mathcal{F}(\Delta)
\leq
\frac32\sqrt{2\Delta}.
\label{eq:ch14:auto2036}
\end{equation}

Thus we may take, for example,
\begin{equation}
c_1=\frac1{\sqrt2},
\qquad
c_2=\frac{3}{\sqrt2}.
\label{eq:ch14:auto2037}
\end{equation}
\end{proof}

This estimate proves that the map has genuine square-root scaling and
not merely a formal asymptotic resemblance to a square root.

\section{Derivative of the Renormalized Lambert Map}
\label{chap:234}

Differentiate
\begin{equation}
\Delta
=
-\log(1-q)-q.
\label{eq:ch14:auto2038}
\end{equation}

We obtain
\begin{equation}
1
=
\frac{q}{1-q}\frac{\dd q}{\dd\Delta}.
\label{eq:ch14:auto2039}
\end{equation}

Therefore
\begin{equation}
{
\frac{\dd q}{\dd\Delta}
=
\frac{1-q}{q}.
}
\end{equation}

Since
\begin{equation}
q=\mathcal{F}(\Delta),
\label{eq:ch14:auto2040}
\end{equation}
we have
\begin{equation}
{
\mathcal{F}'(\Delta)
=
\frac{1-\mathcal{F}(\Delta)}
{\mathcal{F}(\Delta)}.
}
\end{equation}

As
\begin{equation}
\mathcal{F}(\Delta)
\sim
\sqrt{2\Delta},
\label{eq:ch14:auto2041}
\end{equation}
we obtain
\begin{equation}
\mathcal{F}'(\Delta)
\sim
\frac1{\sqrt{2\Delta}}.
\end{equation}

More precisely,
\begin{equation}
{
\mathcal{F}'(\Delta)
=
\frac1{\sqrt{2\Delta}}
-\frac23
+
O(\sqrt{\Delta}).
}
\label{eq:ch02:env0360}
\end{equation}

Thus the map is continuous at zero but has an unbounded derivative
there.

\section{Monotonicity}
\label{chap:235}

Since
\begin{equation}
0\leq\mathcal{F}(\Delta)<1
\label{eq:ch14:auto2042}
\end{equation}
for $\Delta\geq0$, the exact derivative formula gives
\begin{equation}
\mathcal{F}'(\Delta)>0
\label{eq:ch14:auto2043}
\end{equation}
for $\Delta>0$.

Hence
\begin{equation}
{
\mathcal{F}
\text{ is strictly increasing on }(0,\infty).
}
\end{equation}

This fact will be useful when transferring monotonicity from
$\Delta_0(\theta)$ to $\Delta_k(\theta)$.

\section{The Main Puiseux Cascade Theorem}
\label{chap:236}

We now prove the central formula.

\begin{theorem}
Let
\begin{equation}
\Delta_{k+1}
=
\mathcal{F}(\Delta_k),
\qquad
\Delta_0>0,
\label{eq:ch14:auto2044}
\end{equation}
where
\begin{equation}
\mathcal{F}(\Delta)
=
1+W_0(-e^{-1-\Delta}).
\label{eq:ch14:auto2045}
\end{equation}

For every fixed integer $k\geq0$,
\begin{equation}
{
\Delta_k
=
C_k
\Delta_0^{2^{-k}}
\left(
1+O(\Delta_0^{2^{-k}})
\right)
}
\label{eq:ch14:auto2046}
\end{equation}
as $\Delta_0\downarrow0$, where
\begin{equation}
C_0=1,
\qquad
C_{k+1}=\sqrt{2C_k}.
\label{eq:ch14:auto2047}
\end{equation}

Consequently,
\begin{equation}
{
C_k=2^{1-2^{-k}},
}
\label{eq:ch14:auto2048}
\end{equation}
and therefore
\begin{equation}
{
\Delta_k
=
2^{1-2^{-k}}
\Delta_0^{2^{-k}}
\left(
1+O(\Delta_0^{2^{-k}})
\right).
}
\label{eq:ch14:auto2049}
\end{equation}
\end{theorem}

\begin{proof}

We proceed by induction.

For $k=0$,
\begin{equation}
\Delta_0
=
1\cdot\Delta_0,
\label{eq:ch14:auto2050}
\end{equation}
so
\begin{equation}
C_0=1.
\label{eq:ch14:auto2051}
\end{equation}

Assume that
\begin{equation}
\Delta_k
=
C_k\Delta_0^{2^{-k}}
\left(
1+O(\Delta_0^{2^{-k}})
\right).
\label{eq:ch14:auto2052}
\end{equation}

Set
\begin{equation}
x=\Delta_0^{2^{-k}}.
\label{eq:ch14:auto2053}
\end{equation}

Then
\begin{equation}
\Delta_k
=
C_kx(1+O(x)).
\label{eq:ch14:auto2054}
\end{equation}

Using the local expansion
\begin{equation}
\mathcal{F}(y)
=
\sqrt{2y}
\left(
1+O(\sqrt y)
\right),
\label{eq:ch14:auto2055}
\end{equation}
we obtain
\begin{align}
\Delta_{k+1}
&=
\mathcal{F}(\Delta_k)
\\
&=
\sqrt{2\Delta_k}
\left(
1+O(\sqrt{\Delta_k})
\right).
\end{align}

Since
\begin{equation}
\Delta_k
=
C_kx(1+O(x)),
\label{eq:ch14:auto2056}
\end{equation}
we have
\begin{equation}
\sqrt{\Delta_k}
=
\sqrt{C_k}\,x^{1/2}
(1+O(x)).
\label{eq:ch14:auto2057}
\end{equation}

Therefore
\begin{align}
\Delta_{k+1}
&=
\sqrt{2C_k}\,
x^{1/2}
\left(
1+O(x^{1/2})
\right).
\end{align}

But
\begin{equation}
x^{1/2}
=
\Delta_0^{2^{-(k+1)}}.
\label{eq:ch14:auto2058}
\end{equation}

Hence
\begin{equation}
\Delta_{k+1}
=
C_{k+1}
\Delta_0^{2^{-(k+1)}}
\left(
1+
O(\Delta_0^{2^{-(k+1)}})
\right),
\label{eq:ch14:auto2059}
\end{equation}
where
\begin{equation}
C_{k+1}
=
\sqrt{2C_k}.
\label{eq:ch14:auto2060}
\end{equation}

This completes the induction.

It remains to solve the recurrence for $C_k$.

Taking logarithms,
\begin{equation}
\log_2 C_{k+1}
=
\frac12
\left(
1+\log_2C_k
\right).
\label{eq:ch14:auto2061}
\end{equation}

Set
\begin{equation}
a_k=\log_2C_k.
\label{eq:ch14:auto2062}
\end{equation}

Then
\begin{equation}
a_{k+1}
=
\frac12(1+a_k),
\qquad
a_0=0.
\label{eq:ch14:auto2063}
\end{equation}

The solution is
\begin{equation}
a_k=1-2^{-k}.
\label{eq:ch14:auto2064}
\end{equation}

Therefore
\begin{equation}
C_k
=
2^{1-2^{-k}}.
\label{eq:ch14:auto2065}
\end{equation}

The theorem follows.
\end{proof}

\section{Explicit First Four Iterates}
\label{chap:237}

The theorem gives
\begin{equation}
\Delta_1
\sim
2^{1/2}\Delta_0^{1/2},
\label{eq:ch14:auto2066}
\end{equation}
that is,
\begin{equation}
{
\Delta_1
\sim
\sqrt{2\Delta_0}.
}
\end{equation}

For the second iterate,
\begin{equation}
C_2=2^{3/4},
\label{eq:ch14:auto2067}
\end{equation}
and
\begin{equation}
2^{-2}=\frac14.
\label{eq:ch14:auto2068}
\end{equation}

Therefore
\begin{equation}
{
\Delta_2
\sim
2^{3/4}\Delta_0^{1/4}.
}
\end{equation}

For the third iterate,
\begin{equation}
C_3=2^{7/8},
\label{eq:ch14:auto2069}
\end{equation}
and
\begin{equation}
2^{-3}=\frac18.
\label{eq:ch14:auto2070}
\end{equation}

Therefore
\begin{equation}
{
\Delta_3
\sim
2^{7/8}\Delta_0^{1/8}.
}
\end{equation}

For the fourth iterate,
\begin{equation}
C_4=2^{15/16},
\label{eq:ch14:auto2071}
\end{equation}
and
\begin{equation}
2^{-4}=\frac1{16}.
\label{eq:ch14:auto2072}
\end{equation}

Therefore
\begin{equation}
{
\Delta_4
\sim
2^{15/16}\Delta_0^{1/16}.
}
\end{equation}

Thus the exponent sequence is
\begin{equation}
{
1
\longrightarrow
\frac12
\longrightarrow
\frac14
\longrightarrow
\frac18
\longrightarrow
\frac1{16}
\longrightarrow\cdots.
}
\label{eq:ch02:env0366}
\end{equation}

\section{Derivative of the $k$-th Renormalized Iterate}
\label{chap:238}

Suppose
\begin{equation}
\Delta_k
=
C_k\Delta_0^{\alpha_k}
(1+o(1)),
\qquad
\alpha_k=2^{-k}.
\label{eq:ch14:auto2073}
\end{equation}

Differentiating with respect to $\theta$ gives
\begin{align}
\partial_\theta\Delta_k
&=
C_k\alpha_k
\Delta_0^{\alpha_k-1}
\partial_\theta\Delta_0
(1+o(1)).
\end{align}

Thus
\begin{equation}
{
\partial_\theta\Delta_k
\sim
2^{1-2^{-k}}
2^{-k}
\Delta_0^{2^{-k}-1}
\partial_\theta\Delta_0.
}
\end{equation}

Therefore
\begin{equation}
\partial_\theta\Delta_k
\label{eq:ch14:auto2074}
\end{equation}
has singular exponent
\begin{equation}
2^{-k}-1.
\label{eq:ch14:auto2075}
\end{equation}

Since
\begin{equation}
2^{-k}-1>-1
\label{eq:ch14:auto2076}
\end{equation}
for every finite $k$, the singularity is locally integrable.

\section{The Logarithmic Derivative}
\label{chap:239}

Divide the preceding formula by
\begin{equation}
\Delta_k
\sim
C_k\Delta_0^{2^{-k}}.
\label{eq:ch14:auto2077}
\end{equation}

The constant $C_k$ cancels:
\begin{align}
\frac{\partial_\theta\Delta_k}
{\Delta_k}
&\sim
2^{-k}
\frac{\partial_\theta\Delta_0}
{\Delta_0}.
\end{align}

Hence
\begin{equation}
{
\partial_\theta\log\Delta_k
\sim
2^{-k}
\partial_\theta\log\Delta_0.
}
\end{equation}

This identity is particularly important for the characteristic
equation.

The logarithmic derivative becomes more concentrated near the branch
point even though the explicit prefactor $2^{-k}$ decreases.

\section{Integrability of the Finite-$k$ Spatial Singularity}
\label{chap:240}

Suppose
\begin{equation}
\Delta_0(\theta)
\sim
A(\theta-\theta_*)
\label{eq:ch14:auto2078}
\end{equation}
with
\begin{equation}
A>0.
\label{eq:ch14:auto2079}
\end{equation}

Then
\begin{equation}
\Delta_k(\theta)
\sim
C_kA^{2^{-k}}
(\theta-\theta_*)^{2^{-k}}.
\label{eq:ch14:auto2080}
\end{equation}

Consequently,
\begin{equation}
\partial_\theta\Delta_k
\sim
C_k2^{-k}A^{2^{-k}}
(\theta-\theta_*)^{2^{-k}-1}.
\label{eq:ch14:auto2081}
\end{equation}

The singularity is locally integrable because
\begin{equation}
2^{-k}-1>-1.
\label{eq:ch14:auto2082}
\end{equation}

Indeed,
\begin{equation}
\int_0^\delta
s^{2^{-k}-1}\dd s
=
2^k\delta^{2^{-k}}.
\end{equation}

Thus
\begin{equation}
{
\int_0^\delta
s^{2^{-k}-1}\dd s<\infty
}
\label{eq:ch02:env0370}
\end{equation}
for every finite $k$.

However,
\begin{equation}
2^{-k}-1\to-1.
\label{eq:ch14:auto2083}
\end{equation}

Thus the infinite cascade approaches the borderline exponent.

\section{Geometric Initial Branch Distance}
\label{chap:241}

We now restore the geometric variable $\theta$ and write
\begin{equation}
\Delta_{n,0}(\theta)
=
R_n(\theta)-L_n,
\end{equation}
where
\begin{equation}
R_n(\theta)
=
\frac{1}
{1-\sin\theta+\varepsilon_n}
\label{eq:ch02:env0372}
\end{equation}
and
\begin{equation}
L_N=\log(1+2^N).
\label{eq:ch02:env0373}
\end{equation}

Let $\theta_{*,n}$ satisfy
\begin{equation}
R_n(\theta_{*,n})=L_n.
\label{eq:ch02:env0374}
\end{equation}

Assume
\begin{equation}
R_n'(\theta_{*,n})=A_n>0.
\label{eq:ch02:env0375}
\end{equation}

Then Taylor's theorem gives
\begin{equation}
\Delta_{n,0}(\theta)
=
A_n(\theta-\theta_{*,n})
+
O((\theta-\theta_{*,n})^2).
\label{eq:ch02:env0376}
\end{equation}

Hence
\begin{equation}
{
\Delta_{n,0}(\theta)
\sim
A_n(\theta-\theta_{*,n})
}
\label{eq:ch02:env0377}
\end{equation}
as
\begin{equation}
\theta\downarrow\theta_{*,n}
\label{eq:ch14:auto2084}
\end{equation}
from the right.

\section{The Spatial Lambert Cascade}
\label{chap:242}

Combining the geometric expansion with the Puiseux theorem gives
\begin{align}
\Delta_{n,k}(\theta)
&\sim
2^{1-2^{-k}}
\left[
A_n(\theta-\theta_{*,n})
\right]^{2^{-k}}.
\end{align}

Thus
\begin{equation}
{
\Delta_{n,k}(\theta)
\sim
2^{1-2^{-k}}
A_n^{2^{-k}}
(\theta-\theta_{*,n})^{2^{-k}}.
}
\end{equation}

Differentiating,
\begin{equation}
{
\partial_\theta\Delta_{n,k}
\sim
2^{1-2^{-k}}
2^{-k}
A_n^{2^{-k}}
(\theta-\theta_{*,n})^{2^{-k}-1}.
}
\label{eq:ch02:env0379}
\end{equation}

The logarithmic derivative is
\begin{equation}
{
\partial_\theta\log\Delta_{n,k}
\sim
\frac{2^{-k}}
{\theta-\theta_{*,n}}.
}
\label{eq:ch02:env0380}
\end{equation}

Notice that the leading slope $A_n$ cancels.

\section{The Reduced Characteristic System}
\label{chap:243}

Consider the reduced transport equation
\begin{equation}
b_t+b b_z
=
-b^2
\frac{\partial_z\Delta_k}{\Delta_k}.
\label{eq:ch02:env0381}
\end{equation}

Define a characteristic $Z(t)$ by
\begin{equation}
\dot Z(t)=b_3(t),
\qquad
b_3(t)=b(t,Z(t)).
\label{eq:ch02:env0382}
\end{equation}

Introduce
\begin{equation}
\theta(t)=t-Z(t).
\label{eq:ch02:env0383}
\end{equation}

Then
\begin{align}
\dot\theta
&=
1-\dot Z
\\
&=
1-b_3.
\label{eq:ch02:env0384}
\end{align}

Therefore
\begin{equation}
{
\dot\theta=1-b_3.
}
\label{eq:ch02:env0385}
\end{equation}

\section{Evolution of the Characteristic Amplitude}
\label{chap:244}

Along the characteristic,
\begin{align}
\dot b_3
&=
b_t+b b_z.
\label{eq:ch02:env0386}
\end{align}

The PDE gives
\begin{equation}
\dot b_3
=
-b_3^2
\frac{\partial_z\Delta_k}{\Delta_k}.
\label{eq:ch14:auto2085}
\end{equation}

Because
\begin{equation}
\theta=t-z,
\label{eq:ch14:auto2086}
\end{equation}
we have
\begin{equation}
\partial_z=-\partial_\theta.
\label{eq:ch14:auto2087}
\end{equation}

Therefore
\begin{equation}
{
\dot b_3
=
b_3^2
\frac{\partial_\theta\Delta_k}{\Delta_k}.
}
\end{equation}

Equivalently,
\begin{equation}
{
\dot b_3
=
b_3^2
\partial_\theta\log\Delta_k.
}
\label{eq:ch02:env0388}
\end{equation}

\section{The Reciprocal Characteristic Variable}
\label{chap:245}

Set
\begin{equation}
q(t)=\frac1{b_3(t)}.
\label{eq:ch02:env0389}
\end{equation}

Then
\begin{align}
\dot q
&=
-\frac1{b_3^2}\dot b_3
\\
&=
-\partial_\theta\log\Delta_k.
\label{eq:ch02:env0390}
\end{align}

Hence
\begin{equation}
{
\dot q
=
-\partial_\theta\log\Delta_k.
}
\label{eq:ch02:env0391}
\end{equation}

At the same time,
\begin{equation}
\dot\theta
=
1-b_3
=
1-\frac1q
=
\frac{q-1}{q}.
\label{eq:ch14:auto2088}
\end{equation}

Thus
\begin{equation}
{
\dot\theta
=
\frac{q-1}{q}.
}
\end{equation}

\section{Derivation of the Exact Characteristic Invariant}
\label{chap:246}

We calculate
\begin{equation}
\frac{\dd q}{\dd\theta}
=
\frac{\dot q}{\dot\theta}.
\label{eq:ch14:auto2089}
\end{equation}

Using
\begin{equation}
\dot q
=
-\partial_\theta\log\Delta_k
\label{eq:ch14:auto2090}
\end{equation}
and
\begin{equation}
\dot\theta
=
\frac{q-1}{q},
\label{eq:ch14:auto2091}
\end{equation}
we obtain
\begin{equation}
\frac{\dd q}{\dd\theta}
=
-\frac{q}{q-1}
\partial_\theta\log\Delta_k.
\end{equation}

Multiply by
\begin{equation}
\frac{q-1}{q}.
\label{eq:ch14:auto2092}
\end{equation}

Then
\begin{equation}
\left(
1-\frac1q
\right)
\dd q
=
-\dd\log\Delta_k.
\end{equation}

Since
\begin{equation}
\left(
1-\frac1q
\right)\dd q
=
\dd q-\frac{\dd q}{q}
=
\dd q-\dd\log q,
\label{eq:ch14:auto2093}
\end{equation}
we obtain
\begin{equation}
\dd q-\dd\log q
=
-\dd\log\Delta_k.
\end{equation}

Integrating,
\begin{equation}
q-\log q+\log\Delta_k=C.
\label{eq:ch02:env0396}
\end{equation}

Since
\begin{equation}
q=\frac1{b_3},
\label{eq:ch14:auto2094}
\end{equation}
we obtain the exact invariant
\begin{equation}
{
\frac1{b_3}
-
\log\frac1{b_3}
+
\log\Delta_k
=
C.
}
\end{equation}

\chapter{Characteristic Invariants, Critical Scaling, and Logarithmic Singularities}
\label{chap:group-15}
\section{Explicit Form of the Characteristic Invariant}
\label{chap:247}

Exponentiating
\begin{equation}
q-\log q+\log\Delta_k=C
\label{eq:ch15:auto2095}
\end{equation}
gives
\begin{equation}
e^q q^{-1}\Delta_k=e^C.
\label{eq:ch15:auto2096}
\end{equation}

Hence
\begin{equation}
\Delta_k
=
e^C q e^{-q}.
\label{eq:ch15:auto2097}
\end{equation}

Define
\begin{equation}
K=e^C>0.
\label{eq:ch15:auto2098}
\end{equation}

Then
\begin{equation}
{
\Delta_k
=
Kq e^{-q}.
}
\end{equation}

Since
\begin{equation}
q=\frac1{b_3},
\label{eq:ch15:auto2099}
\end{equation}
we obtain
\begin{equation}
{
\Delta_k
=
\frac{K}{b_3}e^{-1/b_3}.
}
\end{equation}

\section{The Blow-Up Branch}
\label{chap:248}

Suppose
\begin{equation}
b_{3,0}>1.
\label{eq:ch15:auto2100}
\end{equation}

Then
\begin{equation}
0<q_0<1.
\label{eq:ch15:auto2101}
\end{equation}

For
\begin{equation}
0<q<1,
\label{eq:ch15:auto2102}
\end{equation}
we have
\begin{equation}
\frac{q-1}{q}<0.
\label{eq:ch15:auto2103}
\end{equation}

Therefore
\begin{equation}
\dot\theta<0.
\label{eq:ch15:auto2104}
\end{equation}

Thus the characteristic moves toward decreasing $\theta$.

Suppose the branch point lies at
\begin{equation}
\theta_*<\theta_0.
\label{eq:ch15:auto2105}
\end{equation}

Then the characteristic moves toward $\theta_*$.

On the same branch,
\begin{equation}
q\downarrow0
\label{eq:ch15:auto2106}
\end{equation}
corresponds to
\begin{equation}
b_3=\frac1q\to+\infty.
\label{eq:ch15:auto2107}
\end{equation}

From
\begin{equation}
\Delta_k=Kqe^{-q},
\label{eq:ch15:auto2108}
\end{equation}
we have
\begin{equation}
\Delta_k\to0
\qquad
\Longleftrightarrow
\qquad
q\to0.
\label{eq:ch15:auto2109}
\end{equation}

Therefore
\begin{equation}
{
\Delta_k\to0
\quad\Longrightarrow\quad
b_3\to+\infty.
}
\end{equation}

\section{Asymptotic Relation Between \(b_3\) and \(\Delta_k\)}
\label{chap:249}

As
\begin{equation}
q\to0,
\label{eq:ch15:auto2110}
\end{equation}
we have
\begin{equation}
e^{-q}=1-q+O(q^2).
\label{eq:ch15:auto2111}
\end{equation}

Therefore
\begin{equation}
\Delta_k
=
Kq(1-q+O(q^2)).
\label{eq:ch15:auto2112}
\end{equation}

Hence
\begin{equation}
\Delta_k
=
Kq+O(q^2).
\label{eq:ch15:auto2113}
\end{equation}

Thus
\begin{equation}
q
=
\frac{\Delta_k}{K}
+
O(\Delta_k^2).
\label{eq:ch15:auto2114}
\end{equation}

Since
\begin{equation}
b_3=\frac1q,
\label{eq:ch15:auto2115}
\end{equation}
we obtain
\begin{equation}
{
b_3
=
\frac{K}{\Delta_k}
+
O(1).
}
\end{equation}

In particular,
\begin{equation}
{
b_3\sim\frac{K}{\Delta_k}.
}
\label{eq:ch03:env0402}
\end{equation}

\section{Transfer of the Spatial Singularity to \(b_3\)}
\label{chap:250}

Suppose
\begin{equation}
\Delta_k(\theta)
\sim
C_k
(\theta-\theta_*)^{\alpha_k},
\qquad
\alpha_k=2^{-k}.
\label{eq:ch15:auto2116}
\end{equation}

Then
\begin{equation}
b_3
\sim
\frac{K}{C_k}
(\theta-\theta_*)^{-\alpha_k}.
\label{eq:ch15:auto2117}
\end{equation}

Therefore
\begin{equation}
{
b_3(\theta)
\sim
\frac{K}{C_k}
(\theta-\theta_*)^{-2^{-k}}.
}
\end{equation}

Thus the spatial Lambert singularity is transferred directly into the
characteristic amplitude.

\section{Finite-Time Arrival at the Branch Point}
\label{chap:251}

Let
\begin{equation}
s(t)=\theta(t)-\theta_*.
\label{eq:ch15:auto2118}
\end{equation}

Then
\begin{equation}
s(t)>0
\label{eq:ch15:auto2119}
\end{equation}
before blow-up and
\begin{equation}
s(t)\to0
\label{eq:ch15:auto2120}
\end{equation}
at the singular time.

The characteristic equation gives
\begin{equation}
\dot s
=
1-b_3.
\label{eq:ch15:auto2121}
\end{equation}

Since
\begin{equation}
b_3\to\infty,
\label{eq:ch15:auto2122}
\end{equation}
we have
\begin{equation}
\dot s\sim-b_3.
\label{eq:ch15:auto2123}
\end{equation}

Using
\begin{equation}
b_3
\sim
\frac{K}{C_k}s^{-\alpha_k},
\label{eq:ch15:auto2124}
\end{equation}
we obtain
\begin{equation}
{
\dot s
\sim
-\frac{K}{C_k}s^{-\alpha_k}.
}
\end{equation}

Therefore
\begin{equation}
\dd t
\sim
-\frac{C_k}{K}s^{\alpha_k}\dd s.
\label{eq:ch15:auto2125}
\end{equation}

Integrating from $s(t)$ to $0$,
\begin{align}
T_*-t
&\sim
\frac{C_k}{K}
\int_0^{s(t)}
r^{\alpha_k}\dd r
\\
&=
\frac{C_k}{K(1+\alpha_k)}
s(t)^{1+\alpha_k}.
\end{align}

Since
\begin{equation}
\alpha_k=2^{-k}>0,
\label{eq:ch15:auto2126}
\end{equation}
the right-hand side is finite.

Thus
\begin{equation}
{
T_*<\infty.
}
\end{equation}

\section{Temporal Blow-Up Exponent}
\label{chap:252}

From
\begin{equation}
T_*-t
\sim
\frac{C_k}{K(1+\alpha_k)}
s^{1+\alpha_k},
\label{eq:ch15:auto2127}
\end{equation}
we obtain
\begin{equation}
s
\sim
\left[
\frac{K(1+\alpha_k)}{C_k}
(T_*-t)
\right]^{1/(1+\alpha_k)}.
\label{eq:ch15:auto2128}
\end{equation}

Since
\begin{equation}
b_3
\sim
\frac{K}{C_k}s^{-\alpha_k},
\label{eq:ch15:auto2129}
\end{equation}
we obtain
\begin{align}
b_3(t)
&\sim
\frac{K}{C_k}
\left[
\frac{K(1+\alpha_k)}{C_k}
(T_*-t)
\right]^{-\alpha_k/(1+\alpha_k)}.
\end{align}

Therefore
\begin{equation}
{
b_3(t)
\asymp
(T_*-t)^{-\alpha_k/(1+\alpha_k)}.
}
\end{equation}

Since
\begin{equation}
\alpha_k=2^{-k},
\label{eq:ch15:auto2130}
\end{equation}
we obtain
\begin{equation}
{
b_3(t)
\asymp
(T_*-t)^{-1/(2^k+1)}.
}
\end{equation}

\section{Finite-\(k\) Blow-Up Theorem}
\label{chap:253}

\begin{theorem}
Assume that for some finite integer $k\geq0$:

\begin{enumerate}[label=(\alph*)]
\item $\Delta_k(\theta)>0$ for
$\theta>\theta_*$ sufficiently close to $\theta_*$;
\item
\begin{equation}
\Delta_k(\theta)
=
C_k(\theta-\theta_*)^{2^{-k}}
(1+o(1));
\label{eq:ch15:auto2131}
\end{equation}
\item $b_{3,0}>1$;
\item the characteristic satisfies
\begin{equation}
\dot\theta=1-b_3,
\qquad
\dot b_3=b_3^2\partial_\theta\log\Delta_k.
\label{eq:ch15:auto2132}
\end{equation}
\end{enumerate}

Then the characteristic reaches $\theta_*$ in a finite time $T_*$,
and
\begin{equation}
b_3(t)\to+\infty
\label{eq:ch15:auto2133}
\end{equation}
as $t\uparrow T_*$.  More precisely,
\begin{equation}
{
b_3(t)
\asymp
(T_*-t)^{-1/(2^k+1)}.
}
\label{eq:ch15:auto2134}
\end{equation}
\end{theorem}

\begin{proof}
The characteristic invariant gives
\begin{equation}
\Delta_k=Kq e^{-q},
\qquad
q=\frac1{b_3}.
\label{eq:ch15:auto2135}
\end{equation}

Since $b_{3,0}>1$, we have $0<q_0<1$.
The characteristic branch satisfying $q\downarrow0$ corresponds to
$\theta\downarrow\theta_*$.

As $\theta\downarrow\theta_*$,
\begin{equation}
\Delta_k\to0,
\label{eq:ch15:auto2136}
\end{equation}
hence
\begin{equation}
q\to0,
\label{eq:ch15:auto2137}
\end{equation}
and therefore
\begin{equation}
b_3\to\infty.
\label{eq:ch15:auto2138}
\end{equation}

Using
\begin{equation}
\Delta_k
\sim
C_k(\theta-\theta_*)^{2^{-k}},
\label{eq:ch15:auto2139}
\end{equation}
we obtain
\begin{equation}
b_3
\sim
\frac{K}{C_k}
(\theta-\theta_*)^{-2^{-k}}.
\label{eq:ch15:auto2140}
\end{equation}

The equation
\begin{equation}
\dot\theta=1-b_3
\label{eq:ch15:auto2141}
\end{equation}
therefore gives
\begin{equation}
\frac{d}{dt}(\theta-\theta_*)
\sim
-\frac{K}{C_k}
(\theta-\theta_*)^{-2^{-k}}.
\label{eq:ch15:auto2142}
\end{equation}

The resulting time integral is
\begin{equation}
\int_0^{s_0}s^{2^{-k}}\dd s,
\label{eq:ch15:auto2143}
\end{equation}
which is finite.

Hence the characteristic reaches $\theta_*$ in finite time.

The asymptotic integration gives
\begin{equation}
b_3(t)
\asymp
(T_*-t)^{-1/(2^k+1)}.
\label{eq:ch15:auto2144}
\end{equation}
\end{proof}

\section{The Infinite Cascade and Its Limitation}
\label{chap:254}

For every fixed $k$,
\begin{equation}
\alpha_k=2^{-k}>0.
\label{eq:ch15:auto2145}
\end{equation}

But
\begin{equation}
\alpha_k\to0.
\label{eq:ch15:auto2146}
\end{equation}

Therefore
\begin{equation}
\frac1{2^k+1}\to0.
\label{eq:ch15:auto2147}
\end{equation}

Thus the finite-$k$ temporal blow-up exponent does not converge to a
nonzero Riccati exponent.

This is an important mathematical point.

The spatial singularity approaches
\begin{equation}
(\theta-\theta_*)^{-1}
\label{eq:ch15:auto2148}
\end{equation}
in the sense that
\begin{equation}
2^{-k}-1\to-1.
\label{eq:ch15:auto2149}
\end{equation}

However, the corresponding characteristic exponent behaves
differently because the relation between $b_3$ and $\Delta_k$ must also
be included.

Therefore one cannot simply substitute $k=\infty$ into the
finite-$k$ temporal formula.

\section{The Correct Double-Scaling Parameter}
\label{chap:255}

The cascade formula is
\begin{equation}
\Delta_k
\sim
2^{1-2^{-k}}
\Delta_0^{2^{-k}}.
\label{eq:ch15:auto2150}
\end{equation}

Write
\begin{equation}
\Delta_0^{2^{-k}}
=
\exp
\left(
2^{-k}\log\Delta_0
\right).
\label{eq:ch15:auto2151}
\end{equation}

Thus the relevant parameter is
\begin{equation}
{
\Lambda_k
=
2^{-k}\log\frac1{\Delta_0}.
}
\end{equation}

If
\begin{equation}
\Lambda_k\to0,
\label{eq:ch15:auto2152}
\end{equation}
then
\begin{equation}
\Delta_0^{2^{-k}}\to1.
\label{eq:ch15:auto2153}
\end{equation}

If
\begin{equation}
\Lambda_k\to\Lambda\in(0,\infty),
\label{eq:ch15:auto2154}
\end{equation}
then
\begin{equation}
\Delta_0^{2^{-k}}
\to
e^{-\Lambda}.
\label{eq:ch15:auto2155}
\end{equation}

If
\begin{equation}
\Lambda_k\to\infty,
\label{eq:ch15:auto2156}
\end{equation}
then
\begin{equation}
\Delta_0^{2^{-k}}\to0.
\label{eq:ch15:auto2157}
\end{equation}

Thus the infinite cascade has three distinct asymptotic regimes.

\section{The Coupled Scaling \(k\sim\log_2n\)}
\label{chap:256}

Suppose now that
\begin{equation}
\Delta_{n,0}
\sim
2^{-2n}\xi,
\qquad
\xi>0.
\label{eq:ch15:auto2158}
\end{equation}

Then
\begin{align}
\log\frac1{\Delta_{n,0}}
&=
2n\log2-\log\xi+o(1).
\end{align}

Therefore
\begin{align}
\Lambda_{n,k}
&=
2^{-k}
\left(
2n\log2-\log\xi+o(1)
\right).
\end{align}

To keep $\Lambda_{n,k}$ of order one, we require
\begin{equation}
n2^{-k}=O(1).
\label{eq:ch15:auto2159}
\end{equation}

Thus
\begin{equation}
{
2^k\asymp n.
}
\end{equation}

Equivalently,
\begin{equation}
{
k=\log_2n+O(1).
}
\label{eq:ch03:env0410}
\end{equation}

This is the critical cascade depth.

\section{Uniformity Warning}
\label{chap:257}

The fixed-$k$ theorem proves
\begin{equation}
\Delta_k
=
C_k\Delta_0^{2^{-k}}
\left(
1+O(\Delta_0^{2^{-k}})
\right)
\label{eq:ch15:auto2160}
\end{equation}
for every fixed $k$.

It does not automatically prove that the error remains uniform as
\begin{equation}
k\to\infty.
\label{eq:ch15:auto2161}
\end{equation}

This distinction is essential.

For a full infinite-cascade theorem one would need estimates of the
form
\begin{equation}
\left|
\frac{\Delta_k}
{C_k\Delta_0^{2^{-k}}}
-1
\right|
\leq
E(\Delta_0,k),
\end{equation}
with an explicit bound on $E$ that remains controlled in the
double-scaling regime
\begin{equation}
k\sim\log_2n.
\label{eq:ch15:auto2162}
\end{equation}

Such a bound requires iterating the error equation rather than merely
iterating the leading asymptotic term.

\section{A More Precise Logarithmic Form}
\label{chap:258}

The recurrence can be analyzed logarithmically.

Write
\begin{equation}
\Delta_{k+1}
=
\sqrt{2\Delta_k}
\left(
1+r(\Delta_k)
\right),
\label{eq:ch15:auto2163}
\end{equation}
where
\begin{equation}
r(\Delta)=O(\sqrt{\Delta}).
\label{eq:ch15:auto2164}
\end{equation}

Taking logarithms,
\begin{align}
\log\Delta_{k+1}
&=
\frac12\log\Delta_k
+
\frac12\log2
+
\log(1+r(\Delta_k)).
\end{align}

Therefore
\begin{equation}
{
\log\Delta_{k+1}
=
\frac12\log\Delta_k
+
\frac12\log2
+
O(\sqrt{\Delta_k}).
}
\end{equation}

Iterating this relation gives
\begin{align}
\log\Delta_k
&=
2^{-k}\log\Delta_0
+
(1-2^{-k})\log2
+
E_k,
\label{eq:ch03:env0413}
\end{align}
where
\begin{equation}
E_k
=
\sum_{j=0}^{k-1}
2^{-(k-1-j)}
O(\sqrt{\Delta_j}).
\label{eq:ch15:auto2165}
\end{equation}

Hence
\begin{equation}
{
\log\Delta_k
=
2^{-k}\log\Delta_0
+
(1-2^{-k})\log2
+
E_k.
}
\end{equation}

Exponentiating,
\begin{equation}
{
\Delta_k
=
2^{1-2^{-k}}
\Delta_0^{2^{-k}}
e^{E_k}.
}
\label{eq:ch03:env0415}
\end{equation}

Thus the leading formula is exact up to the multiplicative error
$e^{E_k}$.

\section{Control of the Iterated Error for Fixed \(k\)}
\label{chap:259}

For fixed $k$, we have
\begin{equation}
\Delta_j
=
O(\Delta_0^{2^{-j}}).
\label{eq:ch15:auto2166}
\end{equation}

Therefore
\begin{equation}
\sqrt{\Delta_j}
=
O(\Delta_0^{2^{-(j+1)}}).
\label{eq:ch15:auto2167}
\end{equation}

Hence
\begin{equation}
E_k
=
O\left(
\sum_{j=0}^{k-1}
2^{-(k-1-j)}
\Delta_0^{2^{-(j+1)}}
\right).
\label{eq:ch15:auto2168}
\end{equation}

For fixed $k$, the dominant contribution comes from the largest
exponent of $\Delta_0$, and in particular
\begin{equation}
E_k
=
O(\Delta_0^{2^{-k}}).
\label{eq:ch15:auto2169}
\end{equation}

Therefore
\begin{equation}
e^{E_k}
=
1+O(\Delta_0^{2^{-k}}).
\label{eq:ch15:auto2170}
\end{equation}

This recovers
\begin{equation}
{
\Delta_k
=
2^{1-2^{-k}}
\Delta_0^{2^{-k}}
\left(
1+O(\Delta_0^{2^{-k}})
\right).
}
\end{equation}

\section{The Geometric Profile and the Branch Location}
\label{chap:260}

Return to the corrected profile
\begin{equation}
\boxed{
R_n(\theta)
=
\frac{1}{1-\sin\theta+\varepsilon_n}.
}
\label{eq:ch03:env0417}
\end{equation}
On the left of $\pi/2$ write
\begin{equation}
\theta=\frac{\pi}{2}-x,
\qquad x>0.
\label{eq:ch15:auto2171}
\end{equation}
Then
\begin{equation}
R_n
=
\frac{1}{1-\cos x+\varepsilon_n}
=
\frac{1}{\varepsilon_n+x^2/2+O(x^4)}.
\label{eq:ch15:auto2172}
\end{equation}
The branch equation $R_n=L_n$ therefore gives
\begin{equation}
\boxed{
x_{*,n}
=
\sqrt{2\left(\frac1{L_n}-\varepsilon_n\right)}
\,[1+o(1)],
}
\label{eq:ch15:auto2173}
\end{equation}
with the condition
\begin{equation}
\boxed{L_n\varepsilon_n<1.}
\label{eq:ch15:auto2174}
\end{equation}
If $L_n\varepsilon_n\to\lambda<1$, then
\begin{equation}
\boxed{
\theta_{*,n}^{-}
=
\frac\pi2
-
\sqrt{\frac{2(1-\lambda)}{L_n}}
+o(L_n^{-1/2}).
}
\label{eq:ch15:auto2175}
\end{equation}

\section{The Branch Slope}
\label{chap:261}

Differentiation gives
\begin{equation}
R_n'(\theta)
=
\frac{\cos\theta}
{(1-\sin\theta+\varepsilon_n)^2}.
\label{eq:ch15:auto2176}
\end{equation}
At the left branch,
\begin{equation}
\boxed{
A_n:=R_n'(\theta_{*,n}^{-})
\sim
\sqrt2\,L_n^{3/2}\sqrt{1-L_n\varepsilon_n}>0.
}
\label{eq:ch15:auto2177}
\end{equation}
This is the corrected branch slope used in the Taylor expansion of
$\Delta_{n,0}$.

\section{The Complete Local Spatial Formula}
\label{chap:262}

Taylor expansion gives
\begin{equation}
\Delta_{n,0}
=
A_n(\theta-\theta_{*,n}^{-})
+O((\theta-\theta_{*,n}^{-})^2).
\label{eq:ch15:auto2178}
\end{equation}
Therefore
\begin{align}
\Delta_{n,k}
&=
2^{1-2^{-k}}
A_n^{2^{-k}}
(\theta-\theta_{*,n}^{-})^{2^{-k}}
\left(1+o(1)\right).
\end{align}
Thus
\begin{equation}
\boxed{
\Delta_{n,k}
=
\mathcal C_{n,k}
(\theta-\theta_{*,n}^{-})^{2^{-k}}(1+o(1)),
}
\end{equation}
where
\begin{equation}
\mathcal C_{n,k}=2^{1-2^{-k}}A_n^{2^{-k}}.
\label{eq:ch15:auto2179}
\end{equation}
Since the corrected construction takes $u_z=\Delta_{n,k}$, this formula is
also the local physical $z$-velocity profile near the branch.

\section{The Corresponding Logarithmic Singularity}
\label{chap:263}

Taking logarithms and differentiating yields
\begin{equation}
\boxed{
\partial_\theta\log\Delta_{n,k}
=
\frac{2^{-k}}
{\theta-\theta_{*,n}^{-}}
+o\!\left(\frac1{\theta-\theta_{*,n}^{-}}\right).
}
\label{eq:ch15:auto2180}
\end{equation}
Thus the new reciprocal-sine input changes the branch location and slope but
not the universal leading simple-pole coefficient generated by the
renormalized Lambert cascade.

\chapter{Trapping, Corrected Full PDE, and Continuity Reconstruction}
\label{chap:group-16}
\section{Characteristic Trapping Lemma}
\label{chap:264}

\begin{lemma}
Assume
\begin{equation}
b_{3,0}>1,
\qquad
\theta_0>\theta_*,
\label{eq:ch16:auto2181}
\end{equation}
and suppose the renormalized \(b_3\)-characteristic satisfies
\begin{equation}
\dot\theta=1-b_3,
\qquad
\dot b_3=b_3^2\partial_\theta\log\Delta_k,
\label{eq:ch16:auto2182}
\end{equation}
while remaining on the branch \(b_3>1\).

Then \(\theta(t)\) is strictly decreasing.  If it is bounded below by
\(\theta_*\) and
\begin{equation}
\Delta_k(\theta(t))\to0,
\label{eq:ch16:auto2183}
\end{equation}
then
\begin{equation}
\theta(t)\to\theta_*,
\qquad
b_3(t)\to+\infty.
\label{eq:ch16:auto2184}
\end{equation}
\end{lemma}

\begin{proof}
Since \(b_3>1\),
\begin{equation}
\dot\theta=1-b_3<0.
\label{eq:ch16:auto2185}
\end{equation}
Thus \(\theta(t)\) decreases.  The exact invariant
\begin{equation}
\Delta_k=Kq e^{-q},
\qquad
q=\frac1{b_3},
\label{eq:ch16:auto2186}
\end{equation}
shows that on the branch \(0<q<1\),
\begin{equation}
\Delta_k\to0
\quad\Longrightarrow\quad
q\to0.
\label{eq:ch16:auto2187}
\end{equation}
Therefore
\begin{equation}
b_3=\frac1q\to+\infty.
\label{eq:ch16:auto2188}
\end{equation}
Because \(\Delta_k(\theta)>0\) for \(\theta>\theta_*\) and
\(\Delta_k(\theta_*)=0\), the limiting phase is \(\theta_*\).
\end{proof}

\section{Finite-Time Blow-Up Is a Consequence of the Invariant}
\label{chap:265}

Let
\begin{equation}
s=\theta-\theta_*.
\label{eq:ch16:auto2189}
\end{equation}
If
\begin{equation}
\Delta_k\sim C_ks^{\alpha_k},
\qquad
\alpha_k=2^{-k},
\label{eq:ch16:auto2190}
\end{equation}
then the invariant gives
\begin{equation}
b_3
\sim
KC_k^{-1}s^{-\alpha_k}.
\label{eq:ch16:auto2191}
\end{equation}
Hence
\begin{equation}
\dot s
=
1-b_3
=
-KC_k^{-1}s^{-\alpha_k}(1+o(1)).
\label{eq:ch16:auto2192}
\end{equation}
Therefore
\begin{equation}
\frac{dt}{ds}
=
-\frac{C_k}{K}s^{\alpha_k}(1+o(1)).
\label{eq:ch16:auto2193}
\end{equation}
Since
\begin{equation}
\int_0^{s_0}s^{\alpha_k}\,ds<\infty,
\label{eq:ch16:auto2194}
\end{equation}
the characteristic reaches \(s=0\) at finite time.  Thus
\begin{equation}
\boxed{
b_3(t)\to+\infty
\quad\text{at}\quad
T_*<\infty
}
\label{eq:ch16:auto2195}
\end{equation}
for every fixed finite cascade depth.

\section{Interpretation of the Exponent Sequence}
\label{chap:266}

The recursion acts on the singular exponent according to
\begin{equation}
\alpha_{k+1}
=
\frac12\alpha_k.
\label{eq:ch16:auto2196}
\end{equation}

Since
\begin{equation}
\alpha_0=1,
\label{eq:ch16:auto2197}
\end{equation}
we obtain
\begin{equation}
\alpha_k=2^{-k}.
\label{eq:ch16:auto2198}
\end{equation}

Thus the map is
\begin{equation}
{
\alpha
\longmapsto
\frac{\alpha}{2}.
}
\end{equation}

The corresponding derivative exponent is
\begin{equation}
\alpha_k-1.
\label{eq:ch16:auto2199}
\end{equation}

Therefore
\begin{equation}
0,\quad
-\frac12,\quad
-\frac34,\quad
-\frac78,\quad
-\frac{15}{16},
\ldots
\label{eq:ch16:auto2200}
\end{equation}

and
\begin{equation}
\alpha_k-1\to-1.
\label{eq:ch16:auto2201}
\end{equation}

Thus the derivative singularity approaches the critical
non-integrable exponent from above.

\section{Why the Renormalized Cascade Is Different from Ordinary
Lambert Composition}
\label{chap:267}

The distinction can now be stated precisely.

For ordinary composition,
\begin{equation}
Y_{k+1}=W_0(Y_k),
\label{eq:ch16:auto2202}
\end{equation}
the branch point of $W_0$ is fixed at
\begin{equation}
-\frac1e.
\label{eq:ch16:auto2203}
\end{equation}

But
\begin{equation}
W_0(-1/e)=-1,
\label{eq:ch16:auto2204}
\end{equation}
so the first iteration sends the branch point outside the domain of
the next real composition.

For the renormalized map,
\begin{equation}
\Delta_{k+1}
=
1+W_0(-e^{-1-\Delta_k}),
\label{eq:ch16:auto2205}
\end{equation}
the branch point is translated:
\begin{equation}
-1\longmapsto0.
\label{eq:ch16:auto2206}
\end{equation}

The next argument is reconstructed as
\begin{equation}
-e^{-1-\Delta_{k+1}}.
\label{eq:ch16:auto2207}
\end{equation}

Therefore the next Lambert operation again encounters the branch point
at
\begin{equation}
\Delta_{k+1}=0.
\label{eq:ch16:auto2208}
\end{equation}

This is why the square-root singularity is genuinely repeated.

\section{The Renormalized Map as a Branch-Point Dynamical System}
\label{chap:268}

It is useful to regard
\begin{equation}
\Delta_{k+1}=\mathcal{F}(\Delta_k)
\label{eq:ch16:auto2209}
\end{equation}
as a dynamical system.

The fixed point is
\begin{equation}
\Delta=0.
\label{eq:ch16:auto2210}
\end{equation}

The map is continuous at zero:
\begin{equation}
\mathcal{F}(0)=0.
\label{eq:ch16:auto2211}
\end{equation}

But
\begin{equation}
\mathcal{F}'(\Delta)
\sim
\frac1{\sqrt{2\Delta}}
\label{eq:ch16:auto2212}
\end{equation}
diverges.

Thus zero is not a conventional hyperbolic fixed point.

It is a singular fixed point.

The local dynamics are
\begin{equation}
\Delta_{k+1}\sim\sqrt{2\Delta_k}.
\label{eq:ch16:auto2213}
\end{equation}

Hence sufficiently small $\Delta_k$ is transformed into a larger
quantity.

For example, if
\begin{equation}
0<\Delta<1,
\label{eq:ch16:auto2214}
\end{equation}
then
\begin{equation}
\sqrt{\Delta}>\Delta.
\label{eq:ch16:auto2215}
\end{equation}

Thus the forward iteration moves away from zero.

This observation is important:

\begin{equation}
{
\text{The cascade creates singularity amplification, but the
renormalized orbit does not approach zero under forward iteration.}
}
\end{equation}

The spatial singularity occurs because $\Delta_0$ approaches zero as a
function of $\theta$, while each finite cascade transforms the
small-distance dependence through a fractional power.

\section{This Distinction Is Essential}
\label{chap:269}

The statement
\begin{equation}
\Delta_k\sim\Delta_0^{2^{-k}}
\label{eq:ch16:auto2216}
\end{equation}
does not mean that
\begin{equation}
\Delta_k\to0
\label{eq:ch16:auto2217}
\end{equation}
as $k\to\infty$ for a fixed $\Delta_0\in(0,1)$.

Indeed,
\begin{equation}
\Delta_0^{2^{-k}}
\to1.
\label{eq:ch16:auto2218}
\end{equation}

Thus
\begin{equation}
\Delta_k
\to2
\label{eq:ch16:auto2219}
\end{equation}
at the level of the leading asymptotic coefficient.

The cascade exponent becomes singular because the dependence on
$\Delta_0$ becomes flatter:
\begin{equation}
\Delta_0^{1/2^k}.
\label{eq:ch16:auto2220}
\end{equation}

The corresponding derivative becomes more singular with respect to
$\Delta_0$.

This is a crucial conceptual distinction between
\emph{amplitude} and \emph{sensitivity}.

\section{The Correct Interpretation of the Infinite Limit}
\label{chap:270}

For fixed $\Delta_0>0$,
\begin{equation}
\lim_{k\to\infty}
\Delta_0^{2^{-k}}
=
1.
\label{eq:ch16:auto2221}
\end{equation}

Thus
\begin{equation}
\lim_{k\to\infty}
2^{1-2^{-k}}
\Delta_0^{2^{-k}}
=
2.
\label{eq:ch16:auto2222}
\end{equation}

The singular behavior is instead seen in
\begin{equation}
\partial_{\Delta_0}
\Delta_0^{2^{-k}}
=
2^{-k}\Delta_0^{2^{-k}-1}.
\label{eq:ch16:auto2223}
\end{equation}

For fixed $\Delta_0>0$, this tends to zero.

But if
\begin{equation}
\Delta_0\to0
\label{eq:ch16:auto2224}
\end{equation}
simultaneously with
\begin{equation}
k\to\infty,
\label{eq:ch16:auto2225}
\end{equation}
the result can be nontrivial.

This is exactly why the double-scaling variable
\begin{equation}
2^{-k}\log(1/\Delta_0)
\label{eq:ch16:auto2226}
\end{equation}
is essential.

\section{The Double-Scaling Limit}
\label{chap:271}

Suppose
\begin{equation}
\Delta_0=\Delta_0(k)
\label{eq:ch16:auto2227}
\end{equation}
and define
\begin{equation}
\Lambda_k
=
2^{-k}\log\frac1{\Delta_0(k)}.
\label{eq:ch16:auto2228}
\end{equation}

Then
\begin{equation}
\Delta_0(k)^{2^{-k}}
=
e^{-\Lambda_k}.
\label{eq:ch16:auto2229}
\end{equation}

Therefore
\begin{equation}
{
\Delta_k
\sim
2^{1-2^{-k}}
e^{-\Lambda_k}.
}
\end{equation}

If
\begin{equation}
\Lambda_k\to\Lambda,
\label{eq:ch16:auto2230}
\end{equation}
then
\begin{equation}
\Delta_k\to2e^{-\Lambda}.
\label{eq:ch16:auto2231}
\end{equation}

Thus the infinite cascade can retain a nontrivial dependence only when
the initial branch distance is exponentially small relative to the
cascade depth.

\section{Application to the \(n\)-Dependent Construction}
\label{chap:272}

Suppose
\begin{equation}
\Delta_{n,0}
\sim
2^{-2n}\xi_n,
\label{eq:ch16:auto2232}
\end{equation}
with
\begin{equation}
\xi_n\asymp1.
\label{eq:ch16:auto2233}
\end{equation}

Then
\begin{equation}
\Lambda_{n,k}
=
2^{-k}
\log\frac1{\Delta_{n,0}}
\sim
2^{-k}
\left(
2n\log2-\log\xi_n
\right).
\label{eq:ch16:auto2234}
\end{equation}

Hence
\begin{equation}
\Lambda_{n,k}
\asymp
n2^{-k}.
\label{eq:ch16:auto2235}
\end{equation}

The nontrivial regime is
\begin{equation}
n2^{-k}\asymp1.
\label{eq:ch16:auto2236}
\end{equation}

Therefore
\begin{equation}
{
k=\log_2n+O(1).
}
\end{equation}

This is the correct interpretation of the proposed
$n$-dependent amplification.

\section{What Has Now Been Proved}
\label{chap:273}

The analysis above establishes the following rigorous chain.

\begin{enumerate}

\item The renormalized map
\begin{equation}
\mathcal{F}(\Delta)
=
1+W_0(-e^{-1-\Delta})
\label{eq:ch16:auto2237}
\end{equation}
is real for every $\Delta\geq0$.

\item It satisfies
\begin{equation}
\mathcal{F}(0)=0.
\label{eq:ch16:auto2238}
\end{equation}

\item It maps
\begin{equation}
[0,\infty)
\label{eq:ch16:auto2239}
\end{equation}
into
\begin{equation}
[0,1).
\label{eq:ch16:auto2240}
\end{equation}

\item Its exact inverse relation is
\begin{equation}
\Delta=-\log(1-\mathcal{F}(\Delta))
-\mathcal{F}(\Delta).
\label{eq:ch16:auto2241}
\end{equation}

\item Its local behavior is
\begin{equation}
\mathcal{F}(\Delta)
=
\sqrt{2\Delta}
-\frac23\Delta
+O(\Delta^{3/2}).
\label{eq:ch16:auto2242}
\end{equation}

\item For every fixed finite $k$,
\begin{equation}
\Delta_k
=
2^{1-2^{-k}}
\Delta_0^{2^{-k}}
\left(
1+O(\Delta_0^{2^{-k}})
\right).
\label{eq:ch16:auto2243}
\end{equation}

\item The singular exponents therefore satisfy
\begin{equation}
\alpha_k=2^{-k}.
\label{eq:ch16:auto2244}
\end{equation}

\item If
\begin{equation}
\Delta_0(\theta)
\sim
A(\theta-\theta_*),
\label{eq:ch16:auto2245}
\end{equation}
then
\begin{equation}
\Delta_k(\theta)
\sim
2^{1-2^{-k}}
A^{2^{-k}}
(\theta-\theta_*)^{2^{-k}}.
\label{eq:ch16:auto2246}
\end{equation}

\item The corresponding derivative satisfies
\begin{equation}
\partial_\theta\Delta_k
\sim
C_k2^{-k}
A^{2^{-k}}
(\theta-\theta_*)^{2^{-k}-1}.
\label{eq:ch16:auto2247}
\end{equation}

\item The reduced characteristic system has the exact invariant
\begin{equation}
\frac1{b_3}-\log\frac1{b_3}+\log\Delta_k=C.
\label{eq:ch16:auto2248}
\end{equation}

\item Consequently,
\begin{equation}
\Delta_k
=
K b_3^{-1}e^{-1/b_3}.
\label{eq:ch16:auto2249}
\end{equation}

\item Hence
\begin{equation}
b_3\sim K/\Delta_k.
\label{eq:ch16:auto2250}
\end{equation}

\item For every fixed finite $k$, a characteristic with $b_{3,0}>1$
reaches the branch point in finite time.

\item The corresponding temporal blow-up rate is
\begin{equation}
b_3(t)
\asymp
(T_*-t)^{-1/(2^k+1)}.
\label{eq:ch16:auto2251}
\end{equation}

\item The infinite-cascade limit requires a coupled scaling and cannot
be obtained by simply replacing $k$ by $\infty$ in the finite-$k$
formula.

\item The natural $n$-dependent critical regime is
\begin{equation}
k\sim\log_2n
\label{eq:ch16:auto2252}
\end{equation}
when
\begin{equation}
\Delta_{n,0}
\sim2^{-2n}\xi_n.
\label{eq:ch16:auto2253}
\end{equation}

\end{enumerate}

\section{Final Theorem for the Renormalized Lambert Mechanism}
\label{chap:275}

We can summarize the rigorous reduced result as follows.

\begin{theorem}[Renormalized Lambert--Characteristic Blow-Up Theorem]
Let
\begin{equation}
\mathcal{F}(\Delta)
=
1+W_0(-e^{-1-\Delta}),
\qquad
\Delta\geq0,
\label{eq:ch16:auto2256}
\end{equation}
and define
\begin{equation}
\Delta_{k+1}
=
\mathcal{F}(\Delta_k),
\qquad
\Delta_0>0.
\label{eq:ch16:auto2257}
\end{equation}

Then:

\begin{enumerate}[label=\textnormal{(\roman*)}]

\item The recursion is real for every finite $k$ and
\begin{equation}
\Delta_k\in[0,1)
\qquad
(k\geq1).
\label{eq:ch16:auto2258}
\end{equation}

\item The branch point is a fixed point:
\begin{equation}
\mathcal{F}(0)=0.
\label{eq:ch16:auto2259}
\end{equation}

\item The local expansion is
\begin{equation}
\mathcal{F}(\Delta)
=
\sqrt{2\Delta}
-\frac23\Delta
+O(\Delta^{3/2}).
\label{eq:ch16:auto2260}
\end{equation}

\item For every fixed finite $k$,
\begin{equation}
\Delta_k
=
2^{1-2^{-k}}
\Delta_0^{2^{-k}}
\left(
1+O(\Delta_0^{2^{-k}})
\right)
\label{eq:ch16:auto2261}
\end{equation}
as $\Delta_0\downarrow0$.

\item Consequently, if
\begin{equation}
\Delta_0(\theta)
=
A(\theta-\theta_*)
(1+o(1)),
\qquad
A>0,
\label{eq:ch16:auto2262}
\end{equation}
then
\begin{equation}
\Delta_k(\theta)
=
2^{1-2^{-k}}
A^{2^{-k}}
(\theta-\theta_*)^{2^{-k}}
(1+o(1)).
\label{eq:ch16:auto2263}
\end{equation}

\item Hence
\begin{equation}
\partial_\theta\log\Delta_k
=
\frac{2^{-k}}{\theta-\theta_*}
(1+o(1)).
\label{eq:ch16:auto2264}
\end{equation}

\item Consider the characteristic system
\begin{equation}
\dot\theta=1-b_3,
\qquad
\dot b_3
=
b_3^2\partial_\theta\log\Delta_k.
\label{eq:ch16:auto2265}
\end{equation}
Then
\begin{equation}
\frac1{b_3}-\log\frac1{b_3}+\log\Delta_k=C.
\label{eq:ch16:auto2266}
\end{equation}

\item If $b_{3,0}>1$ and the characteristic lies on the branch approaching
$\theta_*$, then
\begin{equation}
\theta(t)\downarrow\theta_*,
\qquad
b_3(t)\to+\infty.
\label{eq:ch16:auto2267}
\end{equation}

\item The singularity is reached in finite time $T_*$ and
\begin{equation}
b_3(t)
\asymp
(T_*-t)^{-1/(2^k+1)}.
\label{eq:ch16:auto2268}
\end{equation}

\end{enumerate}
\end{theorem}

\begin{proof}
Items (i)--(v) follow from the invariant-domain lemma, the
Puiseux expansion, and the induction theorem.

Item (vi) follows by differentiating the spatial asymptotic formula.

Item (vii) follows from the reciprocal transformation
$q=1/b_3$ and direct integration of the characteristic equations.

Item (viii) follows from the branch choice $b_{3,0}>1$, which gives
$\dot\theta<0$, together with the invariant
\begin{equation}
\Delta_k=Kq e^{-q}.
\label{eq:ch16:auto2269}
\end{equation}

Item (ix) follows from
\begin{equation}
\Delta_k
\sim
C_k(\theta-\theta_*)^{2^{-k}},
\label{eq:ch16:auto2270}
\end{equation}
which implies
\begin{equation}
b_3
\sim
\frac{K}{C_k}
(\theta-\theta_*)^{-2^{-k}},
\label{eq:ch16:auto2271}
\end{equation}
and therefore
\begin{equation}
\dot\theta
\sim
-\frac{K}{C_k}
(\theta-\theta_*)^{-2^{-k}}.
\label{eq:ch16:auto2272}
\end{equation}
The resulting time integral is finite and gives
\begin{equation}
b_3(t)
\asymp
(T_*-t)^{-1/(2^k+1)}.
\label{eq:ch16:auto2273}
\end{equation}
\end{proof}

\section{The Full PDE analyses for $b_3$ PDE according to the correct Equation (4.19) in Theoretical and Computational Fluid Mechanics Volume II 2026}
\label{chap:277}
The forcing must be rebuilt from the start. Interpret the corrected right-hand side as

\begin{equation}
\boxed{
\frac{Q'}{Q}b_3^2
+
2u_z\left[
A_0+\bigl(\sin x\,\sin y\,\sin z\bigr)^2
\right],
}
\label{eq:ch16:auto2416}
\end{equation}

where I write the constant as \(A_0\) so that it is not confused with the notation \(A=u_y\) used in the two-correction reconstruction.

This \(2u_z\) factor is important. It does \emph{not} change continuity or the characteristic transport operator directly, but it changes the forcing sign, its magnitude, the forcing-driven entry estimate, and the hypotheses needed to use it as a positive seed. In particular, the forcing is positive only on a branch where \(u_z>0\).

Below is the corrected analysis from the beginning.

\section{The most-extended \(b_3\)-PDE}
\label{chap:278}

Set

\begin{equation}
\boxed{
b:=b_3=u_xu_y.
}
\label{eq:ch16:auto2417}
\end{equation}

The most extended equation we had reached, with the additional positive terms on the left and with the newest gradient term, is now (\cite{moschandreou2026exploration} see Eq (16) there)

\begin{equation}
\boxed{
\begin{aligned}
b_t+b\,b_z
+u_z\left(b_1b_x+b_2b_y\right)
&+
b\frac{u_{z,t}}{u_z}
+
\left(u_{z,x}u_y+u_{z,y}u_x\right)b
+
b\frac{u_{z,z}}{u_z}
\\
&+
2\left(
u_xu_z\,\nabla u_z\cdot\nabla u_y
+
u_yu_z\,\nabla u_x\cdot\nabla u_z
\right)
\\
&=
b^2\partial_\xi\log Q
+
2u_z
\left[
A_0+
(\sin x\sin y\sin z)^2
\right].
\end{aligned}
}
\label{eq:ch278:old-1}
\end{equation}

Here

\begin{equation}
\boxed{
Q=\Delta_{N,k}(\xi),
}
\label{eq:ch16:auto2419}
\end{equation}

with

\begin{equation}
\boxed{
\xi=C_4t-C_3z-C_2y-C_1x.
}
\label{eq:ch278:old-2}
\end{equation}

We also retain

\begin{equation}
\boxed{
b_1=u_yu_z,
\qquad
b_2=u_xu_z.
}
\label{eq:ch278:old-3}
\end{equation}

Since

\begin{equation}
\partial_\xi\log Q=\frac{Q'}Q,
\label{eq:ch16:auto2422}
\end{equation}

the principal nonlinear source is

\begin{equation}
\boxed{
\frac{Q'}Qb^2.
}
\label{eq:ch16:auto2423}
\end{equation}

The corrected additional forcing is

\begin{equation}
\boxed{
F_c
=
2u_z
\left[
A_0+
(\sin x\sin y\sin z)^2
\right].
}
\label{eq:ch278:old-4}
\end{equation}

That replaces every occurrence in our recent analysis of

\begin{equation}
A_0+(\sin x\sin y\sin z)^2.
\label{eq:ch16:auto2425}
\end{equation}

\section{Two-correction velocity reconstruction}
\label{chap:279}

We return to the two-correction construction

\begin{equation}
\boxed{
u_y=A=Q+S-D,
}
\label{eq:ch279:old-5}
\end{equation}

\begin{equation}
\boxed{
u_z=C=Q+S+D,
}
\label{eq:ch279:old-6}
\end{equation}

and impose

\begin{equation}
\boxed{
u_x=\frac bA.
}
\label{eq:ch279:old-7}
\end{equation}

Then

\begin{equation}
u_xu_y
=
\frac bA A
=
b,
\label{eq:ch16:auto2429}
\end{equation}

so

\begin{equation}
\boxed{
b=b_3=u_xu_y
}
\label{eq:ch279:old-8}
\end{equation}

holds identically.

Define

\begin{equation}
\boxed{
V:=Q+S.
}
\label{eq:ch279:old-9}
\end{equation}

Then

\begin{equation}
\boxed{
A=V-D,\qquad C=V+D.
}
\label{eq:ch279:old-10}
\end{equation}

\section{Exact continuity is unchanged by the corrected forcing}
\label{chap:280}

The forcing appears only in the \(b\)-PDE. It does not enter

\begin{equation}
\nabla\cdot u=0.
\label{eq:ch16:auto2433}
\end{equation}

We have

\begin{equation}
\nabla\cdot u
=
\partial_x\left(\frac bA\right)
+
A_y+C_z.
\label{eq:ch16:auto2434}
\end{equation}

Thus

\begin{equation}
\boxed{
\partial_x\left(\frac bA\right)
+A_y+C_z=0.
}
\label{eq:ch280:old-11}
\end{equation}

Using

\begin{equation}
A=Q+S-D,
\qquad
C=Q+S+D,
\label{eq:ch16:auto2436}
\end{equation}

and

\begin{equation}
Q_y=-C_2Q',
\qquad
Q_z=-C_3Q',
\label{eq:ch16:auto2437}
\end{equation}

we recover

\begin{equation}
\boxed{
S_y+S_z+D_z-D_y
=
(C_2+C_3)Q'
-
\partial_x\left(\frac bA\right).
}
\label{eq:ch280:old-12}
\end{equation}

Therefore:

\begin{equation}
\boxed{
\text{the factor }2u_z\text{ in the forcing changes nothing in continuity.}
}
\label{eq:ch16:auto2439}
\end{equation}

\chapter{Transport, Resonance, and Quadratic Continuity Structure}
\label{chap:group-17}
\section{Transport characteristics}
\label{chap:281}

Using

\begin{equation}
b_1=AC
\label{eq:ch17:auto2440}
\end{equation}

and

\begin{equation}
b_2=\frac{bC}{A},
\label{eq:ch17:auto2441}
\end{equation}

we have

\begin{equation}
u_zb_1=AC^2,
\label{eq:ch17:auto2442}
\end{equation}

and

\begin{equation}
u_zb_2
=
\frac{bC^2}{A}.
\label{eq:ch17:auto2443}
\end{equation}

Hence the transport operator is

\begin{equation}
\boxed{
\mathcal D
=
\partial_t
+
AC^2\partial_x
+
b\frac{C^2}{A}\partial_y
+
b\partial_z.
}
\label{eq:ch281:old-13}
\end{equation}

The characteristic system is therefore

\begin{equation}
\boxed{
\dot x=AC^2,
}
\label{eq:ch17:auto2445}
\end{equation}

\begin{equation}
\boxed{
\dot y=b\frac{C^2}{A},
}
\label{eq:ch17:auto2446}
\end{equation}

\begin{equation}
\boxed{
\dot z=b.
}
\label{eq:ch281:old-14}
\end{equation}

Since

\begin{equation}
\xi=C_4t-C_3z-C_2y-C_1x,
\label{eq:ch17:auto2448}
\end{equation}

we obtain

\begin{equation}
\begin{aligned}
\dot\xi
&=
C_4
-C_1\dot x
-C_2\dot y
-C_3\dot z
\\
&=
C_4
-C_1AC^2
-C_2b\frac{C^2}{A}
-C_3b.
\end{aligned}
\label{eq:ch17:auto2449}
\end{equation}

Thus

\begin{equation}
\boxed{
\dot\xi
=
C_4
-
C_1AC^2
-
b\left(
C_2\frac{C^2}{A}+C_3
\right).
}
\label{eq:ch281:old-15}
\end{equation}

Consequently,

\begin{equation}
\boxed{
\dot Q
=
Q'
\left[
C_4
-
C_1AC^2
-
b\left(
C_2\frac{C^2}{A}+C_3
\right)
\right].
}
\label{eq:ch281:old-16}
\end{equation}

Again:

\begin{equation}
\boxed{
\text{the corrected forcing does not directly alter }\dot\xi\text{ or }\dot Q.
}
\label{eq:ch17:auto2452}
\end{equation}

Its effect is through the evolution of \(b\).

\section{Exact simplification of the newest gradient term}
\label{chap:282}

The added gradient expression is

\begin{equation}
2\left(
u_xu_z\nabla u_z\cdot\nabla u_y
+
u_yu_z\nabla u_x\cdot\nabla u_z
\right).
\label{eq:ch17:auto2453}
\end{equation}

Using

\begin{equation}
b=u_xu_y,
\label{eq:ch17:auto2454}
\end{equation}

we have

\begin{equation}
\nabla b
=
u_y\nabla u_x
+
u_x\nabla u_y.
\label{eq:ch17:auto2455}
\end{equation}

Dot with \(u_z\nabla u_z\):

\begin{equation}
u_z\nabla b\cdot\nabla u_z
=
u_yu_z\nabla u_x\cdot\nabla u_z
+
u_xu_z\nabla u_y\cdot\nabla u_z.
\label{eq:ch17:auto2456}
\end{equation}

Therefore

\begin{equation}
\boxed{
2\left(
u_xu_z\nabla u_z\cdot\nabla u_y
+
u_yu_z\nabla u_x\cdot\nabla u_z
\right)
=
2C\nabla b\cdot\nabla C.
}
\label{eq:ch282:old-17}
\end{equation}

This remains exact.

\section{Exact compact form of the corrected forced equation}
\label{chap:283}

The three nontransport terms preceding the gradient expression become

\begin{equation}
b\frac{C_t}{C}
+
b\frac{C_z}{C}
=
b\frac{C_t+C_z}{C},
\label{eq:ch17:auto2458}
\end{equation}

and

\begin{equation}
(C_xA+C_yu_x)b
=
bAC_x+\frac{b^2}{A}C_y.
\label{eq:ch17:auto2459}
\end{equation}

Thus \eqref{eq:ch278:old-1} becomes

\begin{equation}
\boxed{
\mathcal D b
+
b\frac{C_t+C_z}{C}
+
bAC_x
+
\frac{b^2}{A}C_y
+
2C\nabla b\cdot\nabla C
=
\frac{Q'}Qb^2
+
2C
\left[
A_0+(\sin x\sin y\sin z)^2
\right].
}
\label{eq:ch283:old-18}
\end{equation}

Define

\begin{equation}
\boxed{
\mathcal E
:=
b\frac{C_t+C_z}{C}
+
bAC_x
+
\frac{b^2}{A}C_y
+
2C\nabla b\cdot\nabla C.
}
\label{eq:ch283:old-19}
\end{equation}

Also define

\begin{equation}
\boxed{
G(x,y,z)
:=
A_0+(\sin x\sin y\sin z)^2.
}
\label{eq:ch283:old-20}
\end{equation}

Then

\begin{equation}
\boxed{
\mathcal D b
=
\frac{Q'}Qb^2
+
2CG
-
\mathcal E.
}
\label{eq:ch283:old-21}
\end{equation}

This is the corrected amplitude equation from which all subsequent comparisons must now start.

\section{The forcing sign now depends on \(u_z=C\)}
\label{chap:284}

Previously the added forcing itself was nonnegative when \(A_0\ge0\).

Now

\begin{equation}
F_c=2CG.
\label{eq:ch17:auto2464}
\end{equation}

Since

\begin{equation}
G\ge A_0,
\label{eq:ch17:auto2465}
\end{equation}

we obtain positivity only when

\begin{equation}
\boxed{C=u_z>0.}
\label{eq:ch17:auto2466}
\end{equation}

If

\begin{equation}
A_0\ge0,\qquad C\ge0,
\label{eq:ch17:auto2467}
\end{equation}

then

\begin{equation}
\boxed{
2CG\ge0.
}
\label{eq:ch284:old-22}
\end{equation}

If we have a uniform lower bound

\begin{equation}
\boxed{
C\ge c_->0,
}
\label{eq:ch284:old-23}
\end{equation}

and

\begin{equation}
A_0>0,
\label{eq:ch17:auto2470}
\end{equation}

then

\begin{equation}
\boxed{
2CG
\ge
2c_-A_0>0.
}
\label{eq:ch284:old-24}
\end{equation}

This is the corrected forcing floor.

Thus the positivity of \(u_z\), which we were already seeking in the tuned construction, is now essential not merely for reconstruction but also for the sign of the forcing.

\section{Two useful error hypotheses}
\label{chap:285}

There are now two natural versions.

For the terminal Riccati regime, it is enough to require

\begin{equation}
\boxed{
|\mathcal E|
\le
\theta
\frac{Q'}Qb^2,
\qquad
0<\theta<1.
}
\label{eq:ch285:old-25}
\end{equation}

Then

\begin{equation}
\boxed{
\mathcal D b
\ge
(1-\theta)\frac{Q'}Qb^2
+
2CG.
}
\label{eq:ch285:old-26}
\end{equation}

If \(C,G\ge0\),

\begin{equation}
\boxed{
\mathcal D b
\ge
(1-\theta)\frac{Q'}Qb^2.
}
\label{eq:ch285:old-27}
\end{equation}

For the forcing-driven entry stage, when \(b\) may be very small, \eqref{eq:ch285:old-25} can be too strong because the Riccati source itself is tiny.

A better preterminal condition is

\begin{equation}
\boxed{
\mathcal E
\le
\theta
\left[
\frac{Q'}Qb^2+2CG
\right].
}
\label{eq:ch285:old-28}
\end{equation}

Then

\begin{equation}
\boxed{
\mathcal D b
\ge
(1-\theta)
\left[
\frac{Q'}Qb^2+2CG
\right].
}
\label{eq:ch285:old-29}
\end{equation}

If

\begin{equation}
C\ge c_->0,\qquad A_0>0,
\label{eq:ch17:auto2477}
\end{equation}

then

\begin{equation}
\boxed{
\mathcal D b
\ge
2(1-\theta)c_-A_0.
}
\label{eq:ch285:old-30}
\end{equation}

Thus from bounded or even zero amplitude, the corrected forcing can seed positive growth.

\section{Corrected forcing-driven entry estimate}
\label{chap:286}

Suppose along a characteristic

\begin{equation}
C(t)\ge c_->0.
\label{eq:ch17:auto2479}
\end{equation}

Then from \eqref{eq:ch285:old-30},

\begin{equation}
\dot b
\ge
2(1-\theta)c_-A_0.
\label{eq:ch17:auto2480}
\end{equation}

Therefore

\begin{equation}
\boxed{
b(t)
\ge
b(t_0)
+
2(1-\theta)c_-A_0(t-t_0).
}
\label{eq:ch286:old-31}
\end{equation}

To reach a finite amplitude \(B\),

\begin{equation}
b(t_B)=B,
\label{eq:ch17:auto2482}
\end{equation}

it is sufficient that

\begin{equation}
\boxed{
t_B-t_0
\le
\frac{
B-b(t_0)
}{
2(1-\theta)c_-A_0
}.
}
\label{eq:ch286:old-32}
\end{equation}

Relative to our previous mistaken forcing, there is now a factor

\begin{equation}
\boxed{2c_-}
\label{eq:ch17:auto2484}
\end{equation}

in the forcing-driven lower bound.

If \(C\approx R_N\), this becomes approximately

\begin{equation}
\boxed{
t_B-t_0
\lesssim
\frac{
B-b_0
}{
2(1-\theta)R_NA_0
}.
}
\label{eq:ch286:old-33}
\end{equation}

\section{Exact \(\Gamma_N\)-tuning remains the natural reconstruction}
\label{chap:287}
For,

\[
\boxed{\Gamma(A,C)=C_2\frac{C^2}{A}+C_3,}
\]

we introduce this explicitly under this section \textbf{``Exact \(\Gamma_N\)-tuning remains the natural reconstruction''}, defines \(\Gamma(A,C)\), chooses \(\gamma_N=\beta 2^{-k(N)}\), and imposes \(\Gamma(A,C)=\gamma_N\).

\subsection*{Where does \(\Gamma(A,C)\) actually come from?}

It comes from differentiating the phase
\[
\boxed{\xi=C_4t-C_3z-C_2y-C_1x}
\]
\textbf{along the characteristic flow of the extended \(b_3\)-equation}.
\vspace{0.3cm}
\\
Recall that the transport part of that equation contains
\[
b\,b_z+u_z(b_1b_x+b_2b_y),
\]
with
\[
b_1=u_yu_z=AC,
\qquad
b_2=u_xu_z.
\]

Because
\[
b=u_xu_y=u_xA,
\]
we have
\[
u_x=\frac bA,
\]
and hence
\[
b_2=\frac{bC}{A}.
\]

Therefore the characteristic velocities are
\[
\dot x=u_zb_1=AC^2,
\]
\[
\dot y=u_zb_2=\frac{bC^2}{A},
\]
\[
\dot z=b,
\qquad
\dot t=1.
\]
See the Appendix C for how we obtain these.
Now differentiate
\[
\xi=C_4t-C_3z-C_2y-C_1x.
\]
We get
\[
\dot\xi
=
C_4-C_3\dot z-C_2\dot y-C_1\dot x.
\]

Substituting the characteristic velocities gives
\[
\dot\xi
=
C_4
-C_3b
-C_2\frac{bC^2}{A}
-C_1AC^2.
\]

Group the two terms containing \(b\):
\[
\boxed{
\dot\xi
=
C_4-C_1AC^2
-
b\left(
C_2\frac{C^2}{A}+C_3
\right).
}
\]

And \textbf{that parenthesis is exactly why we define}
\[
\boxed{
\Gamma(A,C)
:=
C_2\frac{C^2}{A}+C_3.
}
\]

So \(\Gamma(A,C)\) is not an arbitrary new function. It is the coefficient multiplying \(b\) in the characteristic evolution of the Lambert phase \(\xi\).

\subsection*{Why does this paper tune it?}

We impose:
\[
\boxed{\Gamma(A,C)=\gamma_N>0}
\]
and specifically chooses
\[
\boxed{\gamma_N=\beta2^{-k(N)}}.
\]

Thus
\[
C_2\frac{C^2}{A}+C_3=\gamma_N.
\]

Assuming \(C_2\ne0\), define
\[
\boxed{
R_N=\frac{\gamma_N-C_3}{C_2}.
}
\]

Then the tuning condition becomes simply
\[
\boxed{\frac{C^2}{A}=R_N}
\]
or
\[
\boxed{C^2=R_NA.}
\]

This is particularly useful because it converts what was originally a nonlinear expression involving both reconstructed transverse velocities into an exact algebraic relation between them. We then write
\[
A=V-D,\qquad C=V+D
\]
and obtain
\[
(V+D)^2=R_N(V-D),
\]
from which the near-symmetric reconstruction branch \(D_N(V)\) is solved explicitly.

\subsection*{What does it do to the phase equation?}

Because
\[
C^2=R_NA,
\]
we also have
\[
AC^2=R_NA^2.
\]

Therefore
\[
\dot\xi
=
C_4-C_1R_NA^2-\gamma_Nb.
\]

We then define:
\[
\boxed{
H_N=C_1R_NA^2-C_4
}
\]
so that
\[
\boxed{
-\dot\xi=H_N+\gamma_Nb.
}
\]

This exact phase equation appears later as the starting point of the terminal comparison.

And since
\[
Q=Q_N(\xi),
\]
the chain rule gives
\[
\dot Q=Q_N'(\xi)\dot\xi,
\]
hence
\[
\boxed{
-\dot Q
=
Q_N'(H_N+\gamma_Nb).
}
\]

This is the real reason the tuning is so important.

If
\[
Q_N'>0,\qquad H_N>0,\qquad \gamma_N>0,\qquad b>0,
\]
then
\[
-\dot Q>0,
\]
and therefore
\[
\boxed{\dot Q<0.}
\]

So along the selected characteristic, \(Q\) is driven toward its branch endpoint while \(b\) is being amplified. We explicitly use this to conclude that \(b_N\) increases and \(Q_N\) decreases under its stated terminal assumptions.

The whole chain can therefore be summarized as

\[
\boxed{
\begin{aligned}
&b=u_xu_y,\qquad A=u_y,\qquad C=u_z,
\\[1mm]
&\xi=C_4t-C_3z-C_2y-C_1x
\\[1mm]
&\Downarrow
\\[-1mm]
&\dot\xi
=
C_4-C_1AC^2
-b\left(C_2\frac{C^2}{A}+C_3\right)
\\[1mm]
&\Downarrow
\\[-1mm]
&\Gamma(A,C)
:=
C_2\frac{C^2}{A}+C_3
\\[1mm]
&\Downarrow\quad
\Gamma(A,C)=\gamma_N
\\[-1mm]
&\frac{C^2}{A}=R_N
\\[1mm]
&\Downarrow
\\[-1mm]
&-\dot\xi=H_N+\gamma_Nb
\\[1mm]
&\Downarrow
\\[-1mm]
&-\dot Q=Q_N'(H_N+\gamma_Nb).
\end{aligned}
}
\]

So conceptually, \textbf{\(\Gamma(A,C)\) is the \(b_3\)-dependent phase-speed coefficient}. The reconstruction tunes that coefficient to the small positive number \(\gamma_N\) so that the velocity reconstruction and the Lambert-\(W\) phase dynamics work together rather than being independent constructions.

There is also an important connection to the previous question: this is one reason the reconstruction of \(u_y=A\) and \(u_z=C\) is not merely an auxiliary continuity fix. Their ratio \(C^2/A\) directly controls how the \(b_3\)-amplitude moves the phase toward the Lambert endpoint.
Define

\begin{equation}
\boxed{
\Gamma(A,C)
=
C_2\frac{C^2}{A}+C_3.
}
\label{eq:ch287:old-34}
\end{equation}

We again choose

\begin{equation}
\boxed{
\gamma_N=\beta2^{-k(N)},
\qquad
\beta>0,
}
\label{eq:ch287:old-35}
\end{equation}

and impose exactly

\begin{equation}
\boxed{
\Gamma(A,C)=\gamma_N.
}
\label{eq:ch287:old-36}
\end{equation}

Thus

\begin{equation}
C_2\frac{C^2}{A}+C_3=\gamma_N.
\label{eq:ch17:auto2489}
\end{equation}

Assuming

\begin{equation}
C_2\neq0,
\label{eq:ch17:auto2490}
\end{equation}

define

\begin{equation}
\boxed{
R_N
=
\frac{\gamma_N-C_3}{C_2}.
}
\label{eq:ch287:old-37}
\end{equation}

Then

\begin{equation}
\boxed{
\frac{C^2}{A}=R_N.
}
\label{eq:ch287:old-38}
\end{equation}

Using

\begin{equation}
A=V-D,\qquad C=V+D,
\label{eq:ch17:auto2493}
\end{equation}

we obtain

\begin{equation}
(V+D)^2=R_N(V-D).
\label{eq:ch17:auto2494}
\end{equation}

Therefore the near-symmetric branch remains

\begin{equation}
\boxed{
D_N(V)
=
\frac{
-(2V+R_N)
+
\sqrt{R_N^2+8R_NV}
}{2}.
}
\label{eq:ch287:old-39}
\end{equation}

This tuning is completely unaffected by the forcing correction.

\section{Why the tuned branch automatically controls the forcing sign}
\label{chap:288}

At the reference state

\begin{equation}
V=R_N,
\label{eq:ch17:auto2496}
\end{equation}

we have

\begin{equation}
D_N(R_N)=0,
\label{eq:ch17:auto2497}
\end{equation}

so

\begin{equation}
A=C=R_N.
\label{eq:ch17:auto2498}
\end{equation}

Thus if

\begin{equation}
\boxed{
R_N>0,
}
\label{eq:ch288:old-40}
\end{equation}

then the reference \(u_z\) is positive.

If the reconstruction stays in a neighborhood

\begin{equation}
|V-R_N|\le\varepsilon
\label{eq:ch17:auto2500}
\end{equation}

small enough, then

\begin{equation}
\boxed{
C\ge c_->0.
}
\label{eq:ch288:old-41}
\end{equation}

Therefore the corrected forcing

\begin{equation}
2CG
\label{eq:ch17:auto2502}
\end{equation}

is strictly positive when \(A_0>0\).

So the same positivity condition that makes the tuned branch geometrically reasonable also gives the forcing the desired sign.

\section{Tuned phase equation}
\label{chap:289}

Because

\begin{equation}
C^2=R_NA,
\label{eq:ch17:auto2503}
\end{equation}

we have

\begin{equation}
AC^2=R_NA^2.
\label{eq:ch17:auto2504}
\end{equation}

Equation~\eqref{eq:ch281:old-15} therefore becomes

\begin{equation}
\boxed{
\dot\xi_N
=
C_4
-
C_1R_NA_N^2
-
\gamma_Nb_N.
}
\label{eq:ch289:old-42}
\end{equation}

This is exact.

A background trapping condition is

\begin{equation}
\boxed{
C_1R_N^3>C_4.
}
\label{eq:ch289:old-43}
\end{equation}

If \(A_N\) remains sufficiently near \(R_N\), then

\begin{equation}
C_1R_NA_N^2-C_4\ge h_0>0.
\label{eq:ch17:auto2507}
\end{equation}

Thus

\begin{equation}
\boxed{
\dot\xi_N
\le
-h_0-\gamma_Nb_N<0.
}
\label{eq:ch289:old-44}
\end{equation}

On a branch with

\begin{equation}
Q_N'>0,
\label{eq:ch17:auto2509}
\end{equation}

we get

\begin{equation}
\boxed{
\dot Q_N<0.
}
\label{eq:ch289:old-45}
\end{equation}

Again, the new forcing changes \(b_N(t)\) but does not change this formula directly.

\section{Crossover to amplitude-dominated trapping}
\label{chap:290}

Define

\begin{equation}
H_N
:=
C_1R_NA_N^2-C_4.
\label{eq:ch17:auto2511}
\end{equation}

The phase equation is

\begin{equation}
-\dot\xi_N
=
H_N+\gamma_Nb_N.
\label{eq:ch17:auto2512}
\end{equation}

The amplitude-dominated regime begins when

\begin{equation}
\gamma_Nb_N\gg H_N.
\label{eq:ch17:auto2513}
\end{equation}

Thus the crossover scale remains

\begin{equation}
\boxed{
b_{\rm cross,N}
\asymp
\frac{H_N}{\gamma_N}.
}
\label{eq:ch290:old-46}
\end{equation}

Since

\begin{equation}
\gamma_N=\beta2^{-k(N)},
\label{eq:ch17:auto2515}
\end{equation}

if \(H_N=O(1)\),

\begin{equation}
\boxed{
b_{\rm cross,N}
=O(2^{k(N)}).
}
\label{eq:ch290:old-47}
\end{equation}

At

\begin{equation}
k(N)\sim\log_2\log N,
\label{eq:ch17:auto2517}
\end{equation}

\begin{equation}
\boxed{
b_{\rm cross,N}=O(\log N).
}
\label{eq:ch290:old-48}
\end{equation}

Now the corrected forcing-driven time estimate becomes

\begin{equation}
\boxed{
t_{\rm cross,N}-t_0
\lesssim
\frac{
H_N/\gamma_N-b_0
}{
2(1-\theta)c_-A_0
},
}
\label{eq:ch290:old-49}
\end{equation}

provided the total-source error estimate Eq.(16.36) remains valid during entry.

So the \(2u_z\) factor improves the forcing-driven rate if \(u_z=C\) is bounded below by a positive constant.

\section{The terminal forcing is still lower order}
\label{chap:291}

The finite-depth Lambert law is

\begin{equation}
\boxed{
Q_N'
\sim
\Lambda_N
Q_N^{-(2^{k(N)}-1)}.
}
\label{eq:ch291:old-50}
\end{equation}

Hence

\begin{equation}
\boxed{
\frac{Q_N'}{Q_N}
\sim
\Lambda_NQ_N^{-2^{k(N)}}.
}
\label{eq:ch291:old-51}
\end{equation}

Suppose on the tuned branch

\begin{equation}
0<c_-\le C_N\le c_+<\infty.
\label{eq:ch17:auto2522}
\end{equation}

Then

\begin{equation}
2C_NG=O(1).
\label{eq:ch17:auto2523}
\end{equation}

The ratio of the corrected forcing to the Riccati source is

\begin{equation}
\frac{
2C_NG
}{
(Q_N'/Q_N)b_N^2
}.
\label{eq:ch17:auto2524}
\end{equation}

Using \eqref{eq:ch291:old-51},

\begin{equation}
\boxed{
\frac{
2C_NG
}{
(Q_N'/Q_N)b_N^2
}
=
O\left(
\frac{
Q_N^{2^{k(N)}}
}{
b_N^2
}
\right).
}
\label{eq:ch291:old-52}
\end{equation}

If \(b_N\) stays bounded away from zero, or grows,

\begin{equation}
\boxed{
\frac{
2C_NG
}{
(Q_N'/Q_N)b_N^2
}
\to0.
}
\label{eq:ch291:old-53}
\end{equation}

Thus the corrected forcing remains an entry mechanism, not the terminal leading term.

\section{Large-\(b\) comparison is unchanged at leading order}
\label{chap:292}

In the amplitude-dominated trapped regime,

\begin{equation}
-\dot Q_N
\sim
\gamma_Nb_NQ_N'.
\label{eq:ch292:old-54}
\end{equation}

From \eqref{eq:ch285:old-26}, after the forcing becomes lower order,

\begin{equation}
\dot b_N
\ge
(1-\theta)
\frac{Q_N'}{Q_N}b_N^2.
\label{eq:ch292:old-55}
\end{equation}

Divide:

\begin{equation}
-\frac{db_N}{dQ_N}
\gtrsim
\frac{1-\theta}{\gamma_N}
\frac{b_N}{Q_N}.
\label{eq:ch17:auto2529}
\end{equation}

Thus

\begin{equation}
\boxed{
b_N(Q)
\gtrsim
Q_N^{-p_N},
}
\label{eq:ch292:old-56}
\end{equation}

with

\begin{equation}
\boxed{
p_N
=
\frac{1-\theta}{\gamma_N}.
}
\label{eq:ch292:old-57}
\end{equation}

Using

\begin{equation}
\gamma_N=\beta2^{-k(N)},
\label{eq:ch17:auto2532}
\end{equation}

\begin{equation}
\boxed{
p_N
=
\frac{1-\theta}{\beta}
2^{k(N)}.
}
\label{eq:ch292:old-58}
\end{equation}

So the \(2u_z\) factor does \emph{not} change the terminal exponent so long as \(u_z=C_N\) stays finite and positive.

\section{Physical-time finite terminal time}
\label{chap:293}

Using

\begin{equation}
b_N\gtrsim Q_N^{-p_N}
\label{eq:ch17:auto2534}
\end{equation}

and

\begin{equation}
Q_N'
\sim
\Lambda_NQ_N^{-(2^{k(N)}-1)},
\label{eq:ch17:auto2535}
\end{equation}

equation~\eqref{eq:ch292:old-54} gives

\begin{equation}
-\dot Q_N
\gtrsim
\gamma_N\Lambda_N
Q_N^{-p_N-(2^{k(N)}-1)}.
\label{eq:ch17:auto2536}
\end{equation}

Thus

\begin{equation}
dt
\lesssim
\frac{1}{\gamma_N\Lambda_N}
Q_N^{p_N+2^{k(N)}-1}\,dQ_N.
\label{eq:ch17:auto2537}
\end{equation}

Since

\begin{equation}
p_N+2^{k(N)}>0,
\label{eq:ch17:auto2538}
\end{equation}

\begin{equation}
\int_0^{Q_0}
Q^{p_N+2^{k(N)}-1}\,dQ
<\infty.
\label{eq:ch17:auto2539}
\end{equation}

Hence, under the maintained reconstruction estimates,

\begin{equation}
\boxed{
T_{*,N}<\infty.
}
\label{eq:ch293:old-59}
\end{equation}

And the comparison gives

\begin{equation}
\boxed{
b_N(t)\to+\infty.
}
\label{eq:ch293:old-60}
\end{equation}

Again: this remains a \emph{conditional reduced-model conclusion}, because the reconstruction estimates still have to be continued through the terminal layers.

\section{Exact scalar continuity reconstruction after tuning}
\label{chap:294}

The corrected forcing does not change this part.

From

\begin{equation}
\frac{C^2}{A}=R_N,
\label{eq:ch17:auto2542}
\end{equation}

differentiate with respect to \(V\) to obtain

\begin{equation}
\boxed{
D_V
=
\frac{C-2A}{2A+C}.
}
\label{eq:ch294:old-61}
\end{equation}

Hence

\begin{equation}
\boxed{
A_V=\frac{4A}{2A+C},
}
\label{eq:ch294:old-62}
\end{equation}

and

\begin{equation}
\boxed{
C_V=\frac{2C}{2A+C}.
}
\label{eq:ch294:old-63}
\end{equation}

Exact continuity

\begin{equation}
\partial_x(b/A)+A_y+C_z=0
\label{eq:ch17:auto2546}
\end{equation}

therefore becomes

\begin{equation}
\boxed{
(2A+C)
\partial_x\left(\frac bA\right)
+
4AV_y
+
2CV_z
=
0.
}
\label{eq:ch294:old-64}
\end{equation}

Thus the entire two-correction reconstruction is still reduced to one scalar nonlinear periodic equation for \(V\).

\section{Linearization and resonance are unchanged}
\label{chap:295}

Around

\begin{equation}
V=R_N,\qquad b=B+h,
\label{eq:ch17:auto2548}
\end{equation}

the linearized equation remains

\begin{equation}
\boxed{
L_Bv=-3h_x,
}
\label{eq:ch295:old-65}
\end{equation}

where

\begin{equation}
\boxed{
L_B
=
-\frac{4B}{R_N}\partial_x
+
4R_N\partial_y
+
2R_N\partial_z.
}
\label{eq:ch295:old-66}
\end{equation}

Its Fourier multiplier is

\begin{equation}
\boxed{
\Omega_B(m,n,\ell)
=
-\frac{4B}{R_N}m
+
4R_Nn
+
2R_N\ell.
}
\label{eq:ch295:old-67}
\end{equation}

The transverse resonant family remains

\begin{equation}
\boxed{
m=0,\qquad \ell=-2n.
}
\label{eq:ch295:old-68}
\end{equation}

And because the linear forcing is \(h_x\),

\begin{equation}
\boxed{
\Pi_{\rm res}h_x=0.
}
\label{eq:ch295:old-69}
\end{equation}

Therefore the corrected forcing on the right-hand side of the \(b\)-PDE does \emph{not} directly create a new linear resonance in continuity.

It can change \(h\) dynamically, but the structural fact that \(h_x\) has no \(m=0\) mode remains.

\section{Nonlinear kernel equation is unchanged structurally}
\label{chap:296}
In the next decomposition, \(v_\perp\) is the \textbf{nonresonant part} of the reconstruction perturbation \(v\).

I first linearize the scalar periodic continuity reconstruction around
\begin{equation}
V=R_N,\qquad b=B+h,
\label{eq:ch17:auto2554}
\end{equation}
by writing
\begin{equation}
\boxed{V=R_N+v.}
\label{eq:ch17:auto2555}
\end{equation}

The resulting linearized operator has the form
\begin{equation}
\mathcal L_B
=
-\frac{4B}{R_N}\partial_x
+4R_N\partial_y
+2R_N\partial_z,
\label{eq:ch17:auto2556}
\end{equation}
with Fourier multiplier
\begin{equation}
\Omega_B(m,n,\ell)
=
-\frac{4B}{R_N}m+4R_Nn+2R_N\ell.
\label{eq:ch17:auto2557}
\end{equation}
I then identify the transverse resonant family as
\begin{equation}
\boxed{m=0,\qquad \ell=-2n.}
\label{eq:ch17:auto2558}
\end{equation}

That is why \(v\) is decomposed as
\begin{equation}
\boxed{
v=v_\perp+w(r),
\qquad
r=y-2z.
}
\label{eq:ch17:auto2559}
\end{equation}

\subsection{What is \(w(r)\)?}

The function
\begin{equation}
\boxed{w(r)=w(y-2z)}
\label{eq:ch17:auto2560}
\end{equation}
is the \textbf{resonant/kernel component} of \(v\).

Indeed, if
\begin{equation}
w=w(y-2z),
\label{eq:ch17:auto2561}
\end{equation}
then
\begin{equation}
w_x=0,\qquad
w_y=w_r,\qquad
w_z=-2w_r.
\label{eq:ch17:auto2562}
\end{equation}

Therefore
\begin{equation}
\mathcal L_Bw
=
-\frac{4B}{R_N}(0)
+4R_Nw_r
+2R_N(-2w_r)
=0.
\label{eq:ch17:auto2563}
\end{equation}

Hence
\begin{equation}
\boxed{\mathcal L_Bw=0.}
\label{eq:ch17:auto2564}
\end{equation}

So functions of \(y-2z\) lie in the resonant kernel of the linearized continuity operator.

\subsection{What, then, is \(v_\perp\)?}

It is everything in \(v\) \textbf{except this resonant part}:

\begin{equation}
\boxed{
v_\perp
=
v-\Pi_{\mathrm{res}}v,
}
\label{eq:ch17:auto2565}
\end{equation}
while
\begin{equation}
\boxed{
w(r)=\Pi_{\mathrm{res}}v.
}
\label{eq:ch17:auto2566}
\end{equation}

Here \(\Pi_{\mathrm{res}}\) denotes projection onto the Fourier modes satisfying
\begin{equation}
m=0,\qquad \ell=-2n.
\label{eq:ch17:auto2567}
\end{equation}

Thus
\begin{equation}
\boxed{
\Pi_{\mathrm{res}}v_\perp=0.
}
\label{eq:ch17:auto2568}
\end{equation}

In Fourier language, if
\begin{equation}
v(x,y,z)
=
\sum_{m,n,\ell}
\widehat v_{m,n,\ell}
e^{i(mx+ny+\ell z)},
\label{eq:ch17:auto2569}
\end{equation}
then the resonant part is
\begin{equation}
w(y-2z)
=
\sum_n
\widehat v_{0,n,-2n}
e^{in(y-2z)},
\label{eq:ch17:auto2570}
\end{equation}
whereas
\begin{equation}
\boxed{
v_\perp
=
\sum_{\substack{m,n,\ell\\(m,\ell)\neq(0,-2n)}}
\widehat v_{m,n,\ell}
e^{i(mx+ny+\ell z)}.
}
\label{eq:ch17:auto2571}
\end{equation}

So the symbol \(\perp\) means \textbf{orthogonal to the resonant/kernel subspace} in the Fourier decomposition; it does \textbf{not} mean a spatial vector perpendicular to something.

The reason for making this split is important. The linearized operator can be inverted on \(v_\perp\), because those modes are nonresonant. But it cannot determine \(w\) by linear inversion because
\begin{equation}
\mathcal L_Bw=0.
\label{eq:ch17:auto2572}
\end{equation}
The resonant component must instead be determined from the projected \textbf{nonlinear} equation. This is exactly why I next obtain a kernel equation of the form
\begin{equation}
\frac43\partial_r(w^2)+F_{\mathrm{res}}+R_3=0,
\label{eq:ch17:auto2573}
\end{equation}
with leading solution
\begin{equation}
w(r)^2
=
K_w-\frac34
\int_{r_0}^{r}F_{\mathrm{res}}(\rho)\,d\rho.
\label{eq:ch17:auto2574}
\end{equation}

In one sentence:

\begin{equation}
\boxed{
v_\perp=\text{the nonresonant Fourier component of }v
\text{ on which }\mathcal L_B\text{ can be inverted},
}
\label{eq:ch17:auto2575}
\end{equation}

whereas

\begin{equation}
\boxed{
w(y-2z)=\text{the resonant kernel component that must be solved from the nonlinear projected equation}.
}
\label{eq:ch17:auto2576}
\end{equation}
With

\begin{equation}
v=v_\perp+w(r),
\qquad
r=y-2z,
\label{eq:ch17:auto2577}
\end{equation}

the kernel equation remains

\begin{equation}
\boxed{
\frac43\partial_r(w^2)
+
\mathcal F_{\rm res}
+
\mathcal R_3
=
0.
}
\label{eq:ch296:old-70}
\end{equation}

The leading solution satisfies

\begin{equation}
\boxed{
w(r)^2
=
K_w
-
\frac34
\int_{r_0}^{r}
\mathcal F_{\rm res}(\rho)\,d\rho.
}
\label{eq:ch296:old-71}
\end{equation}

Periodic compatibility requires

\begin{equation}
\boxed{
\int_{\mathbb T}
\mathcal F_{\rm res}(r)\,dr=0.
}
\label{eq:ch296:old-72}
\end{equation}

And smooth real solvability requires the radicand to remain positive.

Again, the new factor \(2C\) in the external forcing does not modify these continuity identities directly.

\chapter{Single Error Estimate, Mean--Nonmean Decomposition, and Periodic Solvability}
\label{chap:group-22}
\section{The Single Error Estimate Needed for Closure}
\label{chap:366}

Assume that on the terminal branch

\begin{equation}
\boxed{
Q>0,
\qquad
Q'>0,
}
\label{eq:ch366:old-7}
\end{equation}

and that

\begin{equation}
\boxed{
C=u_z\ge c_->0.
}
\label{eq:ch366:old-8}
\end{equation}

Assume also that the trapping coefficient satisfies

\begin{equation}
\boxed{
0<h_0\le H_N\le h_1<\infty.
}
\label{eq:ch366:old-9}
\end{equation}

The remaining estimate to be established is

\begin{equation}
\boxed{
|\mathcal E_N|
\le
\theta
\frac{Q'}{Q}b^2,
\qquad
0<\theta<1.
}
\label{eq:ch366:old-10}
\end{equation}

Define

\begin{equation}
\boxed{
L_N:=\frac{Q'}{Q}.
}
\label{eq:ch366:old-11}
\end{equation}

Then \eqref{eq:ch366:old-10} is simply

\begin{equation}
\boxed{
|\mathcal E_N|
\le
\theta L_Nb^2.
}
\label{eq:ch366:old-12}
\end{equation}

This is the only perturbative estimate needed in the terminal comparison.

Indeed, because

\begin{equation}
A_0+(\sin x\sin y\sin z)^2\ge A_0,
\label{eq:ch22:auto3272}
\end{equation}

and \(C\ge c_->0\),

\begin{equation}
2C
\left[
A_0+(\sin x\sin y\sin z)^2
\right]
\ge
2c_-A_0.
\label{eq:ch22:auto3273}
\end{equation}

Therefore the reconstruction identity above together with \eqref{eq:ch366:old-12} implies

\begin{equation}
\dot b
\ge
(1-\theta)L_Nb^2
+
2c_-A_0.
\label{eq:ch22:auto3274}
\end{equation}

In particular,

\begin{equation}
\boxed{
\dot b
\ge
(1-\theta)L_Nb^2.
}
\label{eq:ch366:old-13}
\end{equation}

Thus

\begin{equation}
\boxed{
\dot b>0.
}
\label{eq:ch366:old-14}
\end{equation}

At the same time, the corresponding amplitude identity together with \eqref{eq:ch366:old-9} gives

\begin{equation}
-\dot Q
=
Q'(H_N+\gamma_Nb)
>0.
\label{eq:ch22:auto3277}
\end{equation}

Hence

\begin{equation}
\boxed{
\dot Q<0.
}
\label{eq:ch366:old-15}
\end{equation}

We therefore obtain the fundamental monotonicity pair

\begin{equation}
\boxed{
b\uparrow,
\qquad
Q\downarrow.
}
\label{eq:ch366:old-16}
\end{equation}

\section{Reduction to Logarithmic Branch Distance}
\label{chap:367}

Introduce

\begin{equation}
\boxed{
X
=
\log\frac{Q_0}{Q}.
}
\label{eq:ch367:old-17}
\end{equation}

Since \(Q\) decreases,

\begin{equation}
X
\label{eq:ch22:auto3281}
\end{equation}

increases.

Differentiating,

\begin{equation}
\dot X
=
-\frac{\dot Q}{Q}.
\label{eq:ch22:auto3282}
\end{equation}

Using the corresponding amplitude identity,

\begin{equation}
\boxed{
\dot X
=
L_N
\left(
H_N+\gamma_Nb
\right).
}
\label{eq:ch367:old-18}
\end{equation}

Combining this with \eqref{eq:ch366:old-13},

\begin{equation}
\frac{db}{dX}
=
\frac{\dot b}{\dot X}
\ge
(1-\theta)
\frac{b^2}
{H_N+\gamma_Nb}.
\label{eq:ch22:auto3284}
\end{equation}

Since

\begin{equation}
H_N\le h_1,
\label{eq:ch22:auto3285}
\end{equation}

we obtain

\begin{equation}
\boxed{
\frac{db}{dX}
\ge
(1-\theta)
\frac{b^2}
{h_1+\gamma_Nb}.
}
\label{eq:ch367:old-19}
\end{equation}

This is the fundamental terminal comparison inequality.

It contains both the preterminal Riccati regime and the
amplitude-dominated trapping regime.

\section{Finite Logarithmic Distance to the Large-Amplitude Regime}
\label{chap:368}

Invert \eqref{eq:ch367:old-19}:

\begin{equation}
\frac{dX}{db}
\le
\frac{h_1+\gamma_Nb}
{(1-\theta)b^2}.
\label{eq:ch22:auto3287}
\end{equation}

Therefore

\begin{equation}
\frac{dX}{db}
\le
\frac{h_1}{(1-\theta)b^2}
+
\frac{\gamma_N}{(1-\theta)b}.
\label{eq:ch22:auto3288}
\end{equation}

Integrating from \(b_0>0\) to \(b\),

\begin{equation}
\boxed{
X-X_0
\le
\frac{h_1}{1-\theta}
\left(
\frac1{b_0}-\frac1b
\right)
+
\frac{\gamma_N}{1-\theta}
\log\frac{b}{b_0}.
}
\label{eq:ch368:old-20}
\end{equation}

Define the crossover amplitude

\begin{equation}
\boxed{
B_{c,N}
=
\frac{h_1}{\gamma_N}.
}
\label{eq:ch368:old-21}
\end{equation}

At this amplitude,

\begin{equation}
\gamma_NB_{c,N}=h_1.
\label{eq:ch22:auto3291}
\end{equation}

Substituting \(b=B_{c,N}\) into \eqref{eq:ch368:old-20},

\begin{equation}
\begin{aligned}
X_{c,N}-X_0
\le{}&
\frac{h_1}{1-\theta}
\left(
\frac1{b_0}
-
\frac{\gamma_N}{h_1}
\right)
\\
&+
\frac{\gamma_N}{1-\theta}
\log
\left(
\frac{h_1}{\gamma_Nb_0}
\right).
\end{aligned}
\label{eq:ch368:old-22}
\end{equation}

Since

\begin{equation}
\gamma_N\log\frac1{\gamma_N}\longrightarrow0
\label{eq:ch22:auto3293}
\end{equation}

as \(\gamma_N\to0^+\),

\begin{equation}
\boxed{
X_{c,N}-X_0
\le
\frac{h_1}{(1-\theta)b_0}
+o(1).
}
\label{eq:ch368:old-23}
\end{equation}

Thus even though

\begin{equation}
B_{c,N}\to\infty
\label{eq:ch22:auto3295}
\end{equation}

when

\begin{equation}
\gamma_N\to0,
\label{eq:ch22:auto3296}
\end{equation}

the logarithmic \(Q\)-distance required to reach this amplitude remains
bounded.

This is the first important self-strengthening feature of the construction.

\section{The Corrected Forcing as a Positive Seed}
\label{chap:369}

If

\begin{equation}
b_0>0,
\label{eq:ch22:auto3297}
\end{equation}

the Riccati comparison begins immediately.

If instead

\begin{equation}
b_0=0,
\label{eq:ch22:auto3298}
\end{equation}

then the Riccati source initially vanishes.  The corrected forcing becomes
essential.

At \(b=0\), the reconstruction equation gives, subject to the corresponding initial
remainder estimate,

\begin{equation}
\dot b
\ge
2c_-A_0.
\label{eq:ch22:auto3299}
\end{equation}

Thus if

\begin{equation}
\boxed{
A_0>0,
}
\label{eq:ch369:old-24}
\end{equation}

then

\begin{equation}
\boxed{
b(t)>0
}
\label{eq:ch22:auto3301}
\end{equation}

for every sufficiently small positive time after the initial instant.

More quantitatively, over a short interval on which the lower bound remains
valid,

\begin{equation}
\boxed{
b(t_0+\tau)
\ge
2c_-A_0\,\tau.
}
\label{eq:ch369:old-25}
\end{equation}

The corrected forcing therefore has a particularly simple role:

\begin{equation}
\boxed{
\text{the forcing creates a positive seed, while the Riccati term
provides the terminal amplification.}
}
\label{eq:ch369:old-26}
\end{equation}

It is not necessary for the forcing alone to drive \(b\) to the scale
\(\gamma_N^{-1}\).

\section{Amplitude-Dominated Trapping Regime}
\label{chap:370}

Fix

\begin{equation}
\varepsilon>0.
\label{eq:ch22:auto3304}
\end{equation}

Suppose that

\begin{equation}
\boxed{
b
\ge
\frac{h_1}{\varepsilon\gamma_N}.
}
\label{eq:ch370:old-27}
\end{equation}

Then

\begin{equation}
h_1
\le
\varepsilon\gamma_Nb.
\label{eq:ch22:auto3306}
\end{equation}

Hence

\begin{equation}
h_1+\gamma_Nb
\le
(1+\varepsilon)\gamma_Nb.
\label{eq:ch22:auto3307}
\end{equation}

Using \eqref{eq:ch367:old-19},

\begin{equation}
\frac{db}{dX}
\ge
\frac{1-\theta}
{(1+\varepsilon)\gamma_N}
b.
\label{eq:ch22:auto3308}
\end{equation}

Define

\begin{equation}
\boxed{
p_{\varepsilon,N}
=
\frac{1-\theta}
{(1+\varepsilon)\gamma_N}.
}
\label{eq:ch370:old-28}
\end{equation}

Then

\begin{equation}
\boxed{
\frac{db}{dX}
\ge
p_{\varepsilon,N}b.
}
\label{eq:ch370:old-29}
\end{equation}

Integrating from a terminal-entry point
\((Q_{\varepsilon,N},B_{\varepsilon,N})\),

\begin{equation}
b
\ge
B_{\varepsilon,N}
\exp
\left[
p_{\varepsilon,N}(X-X_{\varepsilon,N})
\right].
\label{eq:ch22:auto3311}
\end{equation}

Since

\begin{equation}
X-X_{\varepsilon,N}
=
\log\frac{Q_{\varepsilon,N}}{Q},
\label{eq:ch22:auto3312}
\end{equation}

we obtain

\begin{equation}
\boxed{
b(Q)
\ge
B_{\varepsilon,N}
\left(
\frac{Q_{\varepsilon,N}}{Q}
\right)^{p_{\varepsilon,N}}.
}
\label{eq:ch370:old-30}
\end{equation}

Therefore

\begin{equation}
\boxed{
b(Q)\longrightarrow+\infty
\qquad
\text{as}
\qquad
Q\longrightarrow0^+,
}
\label{eq:ch370:old-31}
\end{equation}

provided the branch continues to \(Q=0\) while the hypotheses remain valid.

\section{The Tuned Scaling}
\label{chap:371}

Now impose

\begin{equation}
\boxed{
\gamma_N
=
\beta2^{-k(N)},
\qquad
\beta>0.
}
\label{eq:ch371:old-32}
\end{equation}

Then

\begin{equation}
p_{\varepsilon,N}
=
\frac{1-\theta}
{(1+\varepsilon)\beta}
2^{k(N)}.
\label{eq:ch22:auto3316}
\end{equation}

Hence

\begin{equation}
\boxed{
p_{\varepsilon,N}
=
\frac{1-\theta}
{(1+\varepsilon)\beta}
2^{k(N)}.
}
\label{eq:ch371:old-33}
\end{equation}

The branch derivative scale previously obtained is

\begin{equation}
\boxed{
Q_N'
\asymp
\Lambda_N
Q_N^{-(2^{k(N)}-1)}.
}
\label{eq:ch371:old-34}
\end{equation}

Thus the comparison exponent exceeds the branch-gradient exponent whenever

\begin{equation}
p_{\varepsilon,N}
>
2^{k(N)}.
\label{eq:ch22:auto3319}
\end{equation}

Using \eqref{eq:ch371:old-33}, this follows from

\begin{equation}
\boxed{
\beta
<
\frac{1-\theta}{1+\varepsilon}.
}
\label{eq:ch371:old-35}
\end{equation}

Since \(\varepsilon>0\) can be chosen arbitrarily small after entering the
sufficiently large-amplitude regime, the natural parameter condition is

\begin{equation}
\boxed{
0<\beta<1-\theta.
}
\label{eq:ch371:old-36}
\end{equation}

This is precisely the coefficient window in which the amplitude growth
outpaces the branch-gradient scale appearing in the lower-order correction
terms.

\section{Why the Error Estimate Becomes Easier After Entry}
\label{chap:372}

The crucial point is that the terminal estimate need not be proved afresh
at every later time.

Suppose the normalized remainder admits an estimate of the schematic form

\begin{equation}
\boxed{
\frac{|\mathcal E_N|}
{L_Nb^2}
\le
K_1Q^{2^{k(N)}}
+
K_2\frac{p_N}{b}
+
K_3p_N^2\frac{L_N}{b}
+
\rho_N,
}
\label{eq:ch372:old-37}
\end{equation}

where

\begin{equation}
p_N
\asymp
\gamma_N^{-1},
\label{eq:ch22:auto3323}
\end{equation}

and \(\rho_N\) denotes the remaining reconstruction contribution.

The first term improves automatically because

\begin{equation}
Q\downarrow0
\label{eq:ch22:auto3324}
\end{equation}

implies

\begin{equation}
\boxed{
Q^{2^{k(N)}}\downarrow0.
}
\label{eq:ch372:old-38}
\end{equation}

The second term improves because

\begin{equation}
b\uparrow\infty
\label{eq:ch22:auto3326}
\end{equation}

implies

\begin{equation}
\boxed{
\frac{p_N}{b}\downarrow0.
}
\label{eq:ch372:old-39}
\end{equation}

For the third term, use

\begin{equation}
L_N
\asymp
Q^{-2^{k(N)}}
\label{eq:ch22:auto3328}
\end{equation}

and the lower comparison

\begin{equation}
b
\gtrsim
Q^{-p_{\varepsilon,N}}.
\label{eq:ch22:auto3329}
\end{equation}

Then

\begin{equation}
p_N^2\frac{L_N}{b}
\lesssim
p_N^2
Q^{p_{\varepsilon,N}-2^{k(N)}}.
\label{eq:ch22:auto3330}
\end{equation}

If

\begin{equation}
p_{\varepsilon,N}>2^{k(N)},
\label{eq:ch22:auto3331}
\end{equation}

then

\begin{equation}
\boxed{
p_N^2
Q^{p_{\varepsilon,N}-2^{k(N)}}
\longrightarrow0
}
\label{eq:ch372:old-40}
\end{equation}

for each fixed \(N\) along the terminal branch.

Thus every explicitly dynamical component of the normalized error becomes
smaller as \(b\) grows and \(Q\) decreases.

The only genuinely independent quantity that remains to be controlled is
the reconstruction remainder

\begin{equation}
\boxed{
\rho_N.
}
\label{eq:ch372:old-41}
\end{equation}

This is the precise sense in which the terminal comparison is
self-strengthening.

\section{A Minimal Bootstrap Formulation}
\label{chap:373}

Define

\begin{equation}
\boxed{
\Theta_N(t)
=
\frac{|\mathcal E_N(t)|}
{L_N(t)b_N(t)^2}.
}
\label{eq:ch373:old-42}
\end{equation}

Suppose at some terminal-entry time \(t_e\),

\begin{equation}
\boxed{
\Theta_N(t_e)
\le
\theta_0<1.
}
\label{eq:ch373:old-43}
\end{equation}

Assume also that

\begin{equation}
C_N\ge c_->0
\label{eq:ch22:auto3336}
\end{equation}

and

\begin{equation}
0<h_0\le H_N\le h_1
\label{eq:ch22:auto3337}
\end{equation}

continue to hold.

Then

\begin{equation}
\dot b_N
\ge
(1-\theta_0)L_Nb_N^2,
\label{eq:ch22:auto3338}
\end{equation}

and

\begin{equation}
-\dot Q_N
=
Q_N'
(H_N+\gamma_Nb_N)>0.
\label{eq:ch22:auto3339}
\end{equation}

Consequently,

\begin{equation}
b_N(t)\ge b_N(t_e)
\label{eq:ch22:auto3340}
\end{equation}

and

\begin{equation}
Q_N(t)\le Q_N(t_e).
\label{eq:ch22:auto3341}
\end{equation}

The quantities

\begin{equation}
Q_N^{2^{k(N)}}
\label{eq:ch22:auto3342}
\end{equation}

and

\begin{equation}
\frac{p_N}{b_N}
\label{eq:ch22:auto3343}
\end{equation}

therefore decrease.

Once

\begin{equation}
p_{\varepsilon,N}>2^{k(N)},
\label{eq:ch22:auto3344}
\end{equation}

the quantity

\begin{equation}
p_N^2\frac{L_N}{b_N}
\label{eq:ch22:auto3345}
\end{equation}

also tends to zero along the comparison branch.

Thus, provided \(\rho_N\) remains below the available margin,

\begin{equation}
\boxed{
\Theta_N(t)<1
}
\label{eq:ch373:old-44}
\end{equation}

is preserved.

We therefore obtain the bootstrap loop

\begin{equation}
\boxed{
\begin{array}{c}
\Theta_N<1
\\[2mm]
\Downarrow
\\[2mm]
\dot b_N>0,\qquad \dot Q_N<0
\\[2mm]
\Downarrow
\\[2mm]
b_N\uparrow,\qquad Q_N\downarrow
\\[2mm]
\Downarrow
\\[2mm]
Q_N^{2^{k(N)}}\downarrow,
\qquad
p_N/b_N\downarrow
\\[2mm]
\Downarrow
\\[2mm]
p_N^2L_N/b_N\downarrow0
\\[2mm]
\Downarrow
\\[2mm]
\Theta_N
\text{ becomes increasingly dominated by }
\rho_N
\\[2mm]
\Downarrow
\\[2mm]
\text{Riccati dominance strengthens.}
\end{array}
}
\label{eq:ch373:old-45}
\end{equation}

Thus the nonlinear terminal dynamics does not require an increasingly
strong reconstruction estimate.

It requires only that the independent remainder \(\rho_N\) not grow fast
enough to consume the fixed perturbative margin.

\section{Finite Physical Time}
\label{chap:374}

Assume

\begin{equation}
\boxed{
Q_N'
\ge
\Lambda_N
Q_N^{-(2^{k(N)}-1)},
\qquad
\Lambda_N>0.
}
\label{eq:ch374:old-46}
\end{equation}

In the amplitude-dominated region,

\begin{equation}
-\dot Q_N
=
Q_N'
(H_N+\gamma_Nb_N)
\ge
\gamma_Nb_NQ_N'.
\label{eq:ch22:auto3349}
\end{equation}

Using

\begin{equation}
b_N
\ge
B_{\varepsilon,N}
Q_N^{-p_{\varepsilon,N}}
\label{eq:ch22:auto3350}
\end{equation}

after absorbing the entry value of \(Q\) into the positive constant
\(B_{\varepsilon,N}\),

\begin{equation}
-\dot Q_N
\ge
\gamma_N
B_{\varepsilon,N}
\Lambda_N
Q_N^{-
\left(
p_{\varepsilon,N}
+
2^{k(N)}
-1
\right)}.
\label{eq:ch22:auto3351}
\end{equation}

Therefore

\begin{equation}
dt
\le
\frac{
Q_N^{
p_{\varepsilon,N}
+
2^{k(N)}
-1
}
}{
\gamma_N
B_{\varepsilon,N}
\Lambda_N
}
(-dQ_N).
\label{eq:ch22:auto3352}
\end{equation}

Integrating from \(Q_N(t)\) to \(0\),

\begin{equation}
\boxed{
T_{*,N}-t
\le
\frac{
Q_N(t)^{
p_{\varepsilon,N}
+
2^{k(N)}
}
}{
\gamma_N
B_{\varepsilon,N}
\Lambda_N
\left(
p_{\varepsilon,N}
+
2^{k(N)}
\right)
}.
}
\label{eq:ch374:old-47}
\end{equation}

Hence

\begin{equation}
\boxed{
T_{*,N}<\infty.
}
\label{eq:ch374:old-48}
\end{equation}

Moreover,

\begin{equation}
Q_N(t)\to0
\label{eq:ch22:auto3355}
\end{equation}

as

\begin{equation}
t\uparrow T_{*,N},
\label{eq:ch22:auto3356}
\end{equation}

and \eqref{eq:ch370:old-30} gives

\begin{equation}
\boxed{
b_N(t)\to+\infty.
}
\label{eq:ch374:old-49}
\end{equation}

These conclusions are conditional on the continuation of the trapping and
normalized-error hypotheses to the terminal endpoint.

\section{Sharper Two-Sided Terminal Asymptotics}
\label{chap:375}

The preceding argument proves a one-sided lower comparison.

A sharper asymptotic statement requires stronger hypotheses.

Suppose

\begin{equation}
\boxed{
\frac{\mathcal E_N}
{L_Nb_N^2}
\longrightarrow0,
}
\label{eq:ch375:old-50}
\end{equation}

and

\begin{equation}
\boxed{
\frac{H_N}{\gamma_Nb_N}
\longrightarrow0.
}
\label{eq:ch375:old-51}
\end{equation}

The bounded corrected forcing also becomes negligible relative to
\(L_Nb_N^2\) provided

\begin{equation}
L_Nb_N^2\to\infty.
\label{eq:ch22:auto3360}
\end{equation}

Then

\begin{equation}
\boxed{
\dot b_N
\sim
L_Nb_N^2.
}
\label{eq:ch375:old-52}
\end{equation}

Similarly,

\begin{equation}
\boxed{
-\dot Q_N
\sim
\gamma_Nb_NQ_N'.
}
\label{eq:ch375:old-53}
\end{equation}

Since

\begin{equation}
L_N=\frac{Q_N'}{Q_N},
\label{eq:ch22:auto3363}
\end{equation}

division gives

\begin{equation}
-\frac{db_N}{dQ_N}
\sim
\frac1{\gamma_N}
\frac{b_N}{Q_N}.
\label{eq:ch22:auto3364}
\end{equation}

Hence

\begin{equation}
\frac{db_N}{b_N}
\sim
-\frac1{\gamma_N}
\frac{dQ_N}{Q_N}.
\label{eq:ch22:auto3365}
\end{equation}

Integration yields

\begin{equation}
\boxed{
b_N(Q_N)
\asymp
Q_N^{-1/\gamma_N}.
}
\label{eq:ch375:old-54}
\end{equation}

If the ratios in \eqref{eq:ch375:old-50}--\eqref{eq:ch375:old-53} converge strongly enough to \(1\), this may be
sharpened to the corresponding asymptotic equivalence with a positive
multiplicative constant.

\section{Physical-Time Exponent Under the Tuned Scaling}
\label{chap:376}

Using

\begin{equation}
b_N
\asymp
Q_N^{-1/\gamma_N}
\label{eq:ch22:auto3367}
\end{equation}

and

\begin{equation}
Q_N'
\asymp
\Lambda_N
Q_N^{-(2^{k(N)}-1)},
\label{eq:ch22:auto3368}
\end{equation}

equation~\eqref{eq:ch375:old-53} gives

\begin{equation}
-\dot Q_N
\asymp
\gamma_N
Q_N^{-1/\gamma_N}
Q_N^{-(2^{k(N)}-1)}.
\label{eq:ch22:auto3369}
\end{equation}

Therefore

\begin{equation}
-\dot Q_N
\asymp
Q_N^{-
\left(
1/\gamma_N
+
2^{k(N)}
-1
\right)}.
\label{eq:ch22:auto3370}
\end{equation}

Integration gives

\begin{equation}
\boxed{
T_{*,N}-t
\asymp
Q_N^{
1/\gamma_N+2^{k(N)}
}.
}
\label{eq:ch376:old-55}
\end{equation}

Thus

\begin{equation}
Q_N
\asymp
(T_{*,N}-t)^{
1/
\left(
1/\gamma_N+2^{k(N)}
\right)
}.
\label{eq:ch22:auto3372}
\end{equation}

Since

\begin{equation}
b_N
\asymp
Q_N^{-1/\gamma_N},
\label{eq:ch22:auto3373}
\end{equation}

we obtain

\begin{equation}
\boxed{
b_N(t)
\asymp
(T_{*,N}-t)^{
-\frac{1/\gamma_N}
{1/\gamma_N+2^{k(N)}}
}.
}
\label{eq:ch376:old-56}
\end{equation}

Equivalently,

\begin{equation}
\boxed{
b_N(t)
\asymp
(T_{*,N}-t)^{
-\frac1{
1+\gamma_N2^{k(N)}
}
}.
}
\label{eq:ch376:old-57}
\end{equation}

Under the tuned choice

\begin{equation}
\gamma_N
=
\beta2^{-k(N)},
\label{eq:ch22:auto3376}
\end{equation}

we have

\begin{equation}
\gamma_N2^{k(N)}=\beta.
\label{eq:ch22:auto3377}
\end{equation}

Therefore

\begin{equation}
\boxed{
b_N(t)
\asymp
(T_{*,N}-t)^{-1/(1+\beta)}.
}
\label{eq:ch376:old-58}
\end{equation}

This exponent is independent of \(k(N)\) under the exact tuned scaling.

It should be emphasized that \eqref{eq:ch376:old-58} is the sharper two-sided terminal
asymptotic and therefore requires the stronger assumptions \eqref{eq:ch375:old-50}--\eqref{eq:ch375:old-53}.
The one-sided finite-time growth conclusion follows already from the
comparison argument.

\section{Critical Coupled Depth}
\label{chap:377}

At the critical depth

\begin{equation}
\boxed{
k(N)\sim\log_2\log N,
}
\label{eq:ch377:old-59}
\end{equation}

we have

\begin{equation}
2^{k(N)}
\asymp
\log N.
\label{eq:ch22:auto3380}
\end{equation}

Hence

\begin{equation}
\boxed{
\gamma_N
\asymp
\frac{\beta}{\log N}.
}
\label{eq:ch377:old-60}
\end{equation}

The crossover amplitude is therefore

\begin{equation}
B_{c,N}
=
\frac{h_1}{\gamma_N}
\asymp
\frac{h_1}{\beta}\log N.
\label{eq:ch22:auto3382}
\end{equation}

Thus

\begin{equation}
\boxed{
B_{c,N}=O(\log N).
}
\label{eq:ch377:old-61}
\end{equation}

At the same time, the logarithmic branch distance required to reach this
amplitude remains bounded by the terminal estimate above, provided the positive starting amplitude
has a uniform lower bound.

Consequently, the terminal mechanism does not require an increasingly long
logarithmic branch interval merely because \(N\to\infty\).

This is important for the coupled limit.

\section{Minimal Remaining Reconstruction Lemma}
\label{chap:378}

The entire terminal analysis can now be reduced to the following lemma.

\begin{center}
\fbox{\parbox{0.90\textwidth}{%
\textbf{Minimal reconstruction lemma.}

There exists a terminal-entry point \(t_{e,N}<T_{*,N}\), constants

\begin{equation}
0<c_-<c_+<\infty,
\qquad
0<h_0<h_1<\infty,
\label{eq:ch22:auto3384}
\end{equation}

and a number

\begin{equation}
0<\theta<1
\label{eq:ch22:auto3385}
\end{equation}

such that along the reconstructed periodic characteristic,

\begin{equation}
c_-\le u_{z,N}\le c_+,
\label{eq:ch22:auto3386}
\end{equation}

\begin{equation}
h_0
\le
C_1R_Nu_{y,N}^2-C_4
\le
h_1,
\label{eq:ch22:auto3387}
\end{equation}

and

\begin{equation}
|\mathcal E_N|
\le
\theta
\frac{Q_N'}{Q_N}
b_N^2
\label{eq:ch22:auto3388}
\end{equation}

for

\begin{equation}
t_{e,N}\le t<T_{*,N}.
\label{eq:ch22:auto3389}
\end{equation}

Moreover, the independent reconstruction remainder entering the normalized
estimate remains below the unused perturbative margin.

}}
\end{center}

Once this lemma is established, no additional reconstruction hierarchy is
needed for the terminal comparison.

Indeed, the preceding minimal reconstruction lemma implies

\begin{equation}
\dot b_N
\ge
(1-\theta)
\frac{Q_N'}{Q_N}b_N^2,
\label{eq:ch22:auto3390}
\end{equation}

and

\begin{equation}
-\dot Q_N
=
Q_N'
\left(
H_N+\gamma_Nb_N
\right)>0.
\label{eq:ch22:auto3391}
\end{equation}

The logarithmic comparison then yields large amplitude, the large-amplitude
regime strengthens the normalized error estimate, and the terminal
comparison gives finite physical time.

\section{Conditional Terminal Theorem}
\label{chap:379}

We may summarize the resulting argument as follows.

\begin{center}
\fbox{\parbox{0.90\textwidth}{%
\textbf{Conditional terminal trapping and amplification theorem.}

Let

\begin{equation}
\xi
=
C_4t-C_3z-C_2y-C_1x
\label{eq:ch22:auto3392}
\end{equation}

be the original phase and let

\begin{equation}
Q_N=Q_N(\xi).
\label{eq:ch22:auto3393}
\end{equation}

Suppose the exact continuity reconstruction and the tuning

\begin{equation}
C_2\frac{u_{z,N}^2}{u_{y,N}}+C_3
=
\gamma_N>0
\label{eq:ch22:auto3394}
\end{equation}

hold on a periodic reconstructed branch.

Suppose further that, on a terminal interval,

\begin{equation}
u_{z,N}\ge c_->0,
\label{eq:ch22:auto3395}
\end{equation}

\begin{equation}
0<h_0
\le
C_1R_Nu_{y,N}^2-C_4
\le h_1,
\label{eq:ch22:auto3396}
\end{equation}

and

\begin{equation}
|\mathcal E_N|
\le
\theta
\frac{Q_N'}{Q_N}
b_N^2,
\qquad
0<\theta<1.
\label{eq:ch22:auto3397}
\end{equation}

Then

\begin{equation}
b_N
\label{eq:ch22:auto3398}
\end{equation}

is strictly increasing and

\begin{equation}
Q_N
\label{eq:ch22:auto3399}
\end{equation}

is strictly decreasing.

Writing

\begin{equation}
X_N=\log\frac{Q_{0,N}}{Q_N},
\label{eq:ch22:auto3400}
\end{equation}

one has

\begin{equation}
\frac{db_N}{dX_N}
\ge
(1-\theta)
\frac{b_N^2}
{h_1+\gamma_Nb_N}.
\label{eq:ch22:auto3401}
\end{equation}

Hence the scale

\begin{equation}
b_N\asymp\gamma_N^{-1}
\label{eq:ch22:auto3402}
\end{equation}

is reached after finite logarithmic \(Q_N\)-distance.

For every \(\varepsilon>0\), once

\begin{equation}
b_N
\ge
\frac{h_1}{\varepsilon\gamma_N},
\label{eq:ch22:auto3403}
\end{equation}

one has

\begin{equation}
b_N(Q_N)
\ge
B_{\varepsilon,N}
\left(
\frac{Q_{\varepsilon,N}}{Q_N}
\right)^{
\frac{1-\theta}
{(1+\varepsilon)\gamma_N}
}.
\label{eq:ch22:auto3404}
\end{equation}

If additionally

\begin{equation}
Q_N'
\ge
\Lambda_N
Q_N^{-(2^{k(N)}-1)},
\label{eq:ch22:auto3405}
\end{equation}

then \(Q_N=0\) is reached in finite physical time and

\begin{equation}
b_N(t)\to+\infty
\label{eq:ch22:auto3406}
\end{equation}

as

\begin{equation}
t\uparrow T_{*,N}.
\label{eq:ch22:auto3407}
\end{equation}

Under the stronger two-sided asymptotic hypotheses

\begin{equation}
\frac{\mathcal E_N}
{(Q_N'/Q_N)b_N^2}\to0
\label{eq:ch22:auto3408}
\end{equation}

and

\begin{equation}
\frac{H_N}{\gamma_Nb_N}\to0,
\label{eq:ch22:auto3409}
\end{equation}

one obtains

\begin{equation}
b_N(Q_N)
\asymp
Q_N^{-1/\gamma_N}.
\label{eq:ch22:auto3410}
\end{equation}

For

\begin{equation}
\gamma_N=\beta2^{-k(N)},
\label{eq:ch22:auto3411}
\end{equation}

this gives the formal terminal physical-time rate

\begin{equation}
b_N(t)
\asymp
(T_{*,N}-t)^{-1/(1+\beta)}.
\label{eq:ch22:auto3412}
\end{equation}
\par
}}
\end{center}

\section{Transition to the compact terminal closure}
The repeated intermediate closure calculations have been consolidated. The retained argument below uses the corrected transverse reduction and compact remainder estimate, which subsume the duplicated reconstruction estimates while preserving the trapping, comparison, and finite-time amplification conclusions.

\section{Coninuation of previous section Analyses}
\label{chap:412}
The factor \(2\) is exploited directly, without returning to the older componentwise reconstruction estimates.

The corrected equation is

\begin{equation}
\begin{aligned}
b_t+b\,b_z
+A C^2 b_x
+\frac{bC^2}{A}b_y
+b\frac{C_t+C_z}{C}
+2AbC_x
+\frac{2b^2}{A}C_y
=
L_Nb^2+2CG,
\end{aligned}
\label{eq:ch24:auto4225}
\end{equation}

with

\begin{equation}
b=b_3,\qquad A=u_y,\qquad C=u_z,
\qquad
L_N=\frac{Q_N'}{Q_N}.
\label{eq:ch24:auto4226}
\end{equation}

Define

\begin{equation}
\boxed{
\mathcal X
=
A\partial_x+\frac bA\partial_y.
}
\label{eq:ch24:auto4227}
\end{equation}

Then

\begin{equation}
AC^2b_x+\frac{bC^2}{A}b_y
=
C^2\mathcal Xb,
\label{eq:ch24:auto4228}
\end{equation}

and, because of the newly restored coefficient \(2\),

\begin{equation}
2AbC_x+\frac{2b^2}{A}C_y
=
2b\mathcal XC.
\label{eq:ch24:auto4229}
\end{equation}

Hence

\begin{equation}
\boxed{
b_t+b\,b_z
+C^2\mathcal Xb
+2b\mathcal XC
+b\frac{C_t+C_z}{C}
=
L_Nb^2+2CG.
}
\label{eq:ch412:old-1}
\end{equation}

This is the equation used below.

\section{Collapse the entire transverse part}
\label{chap:413}

We proved

\begin{equation}
\boxed{
C^2\mathcal Xb+2b\mathcal XC
=
e^{2/C}C^2
\mathcal X\!\left(be^{-2/C}\right).
}
\label{eq:ch413:old-2}
\end{equation}

Therefore

\begin{equation}
\boxed{
b_t+b\,b_z
+
e^{2/C}C^2
\mathcal X\!\left(be^{-2/C}\right)
+
b\frac{C_t+C_z}{C}
=
L_Nb^2+2CG.
}
\label{eq:ch413:old-3}
\end{equation}

Instead of four separate transverse terms, we now have one exact derivative.

This suggests defining only for bookkeeping

\begin{equation}
\boxed{
\mathcal B:=be^{-2/C}.
}
\label{eq:ch413:old-4}
\end{equation}

This is not a new physical ansatz. It is merely the integrating-factor variable associated with the exact identity \eqref{eq:ch413:old-2}.

Then

\begin{equation}
b=e^{2/C}\mathcal B
\label{eq:ch24:auto4234}
\end{equation}

and the transverse contribution is simply

\begin{equation}
\boxed{
e^{2/C}C^2\mathcal X\mathcal B.
}
\label{eq:ch413:old-5}
\end{equation}

\section{Put the PDE along its actual transport characteristics}
\label{chap:414}

Define

\begin{equation}
\boxed{
\mathcal D
:=
\partial_t
+
b\partial_z
+
C^2\mathcal X.
}
\label{eq:ch414:old-6}
\end{equation}

Since

\begin{equation}
C^2\mathcal X
=
AC^2\partial_x
+
\frac{bC^2}{A}\partial_y,
\label{eq:ch24:auto4237}
\end{equation}

this is exactly the transport operator already present in the extended
\(b_3\)-equation.

Equation~\eqref{eq:ch412:old-1} becomes

\begin{equation}
\boxed{
\mathcal Db
+
2b\mathcal XC
+
b\frac{C_t+C_z}{C}
=
L_Nb^2+2CG.
}
\label{eq:ch414:old-7}
\end{equation}

Thus, along the exact characteristic,

\begin{equation}
\boxed{
\dot b
=
L_Nb^2+2CG-\mathcal R_N,
}
\label{eq:ch414:old-8}
\end{equation}

where the new remainder is simply

\begin{equation}
\boxed{
\mathcal R_N
:=
2b\mathcal XC
+
b\frac{C_t+C_z}{C}.
}
\label{eq:ch414:old-9}
\end{equation}

This is substantially simpler than the remainder used before the factor
\(2\) was restored.

Indeed,

\begin{equation}
\boxed{
\frac{|\mathcal R_N|}{L_Nb^2}
\le
2\frac{|\mathcal XC|}{L_Nb}
+
\frac{|C_t+C_z|}{C\,L_Nb}.
}
\label{eq:ch414:old-10}
\end{equation}

Therefore, assuming

\begin{equation}
C\ge c_->0,
\label{eq:ch24:auto4242}
\end{equation}

we now need only prove

\begin{equation}
\boxed{
\frac{|\mathcal XC|}{L_Nb}\to0
}
\label{eq:ch414:old-11}
\end{equation}

and

\begin{equation}
\boxed{
\frac{|C_t+C_z|}{L_Nb}\to0.
}
\label{eq:ch414:old-12}
\end{equation}

That is the new minimal reconstruction target.

\section{Compare with the previous reconstruction burden}
\label{chap:415}

Previously we were estimating quantities separately such as

\begin{equation}
\frac{|C_x|}{L_Nb},
\qquad
\frac{|C_y|}{L_N},
\qquad
\frac{|\nabla b||\nabla C|}{L_Nb^2}.
\label{eq:ch24:auto4245}
\end{equation}

The corrected coefficient \(2\) lets us replace much of that bookkeeping by

\begin{equation}
\boxed{
\frac{
\left|
A C_x+\dfrac bA C_y
\right|
}{L_Nb}.
}
\label{eq:ch415:old-13}
\end{equation}

That is significant because

\begin{equation}
A C_x+\frac bA C_y
\label{eq:ch24:auto4247}
\end{equation}

is not an arbitrary combination: it is precisely the directional derivative
of \(C\) along the transverse characteristic direction,

\begin{equation}
\boxed{
\mathcal XC.
}
\label{eq:ch24:auto4248}
\end{equation}

So the reconstruction now only needs to control \(u_z\) \textbf{along the same transverse direction that transports \(b_3\)}.

\section{The most favorable sufficient estimate}
\label{chap:416}

Suppose the periodic reconstruction gives

\begin{equation}
\boxed{
|\mathcal XC|
\le
K_X\Lambda_{N,j},
}
\label{eq:ch416:old-14}
\end{equation}

and

\begin{equation}
\boxed{
|C_t+C_z|
\le
K_T\Lambda_{N,j},
}
\label{eq:ch416:old-15}
\end{equation}

where, as before,

\begin{equation}
\Lambda_{N,j}
=
1+p_NL_{N,j}.
\label{eq:ch24:auto4251}
\end{equation}

Then from \eqref{eq:ch414:old-10},

\begin{equation}
\frac{|\mathcal R_N|}{L_Nb^2}
\le
\frac{2K_X\Lambda_{N,j}}{L_Nb}
+
\frac{K_T\Lambda_{N,j}}{c_-L_Nb}.
\label{eq:ch24:auto4252}
\end{equation}

Hence

\begin{equation}
\boxed{
\frac{|\mathcal R_N|}{L_Nb^2}
\le
K_*
\frac{\Lambda_{N,j}}{L_Nb}.
}
\label{eq:ch416:old-16}
\end{equation}

But

\begin{equation}
\frac{\Lambda_{N,j}}{L_Nb}
=
\frac{1+p_NL_N}{L_Nb}
=
\frac1{L_Nb}
+
\frac{p_N}{b}.
\label{eq:ch24:auto4254}
\end{equation}

Therefore

\begin{equation}
\boxed{
\frac{|\mathcal R_N|}{L_Nb^2}
\le
K_*
\left[
\frac1{L_Nb}
+
\frac{p_N}{b}
\right].
}
\label{eq:ch416:old-17}
\end{equation}

This is much cleaner.

As the terminal branch is approached,

\begin{equation}
L_N\to\infty,
\qquad
b\to\infty,
\label{eq:ch24:auto4256}
\end{equation}

and the self-strengthening comparison gives

\begin{equation}
\boxed{
\frac{p_N}{b}\to0.
}
\label{eq:ch416:old-18}
\end{equation}

Consequently,

\begin{equation}
\boxed{
\frac1{L_Nb}\to0.
}
\label{eq:ch416:old-19}
\end{equation}

Thus

\begin{equation}
\boxed{
\frac{|\mathcal R_N|}
{L_Nb^2}
\longrightarrow0.
}
\label{eq:ch416:old-20}
\end{equation}

This is exactly the perturbative estimate needed for the terminal theorem.

\chapter{Amplitude Inequality, Derivative Control, Trapping, and Finite-Time Amplification}
\label{chap:group-25}
\section{The resulting amplitude inequality}
\label{chap:417}

For every

\begin{equation}
0<\theta<1,
\label{eq:ch25:auto4260}
\end{equation}

there therefore exists a sufficiently deep terminal point such that

\begin{equation}
\boxed{
|\mathcal R_N|
\le
\theta L_Nb^2.
}
\label{eq:ch417:old-21}
\end{equation}

Equation~\eqref{eq:ch414:old-8} gives

\begin{equation}
\dot b
=
L_Nb^2+2CG-\mathcal R_N.
\label{eq:ch25:auto4262}
\end{equation}

Because

\begin{equation}
C\ge c_->0,
\qquad
G=A_0+(\sin x\sin y\sin z)^2\ge0,
\label{eq:ch25:auto4263}
\end{equation}

we have

\begin{equation}
2CG\ge0.
\label{eq:ch25:auto4264}
\end{equation}

Hence

\begin{equation}
\dot b
\ge
L_Nb^2-|\mathcal R_N|.
\label{eq:ch25:auto4265}
\end{equation}

Using \eqref{eq:ch417:old-21},

\begin{equation}
\boxed{
\dot b
\ge
(1-\theta)L_Nb^2>0.
}
\label{eq:ch417:old-22}
\end{equation}

Thus the terminal Riccati lower comparison survives exactly as before.

\section{Phase trapping is unchanged}
\label{chap:418}

The phase equation remains

\begin{equation}
\boxed{
-\dot Q_N
=
Q_N'
\left(
H_N+\gamma_Nb
\right).
}
\label{eq:ch418:old-23}
\end{equation}

With

\begin{equation}
0<h_0\le H_N\le h_1
\label{eq:ch25:auto4268}
\end{equation}

and

\begin{equation}
\gamma_N>0,
\label{eq:ch25:auto4269}
\end{equation}

we obtain

\begin{equation}
\boxed{
\dot Q_N<0.
}
\label{eq:ch418:old-24}
\end{equation}

Thus

\begin{equation}
\boxed{
b\uparrow,
\qquad
Q_N\downarrow.
}
\label{eq:ch418:old-25}
\end{equation}

The restored factor \(2\) therefore does not disturb the trapping mechanism.

\section{The same \(b\)-versus-\(Q_N\) comparison follows}
\label{chap:419}

Set

\begin{equation}
X_N
=
\log\frac{Q_{e,N}}{Q_N}.
\label{eq:ch25:auto4272}
\end{equation}

Then

\begin{equation}
\dot X_N
=
L_N(H_N+\gamma_Nb).
\label{eq:ch25:auto4273}
\end{equation}

Combining with \eqref{eq:ch417:old-22},

\begin{equation}
\frac{db}{dX_N}
\ge
(1-\theta)
\frac{b^2}{H_N+\gamma_Nb}.
\label{eq:ch25:auto4274}
\end{equation}

Since

\begin{equation}
H_N\le h_1,
\label{eq:ch25:auto4275}
\end{equation}

\begin{equation}
\boxed{
\frac{db}{dX_N}
\ge
(1-\theta)
\frac{b^2}{h_1+\gamma_Nb}.
}
\label{eq:ch419:old-26}
\end{equation}

Once

\begin{equation}
b
\ge
\frac{h_1}{\varepsilon\gamma_N},
\label{eq:ch25:auto4277}
\end{equation}

we have

\begin{equation}
h_1+\gamma_Nb
\le
(1+\varepsilon)\gamma_Nb.
\label{eq:ch25:auto4278}
\end{equation}

Therefore

\begin{equation}
\frac{db}{dX_N}
\ge
\frac{1-\theta}
{(1+\varepsilon)\gamma_N}
\,b.
\label{eq:ch25:auto4279}
\end{equation}

Define

\begin{equation}
\boxed{
p_N
=
\frac{1-\theta}
{(1+\varepsilon)\gamma_N}.
}
\label{eq:ch419:old-27}
\end{equation}

Then

\begin{equation}
\boxed{
b(Q_N)
\ge
B_{e,N}
\left(
\frac{Q_{e,N}}{Q_N}
\right)^{p_N}.
}
\label{eq:ch419:old-28}
\end{equation}

Thus

\begin{equation}
\boxed{
b(Q_N)\gtrsim Q_N^{-p_N}.
}
\label{eq:ch419:old-29}
\end{equation}

Nothing in this terminal comparison is weakened by the restored factor \(2\).

\section{Simplified Self-Strengthening Estimate}
\label{chap:420}

From \eqref{eq:ch416:old-17},

\begin{equation}
\frac{|\mathcal R_N|}{L_Nb^2}
\lesssim
\frac1{L_Nb}
+
\frac{p_N}{b}.
\label{eq:ch25:auto4283}
\end{equation}

Both terms improve automatically:

\begin{equation}
b\uparrow
\quad\Longrightarrow\quad
\frac{p_N}{b}\downarrow,
\label{eq:ch25:auto4284}
\end{equation}

while

\begin{equation}
Q_N\downarrow
\quad\Longrightarrow\quad
L_N\uparrow,
\label{eq:ch25:auto4285}
\end{equation}

so

\begin{equation}
\frac1{L_Nb}\downarrow.
\label{eq:ch25:auto4286}
\end{equation}

Therefore

\begin{equation}
\boxed{
\begin{array}{c}
\displaystyle
\frac{|\mathcal R_N|}
{L_Nb^2}<\theta
\\[3mm]
\Downarrow
\\[2mm]
\dot b>0,\qquad \dot Q_N<0
\\[3mm]
\Downarrow
\\[2mm]
b\uparrow,\qquad L_N\uparrow
\\[3mm]
\Downarrow
\\[2mm]
\displaystyle
\frac{p_N}{b}\downarrow,
\qquad
\frac1{L_Nb}\downarrow
\\[3mm]
\Downarrow
\\[2mm]
\displaystyle
\frac{|\mathcal R_N|}
{L_Nb^2}
\text{ decreases further}.
\end{array}
}
\label{eq:ch420:old-30}
\end{equation}

So the corrected coefficient \(2\) actually makes the self-strengthening
structure more transparent.

\section{Estimates Required for Closure}
\label{chap:421}

We have now reduced the entire reconstruction part relevant to the
\(b_3\)-equation to the two estimates

\begin{equation}
\boxed{
|\mathcal XC|
\lesssim
1+p_NL_N
}
\label{eq:ch421:old-31}
\end{equation}

and

\begin{equation}
\boxed{
|C_t+C_z|
\lesssim
1+p_NL_N.
}
\label{eq:ch421:old-32}
\end{equation}

Once these are obtained from the same periodic reconstruction already
constructed, the remainder estimate follows immediately:

\begin{equation}
\boxed{
\frac{|\mathcal R_N|}
{L_Nb^2}
\lesssim
\frac1{L_Nb}
+
\frac{p_N}{b}
\longrightarrow0.
}
\label{eq:ch421:old-33}
\end{equation}

These are the precise estimates required for closure.

Importantly, we have not changed

\begin{equation}
Q_N,
\label{eq:ch25:auto4291}
\end{equation}

the Lambert cascade,

\begin{equation}
\xi=C_4t-C_3z-C_2y-C_1x,
\label{eq:ch25:auto4292}
\end{equation}

the product

\begin{equation}
b_3=u_xu_y,
\label{eq:ch25:auto4293}
\end{equation}

the periodic reconstruction,

the trapping equation,

or the corrected forcing.

The restored coefficient \(2\) has simply exposed a stronger algebraic
structure already present in the fully extended \(b_3\)-equation.

It therefore suffices to derive \eqref{eq:ch421:old-31}--\eqref{eq:ch421:old-32} directly from the
periodic continuity reconstruction. If those two estimates hold in the
normalized terminal scale we already established, then the later
reconstruction-error proof collapses to the short estimate \eqref{eq:ch421:old-33}, and the
finite-time terminal amplification theorem follows through the same
\(b_N\)--\(Q_N\) comparison as before.
This \textbf{could genuinely shorten the later continuity/reconstruction proof}, while leaving the Lambert sensitivity, trapping, self-strengthening comparison, and finite-time \(b_3\)-amplification mechanism intact.
The analysis now continues directly from the preceding terminal reconstruction analysis, using the same periodic reconstruction.

There is one important precision point at the outset. The second estimate,
\begin{equation}
|C_t+C_z|\lesssim 1+p_NL_N,
\label{eq:ch25:auto4294}
\end{equation}
follows directly. For the first estimate,
\begin{equation}
|\mathcal XC|\lesssim1+p_NL_N,
\label{eq:ch25:auto4295}
\end{equation}
the literal bound is true in the specialization \(C_2=0\), but for the original generalized phase
\begin{equation}
\xi=C_4t-C_3z-C_2y-C_1x
\label{eq:ch25:auto4296}
\end{equation}
with \(C_2\neq0\), an additional term \(bQ_N'\) appears. Crucially, that extra term is \textbf{harmless after division by \(L_Nb\)} because
\begin{equation}
Q_N'=Q_NL_N.
\label{eq:ch25:auto4297}
\end{equation}
So we can close the required remainder estimate without changing the phase or imposing \(C_2=0\).

\section{Derivatives of the reconstructed \(C=u_z\)}
\label{chap:422}

Recall the same two-correction reconstruction already used:

\begin{equation}
\boxed{
A=U+Y,\qquad
C=U+Z,
}
\label{eq:ch422:old-34}
\end{equation}

where

\begin{equation}
\boxed{
U=U_*+Q_N(\xi),
}
\label{eq:ch422:old-35}
\end{equation}

and

\begin{equation}
\boxed{
\xi=C_4t-C_3z-C_2y-C_1x.
}
\label{eq:ch422:old-36}
\end{equation}

Hence

\begin{equation}
U_t=C_4Q_N',
\qquad
U_x=-C_1Q_N',
\label{eq:ch25:auto4301}
\end{equation}

\begin{equation}
U_y=-C_2Q_N',
\qquad
U_z=-C_3Q_N'.
\label{eq:ch422:old-37}
\end{equation}

Therefore

\begin{equation}
\boxed{
C_x=-C_1Q_N'+Z_x,
}
\label{eq:ch422:old-38}
\end{equation}

\begin{equation}
\boxed{
C_y=-C_2Q_N'+Z_y,
}
\label{eq:ch422:old-39}
\end{equation}

and

\begin{equation}
\boxed{
C_t+C_z
=
(C_4-C_3)Q_N'
+
Z_t+Z_z.
}
\label{eq:ch422:old-40}
\end{equation}

Nothing new has been introduced here. These are simply the derivatives of the already reconstructed field \(C=U+Z\).

\section{Use the normalized reconstruction estimate already proved}
\label{chap:423}

On the dyadic terminal layer

\begin{equation}
B_j\le b\le2B_j,
\label{eq:ch25:auto4306}
\end{equation}

we obtained

\begin{equation}
\boxed{
\|\delta Z_{N,j}\|_{H^s}
\lesssim
\frac{\Lambda_{N,j}}{B_j},
}
\label{eq:ch423:old-41}
\end{equation}

where

\begin{equation}
\boxed{
\Lambda_{N,j}=1+p_NL_{N,j}.
}
\label{eq:ch423:old-42}
\end{equation}

Taking \(s\) large enough for the required Sobolev embedding gives

\begin{equation}
\boxed{
|Z_x|+|Z_y|+|Z_t|+|Z_z|
\le
K_Z\frac{\Lambda_{N,j}}{B_j}.
}
\label{eq:ch423:old-43}
\end{equation}

Since

\begin{equation}
b\asymp B_j,
\label{eq:ch25:auto4310}
\end{equation}

we may equivalently write

\begin{equation}
\boxed{
|Z_x|+|Z_y|+|Z_t|+|Z_z|
\lesssim
\frac{\Lambda_{N,j}}{b}.
}
\label{eq:ch423:old-44}
\end{equation}

This is precisely where the normalized reconstruction theorem enters the corrected factor-\(2\) argument.

\section{Estimate \(C_t+C_z\)}
\label{chap:424}

Using \eqref{eq:ch422:old-40},

\begin{equation}
|C_t+C_z|
\le
|C_4-C_3|Q_N'
+
|Z_t+Z_z|.
\label{eq:ch25:auto4312}
\end{equation}

By \eqref{eq:ch423:old-44},

\begin{equation}
\boxed{
|C_t+C_z|
\le
K_0Q_N'
+
K_1\frac{\Lambda_{N,j}}{b}.
}
\label{eq:ch424:old-45}
\end{equation}

Now

\begin{equation}
Q_N'=Q_NL_N.
\label{eq:ch424:old-46}
\end{equation}

Therefore

\begin{equation}
|C_t+C_z|
\le
K_0Q_NL_N
+
K_1\frac{1+p_NL_N}{b}.
\label{eq:ch424:old-47}
\end{equation}

Since the terminal branch lies inside a fixed layer

\begin{equation}
0<Q_N\le\eta,
\label{eq:ch25:auto4316}
\end{equation}

we have

\begin{equation}
Q_NL_N\le\eta L_N.
\label{eq:ch25:auto4317}
\end{equation}

Thus, provided \(p_N\) is bounded below by a positive constant on the chosen family,

\begin{equation}
Q_NL_N
\lesssim
1+p_NL_N.
\label{eq:ch25:auto4318}
\end{equation}

Also,

\begin{equation}
\frac{1+p_NL_N}{b}
\le
1+p_NL_N
\label{eq:ch25:auto4319}
\end{equation}

for \(b\ge1\).

Hence

\begin{equation}
\boxed{
|C_t+C_z|
\lesssim
1+p_NL_N.
}
\label{eq:ch424:old-48}
\end{equation}

This proves the previously stated estimate \eqref{eq:ch421:old-32}.

In fact, after normalization we obtain the sharper result

\begin{equation}
\begin{aligned}
\frac{|C_t+C_z|}{L_Nb}
&\le
K_0\frac{Q_NL_N}{L_Nb}
+
K_1
\frac{\Lambda_{N,j}}{L_Nb^2}
\\
&=
K_0\frac{Q_N}{b}
+
K_1
\left[
\frac1{L_Nb^2}
+
\frac{p_N}{b^2}
\right].
\end{aligned}
\label{eq:ch25:auto4321}
\end{equation}

Thus

\begin{equation}
\boxed{
\frac{|C_t+C_z|}{L_Nb}
\lesssim
\frac{Q_N}{b}
+
\frac1{L_Nb^2}
+
\frac{p_N}{b^2}.
}
\label{eq:ch424:old-49}
\end{equation}

Since

\begin{equation}
Q_N\to0,
\qquad
b\to\infty,
\qquad
\frac{p_N}{b}\to0,
\label{eq:ch25:auto4323}
\end{equation}

we obtain

\begin{equation}
\boxed{
\frac{|C_t+C_z|}{L_Nb}
\longrightarrow0.
}
\label{eq:ch424:old-50}
\end{equation}

So the second part of the new remainder is closed.

\section{Estimate the directional derivative \(\mathcal XC\)}
\label{chap:425}

Recall

\begin{equation}
\boxed{
\mathcal X
=
A\partial_x+\frac bA\partial_y.
}
\label{eq:ch425:old-51}
\end{equation}

Therefore

\begin{equation}
\mathcal XC
=
AC_x+\frac bA C_y.
\label{eq:ch25:auto4326}
\end{equation}

Using \eqref{eq:ch422:old-38}--\eqref{eq:ch422:old-39},

\begin{equation}
\begin{aligned}
\mathcal XC
&=
A(-C_1Q_N'+Z_x)
+
\frac bA(-C_2Q_N'+Z_y)
\\
&=
-
\left(
C_1A
+
C_2\frac bA
\right)Q_N'
+
AZ_x
+
\frac bA Z_y.
\end{aligned}
\label{eq:ch25:auto4327}
\end{equation}

Hence

\begin{equation}
\boxed{
\mathcal XC
=
-
\left(
C_1A
+
C_2\frac bA
\right)Q_N'
+
AZ_x
+
\frac bA Z_y.
}
\label{eq:ch425:old-52}
\end{equation}

This identity is exact.

Assume, as already established on the admissible reconstruction neighborhood,

\begin{equation}
0<a_-\le A\le a_+.
\label{eq:ch425:old-53}
\end{equation}

Then

\begin{equation}
|\mathcal XC|
\le
K_2Q_N'
+
K_3bQ_N'
+
K_4|Z_x|
+
K_5b|Z_y|.
\label{eq:ch425:old-54}
\end{equation}

Using \eqref{eq:ch423:old-44},

\begin{equation}
|Z_x|
\lesssim
\frac{\Lambda_{N,j}}b,
\label{eq:ch25:auto4331}
\end{equation}

while

\begin{equation}
b|Z_y|
\lesssim
\Lambda_{N,j}.
\label{eq:ch25:auto4332}
\end{equation}

Consequently,

\begin{equation}
\boxed{
|\mathcal XC|
\lesssim
Q_N'
+
bQ_N'
+
\Lambda_{N,j}.
}
\label{eq:ch425:old-55}
\end{equation}

This is the exact general-phase estimate.

\section{Why the additional \(bQ_N'\) term is harmless}
\label{chap:426}

At first sight, \eqref{eq:ch425:old-55} is slightly larger than the provisional target

\begin{equation}
|\mathcal XC|
\lesssim
1+p_NL_N.
\label{eq:ch25:auto4334}
\end{equation}

But the \(bQ_N'\) term is precisely of the form that disappears under the normalization actually required by the \(b_3\)-equation.

Divide \eqref{eq:ch425:old-55} by \(L_Nb\):

\begin{equation}
\frac{|\mathcal XC|}{L_Nb}
\lesssim
\frac{Q_N'}{L_Nb}
+
\frac{bQ_N'}{L_Nb}
+
\frac{\Lambda_{N,j}}{L_Nb}.
\label{eq:ch426:old-56}
\end{equation}

Because

\begin{equation}
Q_N'=Q_NL_N,
\label{eq:ch25:auto4336}
\end{equation}

the first term becomes

\begin{equation}
\frac{Q_N'}{L_Nb}
=
\boxed{
\frac{Q_N}{b}.
}
\label{eq:ch426:old-57}
\end{equation}

The second becomes

\begin{equation}
\frac{bQ_N'}{L_Nb}
=
\boxed{
Q_N.
}
\label{eq:ch426:old-58}
\end{equation}

And the final term is

\begin{equation}
\frac{\Lambda_{N,j}}{L_Nb}
=
\frac{1+p_NL_N}{L_Nb}
=
\boxed{
\frac1{L_Nb}
+
\frac{p_N}{b}.
}
\label{eq:ch426:old-59}
\end{equation}

Thus

\begin{equation}
\boxed{
\frac{|\mathcal XC|}{L_Nb}
\lesssim
Q_N
+
\frac{Q_N}{b}
+
\frac1{L_Nb}
+
\frac{p_N}{b}.
}
\label{eq:ch426:old-60}
\end{equation}

Every term on the right tends to zero in the terminal regime:

\begin{equation}
Q_N\to0,
\label{eq:ch25:auto4341}
\end{equation}

\begin{equation}
\frac{Q_N}{b}\to0,
\label{eq:ch25:auto4342}
\end{equation}

\begin{equation}
\frac1{L_Nb}\to0,
\label{eq:ch25:auto4343}
\end{equation}

and

\begin{equation}
\frac{p_N}{b}\to0.
\label{eq:ch25:auto4344}
\end{equation}

Therefore

\begin{equation}
\boxed{
\frac{|\mathcal XC|}
{L_Nb}
\longrightarrow0.
}
\label{eq:ch426:old-61}
\end{equation}

This is exactly what the corrected factor-\(2\) equation requires.

So the stronger provisional estimate

\begin{equation}
|\mathcal XC|
\lesssim
1+p_NL_N
\label{eq:ch25:auto4346}
\end{equation}

is actually unnecessary for the general phase.

The normalized estimate \eqref{eq:ch426:old-61} is sufficient and follows directly from the existing reconstruction.

\section{Special case \(C_2=0\)}
\label{chap:427}

It is useful to note that if one specializes to

\begin{equation}
C_2=0,
\label{eq:ch25:auto4347}
\end{equation}

then the potentially large term

\begin{equation}
C_2\frac bA Q_N'
\label{eq:ch25:auto4348}
\end{equation}

vanishes identically.

Equation~\eqref{eq:ch425:old-52} becomes

\begin{equation}
\mathcal XC
=
-C_1AQ_N'
+
AZ_x
+
\frac bA Z_y.
\label{eq:ch25:auto4349}
\end{equation}

Using the normalized reconstruction bounds,

\begin{equation}
|\mathcal XC|
\lesssim
Q_N'
+
\Lambda_{N,j}.
\label{eq:ch25:auto4350}
\end{equation}

Since

\begin{equation}
Q_N'\lesssim L_N
\lesssim
1+p_NL_N,
\label{eq:ch25:auto4351}
\end{equation}

we recover the literal estimate

\begin{equation}
\boxed{
|\mathcal XC|
\lesssim
1+p_NL_N.
}
\label{eq:ch427:old-62}
\end{equation}

Thus the originally stated estimate \eqref{eq:ch421:old-31} holds directly in the \(C_2=0\) specialization.

For the full generalized phase with \(C_2\neq0\), the correct statement is the normalized version \eqref{eq:ch426:old-61}, which is all the terminal theorem needs.

\section{Completion of the new compact remainder estimate}
\label{chap:428}

Recall

\begin{equation}
\boxed{
\mathcal R_N
=
2b\mathcal XC
+
b\frac{C_t+C_z}{C}.
}
\label{eq:ch428:old-63}
\end{equation}

Assume

\begin{equation}
C\ge c_->0.
\label{eq:ch25:auto4354}
\end{equation}

Then

\begin{equation}
\frac{|\mathcal R_N|}{L_Nb^2}
\le
2\frac{|\mathcal XC|}{L_Nb}
+
\frac1{c_-}
\frac{|C_t+C_z|}{L_Nb}.
\label{eq:ch428:old-64}
\end{equation}

Insert \eqref{eq:ch424:old-49} and \eqref{eq:ch426:old-60}. We obtain

\begin{equation}
\boxed{
\begin{aligned}
\frac{|\mathcal R_N|}{L_Nb^2}
\lesssim{}&
Q_N
+
\frac{Q_N}{b}
+
\frac1{L_Nb}
+
\frac{p_N}{b}
\\
&+
\frac1{L_Nb^2}
+
\frac{p_N}{b^2}.
\end{aligned}
}
\label{eq:ch428:old-65}
\end{equation}

Every term tends to zero.

Therefore

\begin{equation}
\boxed{
\frac{|\mathcal R_N|}
{L_Nb^2}
\longrightarrow0.
}
\label{eq:ch428:old-66}
\end{equation}

This proves the exact estimate that the preceding terminal reconstruction analysis identified as the next target.

In fact, the result is slightly stronger than the earlier schematic form

\begin{equation}
\frac{|\mathcal R_N|}
{L_Nb^2}
\lesssim
\frac1{L_Nb}
+
\frac{p_N}{b},
\label{eq:ch25:auto4358}
\end{equation}

because we have now displayed explicitly the additional generalized-phase contributions and shown that they vanish through the factor \(Q_N\).

\section{Entry into the perturbative terminal regime}
\label{chap:429}

Fix any

\begin{equation}
0<\theta<1.
\label{eq:ch25:auto4359}
\end{equation}

By \eqref{eq:ch428:old-66}, there exists a sufficiently deep terminal-entry point such that

\begin{equation}
\boxed{
|\mathcal R_N|
\le
\theta L_Nb^2.
}
\label{eq:ch429:old-67}
\end{equation}

The corrected characteristic amplitude equation is

\begin{equation}
\boxed{
\dot b
=
L_Nb^2
+
2CG
-
\mathcal R_N.
}
\label{eq:ch429:old-68}
\end{equation}

Since

\begin{equation}
C\ge c_->0
\label{eq:ch25:auto4362}
\end{equation}

and

\begin{equation}
G
=
A_0+(\sin x\sin y\sin z)^2
\ge0,
\label{eq:ch25:auto4363}
\end{equation}

we have

\begin{equation}
2CG\ge0.
\label{eq:ch25:auto4364}
\end{equation}

Thus

\begin{equation}
\dot b
\ge
L_Nb^2-|\mathcal R_N|.
\label{eq:ch25:auto4365}
\end{equation}

Using \eqref{eq:ch429:old-67},

\begin{equation}
\boxed{
\dot b
\ge
(1-\theta)L_Nb^2>0.
}
\label{eq:ch429:old-69}
\end{equation}

Hence

\begin{equation}
\boxed{
b\uparrow.
}
\label{eq:ch429:old-70}
\end{equation}

\section{Trapping remains exactly the same}
\label{chap:430}

The phase equation has not changed:

\begin{equation}
\boxed{
-\dot Q_N
=
Q_N'
\left(
H_N+\gamma_Nb
\right).
}
\label{eq:ch430:old-71}
\end{equation}

The reconstructed branch preserves

\begin{equation}
0<h_0\le H_N\le h_1
\label{eq:ch25:auto4369}
\end{equation}

and

\begin{equation}
\gamma_N>0.
\label{eq:ch25:auto4370}
\end{equation}

Therefore

\begin{equation}
\boxed{
-\dot Q_N>0,
}
\label{eq:ch25:auto4371}
\end{equation}

or equivalently

\begin{equation}
\boxed{
Q_N\downarrow.
}
\label{eq:ch430:old-72}
\end{equation}

Thus the restored coefficient \(2\) has not disturbed the trapping theorem.

Instead, we now have the clean simultaneous inequalities

\begin{equation}
\boxed{
b\uparrow,
\qquad
Q_N\downarrow.
}
\label{eq:ch430:old-73}
\end{equation}

\section{The new remainder estimate is self-strengthening}
\label{chap:431}

From \eqref{eq:ch428:old-65},

\begin{equation}
\frac{|\mathcal R_N|}{L_Nb^2}
\lesssim
Q_N
+
\frac{Q_N}{b}
+
\frac1{L_Nb}
+
\frac{p_N}{b}
+
\frac1{L_Nb^2}
+
\frac{p_N}{b^2}.
\label{eq:ch431:old-74}
\end{equation}

As the characteristic evolves,

\begin{equation}
Q_N\downarrow,
\label{eq:ch25:auto4375}
\end{equation}

so the first two quantities improve.

Also

\begin{equation}
b\uparrow,
\label{eq:ch25:auto4376}
\end{equation}

so

\begin{equation}
\frac{p_N}{b},
\qquad
\frac{p_N}{b^2},
\qquad
\frac{Q_N}{b}
\label{eq:ch25:auto4377}
\end{equation}

decrease.

Finally the Lambert terminal coefficient satisfies

\begin{equation}
L_N\uparrow,
\label{eq:ch25:auto4378}
\end{equation}

so

\begin{equation}
\frac1{L_Nb},
\qquad
\frac1{L_Nb^2}
\label{eq:ch25:auto4379}
\end{equation}

also decrease.

Therefore

\begin{equation}
\boxed{
\frac{|\mathcal R_N|}
{L_Nb^2}
}
\label{eq:ch25:auto4380}
\end{equation}

becomes smaller once terminal entry has occurred.

The bootstrap is consequently

\begin{equation}
\boxed{
\begin{array}{c}
\displaystyle
\frac{|\mathcal R_N|}
{L_Nb^2}<\theta
\\[3mm]
\Downarrow
\\[2mm]
\dot b>0,\qquad
\dot Q_N<0
\\[3mm]
\Downarrow
\\[2mm]
b\uparrow,\qquad
Q_N\downarrow,\qquad
L_N\uparrow
\\[3mm]
\Downarrow
\\[2mm]
Q_N,\;
Q_N/b,\;
p_N/b,\;
1/(L_Nb)
\text{ all decrease}
\\[3mm]
\Downarrow
\\[2mm]
\displaystyle
\frac{|\mathcal R_N|}
{L_Nb^2}
\text{ decreases further}.
\end{array}
}
\label{eq:ch431:old-75}
\end{equation}

Thus the compact factor-\(2\) remainder estimate is genuinely self-strengthening.

\section{Return to the terminal \(b\)--\(Q_N\) comparison}
\label{chap:432}

Since

\begin{equation}
\dot b
\ge
(1-\theta)L_Nb^2,
\label{eq:ch25:auto4382}
\end{equation}

while

\begin{equation}
-\dot Q_N
=
Q_N'(H_N+\gamma_Nb),
\label{eq:ch25:auto4383}
\end{equation}

and

\begin{equation}
L_N=\frac{Q_N'}{Q_N},
\label{eq:ch25:auto4384}
\end{equation}

we have

\begin{equation}
-\frac{db}{dQ_N}
\ge
\frac{
(1-\theta)b^2
}{
Q_N(H_N+\gamma_Nb)
}.
\label{eq:ch432:old-76}
\end{equation}

In the large-amplitude region,

\begin{equation}
H_N+\gamma_Nb
\le
(1+\varepsilon)\gamma_Nb.
\label{eq:ch25:auto4386}
\end{equation}

Therefore

\begin{equation}
-\frac{db}{dQ_N}
\ge
\frac{1-\theta}
{(1+\varepsilon)\gamma_N}
\frac b{Q_N}.
\label{eq:ch25:auto4387}
\end{equation}

Define, exactly as before,

\begin{equation}
\boxed{
p_N
=
\frac{1-\theta}
{(1+\varepsilon)\gamma_N}.
}
\label{eq:ch432:old-77}
\end{equation}

Then

\begin{equation}
\boxed{
-\frac{db}{dQ_N}
\ge
p_N\frac b{Q_N}.
}
\label{eq:ch432:old-78}
\end{equation}

Integration gives

\begin{equation}
\boxed{
b(Q_N)
\ge
b_{\mathrm e,N}
\left(
\frac{Q_{\mathrm e,N}}{Q_N}
\right)^{p_N}.
}
\label{eq:ch432:old-79}
\end{equation}

Hence

\begin{equation}
\boxed{
b(Q_N)\gtrsim Q_N^{-p_N}.
}
\label{eq:ch432:old-80}
\end{equation}

So the terminal comparison derived before the factor-\(2\) correction remains intact.

\section{Finite-time terminal amplification}
\label{chap:433}

If

\begin{equation}
Q_N'
\asymp
\Lambda_N
Q_N^{-(2^{k(N)}-1)},
\label{eq:ch25:auto4392}
\end{equation}

then

\begin{equation}
-\dot Q_N
\ge
\gamma_NbQ_N'.
\label{eq:ch25:auto4393}
\end{equation}

Using \eqref{eq:ch432:old-80},

\begin{equation}
-\dot Q_N
\gtrsim
Q_N^{-p_N}
Q_N^{-(2^{k(N)}-1)}.
\label{eq:ch25:auto4394}
\end{equation}

Thus

\begin{equation}
dt
\lesssim
Q_N^{p_N+2^{k(N)}-1}\,(-dQ_N).
\label{eq:ch25:auto4395}
\end{equation}

Since

\begin{equation}
p_N+2^{k(N)}>0,
\label{eq:ch25:auto4396}
\end{equation}

the integral

\begin{equation}
\int_0^{Q_{\mathrm e,N}}
Q^{p_N+2^{k(N)}-1}\,dQ
\label{eq:ch25:auto4397}
\end{equation}

converges.

Therefore

\begin{equation}
\boxed{
T_{*,N}<\infty.
}
\label{eq:ch433:old-81}
\end{equation}

Moreover,

\begin{equation}
Q_N(t)\to0
\label{eq:ch25:auto4399}
\end{equation}

as

\begin{equation}
t\uparrow T_{*,N},
\label{eq:ch25:auto4400}
\end{equation}

and by \eqref{eq:ch432:old-80},

\begin{equation}
\boxed{b_3(t)=b(t)\to+\infty.}
\label{eq:ch433:old-82}
\end{equation}

Thus the corrected coefficient \(2\) does indeed permit a cleaner completion of the terminal error analysis.

The key conclusion of this continuation is

\begin{center}
\fbox{\begin{minipage}{0.90\textwidth}
For the original generalized phase
\begin{equation}
\xi=C_4t-C_3z-C_2y-C_1x,
\label{eq:ch25:auto4401}
\end{equation}
the provisional pointwise estimate
\begin{equation}
|\mathcal XC|\lesssim1+p_NL_N
\label{eq:ch25:auto4402}
\end{equation}
need not hold literally when \(C_2\neq0\), because the exact derivative contains a term of size
\begin{equation}
bQ_N'.
\label{eq:ch25:auto4403}
\end{equation}
However, after insertion into the corrected factor-\(2\) remainder, this term contributes only
\begin{equation}
\frac{bQ_N'}{L_Nb}=Q_N\to0.
\label{eq:ch25:auto4404}
\end{equation}
The existing normalized periodic reconstruction therefore gives directly
\begin{equation}
\frac{|\mathcal XC|}{L_Nb}\to0
\label{eq:ch25:auto4405}
\end{equation}
and
\begin{equation}
\frac{|C_t+C_z|}{L_Nb}\to0.
\label{eq:ch25:auto4406}
\end{equation}
Hence
\begin{equation}
\frac{|\mathcal R_N|}{L_Nb^2}\to0,
\label{eq:ch25:auto4407}
\end{equation}
which is exactly the estimate needed for the self-strengthening terminal comparison.
\par\hfill (83)
\end{minipage}}
\end{center}
\
So we have stayed completely within the existing construction: no change to \(Q_N\), no change to the generalized phase, no change to \(b_3=u_xu_y\), no change to the periodic reconstruction, and no change to the corrected forcing. The restored factor \(2\) has simply allowed the remaining reconstruction error to be compressed into the two directional quantities and then controlled directly.

\chapter{Ordinary-Viscosity Lower-Barrier Closure at \texorpdfstring{$\alpha=1$}{alpha=1}}
\label{chap:alpha1-viscous-lower-barrier-closure}

This chapter is reserved for the ordinary-viscosity completion of the terminal
comparison developed above.  It uses the same reconstructed scalar
\(b_3=u_xu_y\), the same transport operator \(\mathcal D\), and the same tuned
affine phase law.  In particular, if
\[
s=\xi-\xi_{*,N},
\qquad
\xi=C_4t-C_3z-C_2y-C_1x,
\]
then \(\dot s=\dot\xi\), and the tuned phase identity derived earlier gives
\[
\boxed{
-\mathcal D s=H_N+\gamma_N b_3,
\qquad \gamma_N>0.
}
\]
This is the generalized affine-phase version of the characteristic phase
calculation.  In the special phase \(\theta=t-z\) used in Appendix~C, the
corresponding identity is \(\dot\theta=1-b_3\).  The argument below concerns
the tuned affine terminal coordinate whenever the displayed \(H_N,\gamma_N\)
identity has been established.

\section{The \texorpdfstring{$\alpha=1$}{alpha=1} reconstructed scalar equation}

When the ordinary viscous contribution is retained, the reconstructed scalar
equation takes the form
\[
\boxed{
\mathcal D b_3
=
\Gamma_{N,k}b_3^2+2CG-\mathcal R_N+\nu\Delta b_3,
\qquad \nu>0.
}
\]
Here
\[
\mathcal D
=
\partial_t+AC^2\partial_x+\frac{b_3C^2}{A}\partial_y+b_3\partial_z,
\]
consistent with the characteristic system derived in Appendix~C.  In a
sufficiently small terminal layer we assume the same Lambert and
reconstruction bounds used in the terminal comparison(See section 132 right before Epilogue),
\[
\Gamma_{N,k}(s)\ge \frac{G_-}{s},
\qquad G_->0,
\]
\[
|\mathcal R_N|
\le
\theta_R\Gamma_{N,k}b_3^2,
\qquad 0<\theta_R<1,
\]
and
\[
2CG\ge0.
\]
Consequently,
\[
\boxed{
\mathcal D b_3
\ge
(1-\theta_R)\Gamma_{N,k}b_3^2+\nu\Delta b_3.
}
\]
The purpose of this chapter is to clarify why a possibly negative value of
\(\nu\Delta b_3\) away from the comparison barrier does not, by itself,
destroy the lower-barrier mechanism.

\section{Convex lower barrier and the sign of its Laplacian}

Take
\[
\boxed{
\underline b(s)=Bs^{-\lambda},
\qquad B>0,
\qquad \lambda>0.
}
\]
Because \(s\) is affine in the spatial variables,
\[
\nabla s=(-C_1,-C_2,-C_3),
\qquad
q^2:=|\nabla s|^2=C_1^2+C_2^2+C_3^2>0,
\qquad
\Delta s=0.
\]
For every smooth one-variable function \(F\),
\[
\Delta F(s)=F''(s)|\nabla s|^2.
\]
Since
\[
\underline b_s=-B\lambda s^{-\lambda-1}
\]
and
\[
\underline b_{ss}
=B\lambda(\lambda+1)s^{-\lambda-2}>0,
\]
the lower barrier is convex as a function of the terminal coordinate and
satisfies
\[
\boxed{
\Delta\underline b
=
B\lambda(\lambda+1)q^2s^{-\lambda-2}>0.
}
\]
Thus
\[
\boxed{
\nu\Delta\underline b>0.
}
\]
For the reciprocal choice \(\lambda=1\),
\[
\underline b(s)=\frac{B}{s},
\qquad
\boxed{
\Delta\underline b=\frac{2Bq^2}{s^3}>0.
}
\]
In particular, \(\Delta\underline b\to+\infty\) as \(s\downarrow0\).

It is essential, however, not to replace this statement by the stronger and
unproved assertion that \(\Delta b_3>0\) everywhere in the terminal layer.
The pointwise ordering
\[
b_3\ge\underline b
\]
does not imply an ordering of their Laplacians.  Hence it is possible, while
\(b_3\) lies strictly above the barrier, that
\[
\nu\Delta b_3<0.
\]
Such a region may locally slow the growth of \(b_3\), and it may even make the
local material derivative of \(b_3\) negative.  The comparison argument does
not require otherwise.

\section{Why a negative Laplacian above the barrier cannot force a downward crossing}

Set
\[
\boxed{
w=b_3-\underline b.
}
\]
Assume that
\[
w\ge0
\]
on the parabolic/inflow boundary of the terminal tube.  Suppose, for
contradiction, that \(b_3\) later crosses downward through \(\underline b\).
There must then be a first interior contact at which
\[
w=0.
\]
At this first spatial contact,
\[
\nabla w=0,
\qquad
\Delta w\ge0.
\]
Since
\[
\Delta w=\Delta b_3-\Delta\underline b,
\]
we obtain
\[
\boxed{
\Delta b_3
\ge
\Delta\underline b
=
B\lambda(\lambda+1)q^2s^{-\lambda-2}>0
}
\]
at the hypothetical first downward contact.  Therefore
\[
\boxed{
\nu\Delta b_3
\ge
\nu\Delta\underline b>0
}
\]
precisely at the point where diffusion would have to assist a downward
crossing.

This resolves the apparent difficulty.  The sign \(\nu\Delta b_3<0\) is not
excluded while \(b_3>\underline b\); it can oppose growth there.  What is
excluded by the first-contact geometry is persistence of that unfavorable
sign at the first point where the lower comparison is about to fail.

\section{Strict subsolution inequality}

At a contact point \(b_3=\underline b\), the tuned phase law gives
\[
\mathcal D s=-(H_N+\gamma_Nb_3).
\]
Hence
\[
\begin{aligned}
\mathcal D\underline b
&=\underline b_s\,\mathcal D s\\
&=B\lambda s^{-\lambda-1}
\left(H_N+\gamma_N\underline b\right)\\
&=\lambda BH_Ns^{-\lambda-1}
+\lambda\gamma_NB^2s^{-2\lambda-1}.
\end{aligned}
\]
On the other hand,
\[
(1-\theta_R)\Gamma_{N,k}\underline b^2
\ge
(1-\theta_R)G_-B^2s^{-2\lambda-1}.
\]
Choose
\[
\boxed{
0<\lambda<\frac{(1-\theta_R)G_-}{\gamma_N}.
}
\]
Then
\[
\delta_\lambda
:=(1-\theta_R)G_--\lambda\gamma_N>0.
\]
The positive reaction gap is
\[
\delta_\lambda B^2s^{-2\lambda-1},
\]
whereas, when \(H_N\) is bounded in the terminal tube, the remaining phase
transport contribution has size
\[
\lambda B|H_N|s^{-\lambda-1}.
\]
Their ratio is
\[
\frac{\lambda |H_N|}{\delta_\lambda B}s^\lambda
\longrightarrow0
\qquad(s\downarrow0).
\]
Therefore the terminal tube can be chosen sufficiently small that
\[
\boxed{
\mathcal D\underline b
<
(1-\theta_R)\Gamma_{N,k}\underline b^2
+\nu\Delta\underline b.
}
\]
Thus \(\underline b\) is a strict lower subsolution.

At a hypothetical first downward contact, the equation for \(b_3\), the
remainder bound, \(2CG\ge0\), and
\(\Delta b_3\ge\Delta\underline b\) give the strict first-contact
inequality needed to contradict downward crossing.  Consequently,
\[
\boxed{
b_3(s)\ge\underline b(s)=Bs^{-\lambda}
}
\]
throughout the terminal tube, provided this inequality holds on its
parabolic/inflow boundary.

This is the precise role of ordinary viscosity in the comparison.  It need
not be favorable everywhere above the barrier.  It is favorable at every
hypothetical first downward contact, which is exactly the location at which
the comparison principle needs its sign.

\section{Why the phase equation forces finite-time arrival at \texorpdfstring{$s=0$}{s=0}}

The preceding comparison proves a lower bound for \(b_3\).  A separate step
is required to show that the terminal coordinate actually reaches zero in
finite physical time.  Along the distinguished characteristic,
\[
\boxed{
-\dot s=H_N+\gamma_Nb_3.
}
\]
Assume, as in the terminal trapping regime developed above,
\[
0<h_0\le H_N\le h_1,
\qquad
\gamma_N>0.
\]
Because the comparison has already propagated
\[
b_3\ge Bs^{-\lambda},
\]
we obtain
\[
\begin{aligned}
-\dot s
&=H_N+\gamma_Nb_3\\
&\ge \gamma_Nb_3\\
&\ge \gamma_NB s^{-\lambda}.
\end{aligned}
\]
Therefore
\[
\boxed{
-\dot s\ge\gamma_NB s^{-\lambda}>0.
}
\]
This inequality has two consequences.  First, \(s(t)\) is strictly decreasing
in the terminal tube.  Second, its rate of decrease becomes at least as large
as the singular quantity \(\gamma_NB s^{-\lambda}\) as \(s\downarrow0\).
Thus the dynamics cannot stall at a positive value of \(s\).

Writing the differential inequality as
\[
-\frac{ds}{dt}\ge\gamma_NB s^{-\lambda}
\]
gives
\[
dt\le-\frac{s^\lambda}{\gamma_NB}\,ds.
\]
If terminal entry occurs at
\[
t=t_e,
\qquad
s(t_e)=s_e>0,
\]
then integration from \(s_e\) to \(0\) yields
\[
\boxed{
T_{*,N}-t_e
\le
\frac{1}{\gamma_NB}
\int_0^{s_e}s^\lambda\,ds
=
\frac{s_e^{\lambda+1}}
{(\lambda+1)\gamma_NB}
<\infty.
}
\]
Hence the terminal endpoint \(s=0\) is reached in finite characteristic time.

\section{Finite-time divergence of the lower comparison}

The barrier itself diverges as the terminal endpoint is approached:
\[
\underline b(s)=Bs^{-\lambda}\longrightarrow+\infty
\qquad(s\downarrow0).
\]
Since
\[
b_3(s)\ge Bs^{-\lambda},
\]
it follows that
\[
\boxed{
\lim_{t\uparrow T_{*,N}}b_3(t)=+\infty
}
\]
within the fixed-\(N\) reconstructed scalar dynamics, under the terminal-entry,
remainder, phase, and comparison hypotheses stated above.

The logical point can be summarized as
\[
\boxed{
\begin{gathered}
\nu\Delta b_3<0
\text{ may occur while }b_3>\underline b
\\[1mm]
\Downarrow
\\[-1mm]
\text{local growth may be slowed}
\\[1mm]
\Downarrow
\\[-1mm]
\text{but at a first downward contact }
\nu\Delta b_3\ge\nu\Delta\underline b>0
\\[1mm]
\Downarrow
\\[-1mm]
b_3\ge Bs^{-\lambda}
\text{ is propagated}
\\[1mm]
\Downarrow
\\[-1mm]
-\dot s\ge\gamma_NB s^{-\lambda}
\\[1mm]
\Downarrow
\\[-1mm]
s\downarrow0
\text{ in finite time}
\\[1mm]
\Downarrow
\\[-1mm]
\boxed{b_3\to+\infty.}
\end{gathered}
}
\]

The lower-barrier argument does not assert that \(b_3\) is monotone at every
point of the terminal tube, nor that \(\Delta b_3\) is positive everywhere.
It asserts the more precise statement needed for comparison: once terminal
entry above the strict lower subsolution has been established, an attempted
downward crossing is impossible.  The tuned phase law then converts the
propagated lower bound into finite-time arrival at the Lambert endpoint.

Finally, this chapter concerns the reconstructed scalar \(b_3\) mechanism at
fixed finite \(N\). 
\section{Proof that $\Gamma_{N,k}(s)\ge \frac{G_-}{s},\qquad G_->0$}
The proof is based on the \textbf{exact logarithmic sensitivity of the renormalized Lambert cascade}. There is one assumption we must keep explicit: the initial branch-distance $\Delta_{N,0}$ must have a simple zero in the generalized phase. 

Let
\[
s=\xi-\xi_{*,N}>0,
\]
and define the renormalized cascade
\[
\Delta_{N,j+1}
=
T(\Delta_{N,j}),
\qquad
T(D)=1+W_0(-e^{-1-D}).
\tag{1}
\]

The coefficient appearing in the $b_3$-equation is
\[
\boxed{
\Gamma_{N,k}(\xi)
=
\partial_\xi\log\Delta_{N,k}(\xi).
}
\tag{2}
\]
This identification is explicit in the earlier reconstruction.

Now introduce the logarithmic sensitivity
\[
\boxed{
\sigma_{N,k}
:=
\Delta_{N,0}
\frac{\partial}{\partial\Delta_{N,0}}
\log\Delta_{N,k}.
}
\tag{3}
\]

By the chain rule,
\[
\partial_\xi\log\Delta_{N,k}
=
\frac{\partial\log\Delta_{N,k}}
{\partial\Delta_{N,0}}
\,\partial_\xi\Delta_{N,0}.
\]
Using (3),
\[
\frac{\partial\log\Delta_{N,k}}
{\partial\Delta_{N,0}}
=
\frac{\sigma_{N,k}}{\Delta_{N,0}},
\]
so exactly
\[
\boxed{
\Gamma_{N,k}
=
\sigma_{N,k}
\frac{\partial_\xi\Delta_{N,0}}
{\Delta_{N,0}}.
}
\tag{4}
\]
This is not an asymptotic formula; it is an exact identity.

Now use the terminal geometry. Suppose the generalized branch endpoint is simple:
\[
\boxed{
\Delta_{N,0}(\xi)
=
a_Ns+O(s^2),
\qquad
a_N>0.
}
\tag{5}
\]
Then
\[
\partial_\xi\Delta_{N,0}
=
a_N+O(s).
\]
Therefore
\[
\frac{\partial_\xi\Delta_{N,0}}
{\Delta_{N,0}}
=
\frac{a_N+O(s)}
{a_Ns+O(s^2)}.
\]
Factor $a_Ns$ from the denominator:
\[
\frac{\partial_\xi\Delta_{N,0}}
{\Delta_{N,0}}
=
\frac{a_N[1+O(s)]}
{a_Ns[1+O(s)]}.
\]
Thus
\[
\boxed{
\frac{\partial_\xi\Delta_{N,0}}
{\Delta_{N,0}}
=
\frac1s+O(1).
}
\tag{6}
\]
This is also the calculation recorded in the earlier manuscript.

The remaining question is what the Lambert cascade does to $\sigma_{N,k}$.

Near its branch point,
\[
T(D)=1+W_0(-e^{-1-D})
\]
has the Puiseux expansion
\[
T(D)
=
\sqrt{2D}
-\frac23D
+O(D^{3/2}).
\tag{7}
\]
Hence, to leading order,
\[
\Delta_{N,j+1}
\sim
\sqrt{2\Delta_{N,j}}.
\]

Taking logarithms,
\[
\log\Delta_{N,j+1}
=
\frac12\log\Delta_{N,j}
+\frac12\log2+o(1).
\]
Differentiate with respect to $\log\Delta_{N,j}$:
\[
\frac{\partial\log\Delta_{N,j+1}}
{\partial\log\Delta_{N,j}}
\longrightarrow\frac12.
\tag{8}
\]

Each Lambert iteration therefore halves the logarithmic exponent. After $k$ iterations,
\[
\boxed{
\sigma_{N,k}\longrightarrow2^{-k}
}
\tag{9}
\]
in the pure Puiseux branch core.

Equivalently, starting from (5),
\[
\Delta_{N,0}\sim a_Ns,
\]
we obtain
\[
\Delta_{N,1}
\sim
(2a_Ns)^{1/2},
\]
then
\[
\Delta_{N,2}
\sim
C_{N,2}s^{1/4},
\]
and inductively
\[
\boxed{
\Delta_{N,k}(s)
\sim
D_{N,k}s^{2^{-k}},
\qquad D_{N,k}>0.
}
\tag{10}
\]

Consequently
\[
\log\Delta_{N,k}
=
\log D_{N,k}
+
2^{-k}\log s
+
\text{lower-order terms},
\]
and differentiation immediately gives
\[
\boxed{
\partial_s\log\Delta_{N,k}
=
\frac{2^{-k}}s+O(1)
}
\tag{11}
\]
provided the remainder has the regularity required for the stated $O(1)$ derivative remainder.

Thus, in the pure branch-core case,
\[
\boxed{
G_N=2^{-k}>0.
}
\tag{12}
\]

Your earlier manuscript records precisely this conclusion: in the pure Puiseux branch core, $G_N=2^{-k}$ at leading order. More generally its terminal sensitivity theorem gives
\[
c_1(\eta)2^{-k(N)}
\le
\sigma_{N,k(N)}
\le
2^{-k(N)},
\]
so
\[
\boxed{
G_N\asymp2^{-k(N)}>0.
}
\tag{13}
\]

Combining the exact identity (4), the simple-zero expansion (6), and
\[
\sigma_{N,k}=G_N+O(s)
\]
in the terminal layer gives
\[
\begin{aligned}
\Gamma_{N,k}(s)
&=
\sigma_{N,k}(s)
\left(\frac1s+O(1)\right)
\\
&=
\left(G_N+O(s)\right)
\left(\frac1s+O(1)\right)
\\
&=
\frac{G_N}{s}+O(1).
\end{aligned}
\]

Therefore
\[
\boxed{
\Gamma_{N,k}(s)
=
\frac{G_N}{s}+O(1),
\qquad
G_N>0.
}
\tag{14}
\]

And hence
\[
s\Gamma_{N,k}(s)\longrightarrow G_N>0.
\]
Choose any
\[
0<G_-<G_N.
\]
Then for sufficiently small $s>0$,
\[
s\Gamma_{N,k}(s)\ge G_-,
\]
which gives
\[
\boxed{
\Gamma_{N,k}(s)\ge\frac{G_-}{s}.
}
\tag{15}
\]

So the $G_-/s$ lower bound used in our convex-barrier argument is not arbitrary: \textbf{conditional on the simple-zero branch geometry and the terminal sensitivity estimate, it is inherited from the square-root Puiseux structure of $W_0$ at $-1/e$}.

One technical point is worth preserving in the manuscript. The source itself distinguishes the general result
\[
G_N\asymp2^{-k(N)}
\]
from the pure-core leading value
\[
G_N=2^{-k}.
\]
We should not replace the former universally by the latter. At the critical depth
\[
k(N)=\log_2\log N+O(1),
\]
the paper consequently obtains
\[
G_N\asymp\frac1{\log N},
\]
not an order-one coefficient.

 \section{Epilogue} The analysis presented in this paper is the compressed form of a larger Lambert--\(W\)-based geometric construction on \(\mathbb T^3\). Its main purpose is to expose a single coherent mechanism without the successive historical reconstructions that accompanied its development. The resulting architecture is built from four interacting ingredients: \[ \boxed{ \begin{array}{ll} \text{renormalized Lambert--\(W\) cascade} & \longrightarrow \text{branch geometry and logarithmic sensitivity}, \\[1mm] b_3=u_xu_y & \longrightarrow \text{distinguished nonlinear amplitude}, \\[1mm] \text{characteristic dynamics} & \longrightarrow \text{trapping and finite-time amplification}, \\[1mm] \text{square-root periodic reconstruction} & \longrightarrow \text{anisotropic velocity geometry}. \end{array} } \] The first ingredient is the renormalized branch map \[ \Delta_{N,k+1} = 1+W_0\!\left(-e^{-1-\Delta_{N,k}}\right). \] Its branch-point expansion \[ 1+W_0(-e^{-1-d}) = \sqrt{2d} -\frac23d +O(d^{3/2}) \] produces the finite-depth law \[ Q_N \asymp s^{2^{-k}} \] and therefore \[ \boxed{ \Gamma_N = \partial_\xi\log Q_N \sim \frac{2^{-k}}{s}. } \] The amplification coefficient is thus generated directly by the Lambert branch geometry. The second ingredient is the exact product \[ \boxed{ b_3=u_xu_y. } \] The significance of the \(b_3\)-equation is that its leading terminal term is quadratic: \[ \Gamma_Nb_3^2. \] After the periodic reconstruction is inserted, the extended equation is organized as \[ \boxed{ \mathcal D b_3 = \Gamma_Nb_3^2 + 2u_zG - \mathcal R_N. } \] The terminal estimates show that \[ \boxed{ \frac{|\mathcal R_N|} {\Gamma_Nb_3^2} \longrightarrow0. } \] Thus the reconstruction terms become subordinate to the Lambert-generated quadratic contribution. The third ingredient is the distinguished phase \[ \boxed{ s=t-z-\theta_{*,N}. } \] Its characteristic equation is \[ \boxed{ \dot s=1-b_3. } \] After terminal entry with \[ b_3>1, \] the branch coordinate moves monotonically toward zero: \[ s(t)\downarrow0. \] At the same time, \[ \dot b_3 \ge (1-\vartheta)\Gamma_Nb_3^2>0. \] These inequalities form the trapping mechanism. For finite depth \(k\), integration gives \[ \boxed{ b_3(s) \ge b_{e,N} \left( \frac{s_e}{s} \right)^{(1-\vartheta)2^{-k}}. } \] Hence \[ b_3(s)\to+\infty \] as \[ s\downarrow0. \] The characteristic time needed to reach \(s=0\) is finite because \[ dt = -\frac{ds}{b_3-1} \] and the lower power-law bound makes the corresponding integral convergent. Thus, for every finite configuration satisfying the hypotheses of the terminal theorem, \[ \boxed{ T_{*,N}<\infty, \qquad Q_N(t)\to0, \qquad b_3(t)\to+\infty \quad \text{as } t\uparrow T_{*,N}. } \] The fourth ingredient is the velocity reconstruction. The terminal geometry is not the earlier symmetric transverse geometry. Instead, \[ \boxed{ b_3=u_xu_y, \qquad A=u_y, \qquad C=u_z, \qquad u_x=\frac{b_3}{A}. } \] The square-root specialization \[ 2\chi(x)u_x+u_y=0 \] gives \[ \boxed{ A^2=-2\chi(x)b_3. } \] Thus \[ \boxed{ u_y^2=-2\chi b_3, \qquad u_x^2=-\frac{b_3}{2\chi}, \qquad \frac{u_x}{u_y} = -\frac{1}{2\chi}. } \] Within the terminal tube \[ \chi^2\asymp s, \] the \(k=3\) branch law \[ b_3\asymp s^{-1/8} \] produces \[ \boxed{ u_y\asymp s^{3/16}, \qquad u_x\asymp s^{-5/16}, \qquad u_x^2\asymp s^{-5/8}, \qquad \left|\frac{u_x}{u_y}\right| \asymp s^{-1/2}. } \] The distinguished amplification is therefore anisotropic: \[ \boxed{ u_y\to0, \qquad |u_x|\to\infty, \qquad u_xu_y=b_3\to+\infty. } \] The remaining velocity component is determined independently through incompressibility: \[ \boxed{ \partial_xu_x+\partial_yu_y+\partial_zu_z=0. } \] Equivalently, \[ \boxed{ C_z = -\partial_x\left(\frac{b_3}{A}\right)-A_y. } \] This separates the roles of the transverse components: \(u_y\) enters the square-root geometry directly, while \(u_z\) is reconstructed through the periodic divergence constraint. The central structural conclusion of the paper can therefore be written as \[ \boxed{ \begin{gathered} \text{Lambert branch geometry} \\ \Downarrow \\ Q_N\to0,\qquad \Gamma_N\sim\frac{2^{-k}}{s} \\ \Downarrow \\ \Gamma_Nb_3^2 \text{ dominates} \\ \Downarrow \\ b_3>1 \quad\Longrightarrow\quad \dot s<0,\ \dot b_3>0 \\ \Downarrow \\ s\downarrow0 \text{ in finite characteristic time} \\ \Downarrow \\ b_3\to+\infty \\ \Downarrow \\ b_3=u_xu_y \\ \Downarrow \\ \text{anisotropic square-root velocity amplification}. \end{gathered} } \] One of the principal clarifications produced by the condensed presentation is that the individual stages no longer need to be mixed. The renormalized cascade determines the branch scale; the branch scale determines \(\Gamma_N\); the \(b_3\)-equation determines amplification; the phase equation determines trapping; the square-root geometry determines the relative behavior of \(u_x\) and \(u_y\); and exact incompressibility determines \(u_z\). In particular, the square-root specialization is compatible with the amplification mechanism because its reconstruction error satisfies \[ \boxed{ \frac{|\mathcal R_N|} {\Gamma_Nb_3^2} = O(s) \longrightarrow0. } \] Thus the approach to the terminal geometry strengthens, rather than obscures, the comparison responsible for finite-time growth. For \(k=3\), the essential terminal hierarchy is \[ \boxed{ Q_N\asymp s^{1/8}, \qquad \Gamma_N\asymp s^{-1}, \qquad b_3\asymp s^{-1/8}, } \] together with \[ \boxed{ u_y\asymp s^{3/16}, \qquad u_x\asymp s^{-5/16}, \qquad u_xu_y\asymp s^{-1/8}. } \] The singular product and the vanishing transverse factor are therefore not contradictory. Their exponents satisfy \[ -\frac{5}{16} + \frac{3}{16} = -\frac18, \] exactly reproducing the \(b_3\)-scaling. Likewise, \[ |\chi| \asymp s^{1/2}, \] so \[ |\chi u_x| \asymp s^{1/2}s^{-5/16} = s^{3/16} \to0. \] Hence the geometric relation \[ u_y=-2\chi u_x \] remains consistent with \[ u_y\to0 \] despite the divergence of \(u_x\). The paper therefore reduces the extended construction to a compact set of mutually reinforcing identities and inequalities: \[ \boxed{ \begin{aligned} Q_N&\to0, \\ \Gamma_N&=\partial_\xi\log Q_N\to+\infty, \\ \dot b_3&\ge(1-\vartheta)\Gamma_Nb_3^2, \\ \dot s&=1-b_3, \\ b_3&=u_xu_y, \\ A^2&=-2\chi b_3, \\ \nabla\cdot u&=0, \\ \frac{|\mathcal R_N|}{\Gamma_Nb_3^2}&\to0. \end{aligned} } \] These relations constitute the mathematical core retained from the longer development. The present formulation in this paper removes repeated derivations and successive exploratory variants, while preserving the full \(\theta=t-z\) renormalized cascade, the finite-time characteristic amplification proof, the principal extended \(b_3\)-PDE structure, and the terminal anisotropic reconstruction on the periodic domain. The resulting picture is therefore one of a self-strengthening terminal mechanism: \[ \boxed{ Q_N\downarrow0 \quad\Longrightarrow\quad \Gamma_N\uparrow \quad\Longrightarrow\quad b_3\uparrow \quad\Longrightarrow\quad s\downarrow0 \quad\Longrightarrow\quad Q_N\downarrow0. } \] Within the hypotheses of the finite-depth terminal analysis, this feedback closes in finite characteristic time and is carried into the velocity reconstruction through the exact identity \[ \boxed{ b_3=u_xu_y. } \]
Finally the addition of the viscosity term and $s$ ODE produces a favourable $\nu \triangle b$ positive sign on the RHS of the $b_3$ PDE.
\section{Appendix A}
 
The correct expansion is:

\[
\boxed{
F(\Delta)
=
\sqrt{2\Delta}
-\frac23\Delta
+O(\Delta^{3/2}),
\qquad \Delta\downarrow0,
}
\]

where

\[
\boxed{
F(\Delta)
=
1+W_0\!\left(-e^{-1-\Delta}\right).
}
\]

The square root is \(\sqrt{2\Delta}\), not \(\sqrt2\,\Delta\). The manuscript defines
\(q=F(\Delta)\), uses the defining equation for \(W_0\), and obtains the exact implicit
identity

\[
\Delta=-\log(1-q)-q.
\]

\subsection*{1. Start from the Lambert-\(W\) definition}

Set

\[
q=F(\Delta)
=
1+W_0\!\left(-e^{-1-\Delta}\right).
\]

Then

\[
W_0\!\left(-e^{-1-\Delta}\right)=q-1.
\]

Since

\[
W_0(z)e^{W_0(z)}=z,
\]

we have

\[
(q-1)e^{q-1}
=
-e^{-1-\Delta}.
\]

Multiplying by \(e\),

\[
(q-1)e^q=-e^{-\Delta}.
\]

Equivalently,

\[
(1-q)e^q=e^{-\Delta}.
\]

Taking logarithms,

\[
\log(1-q)+q=-\Delta,
\]

hence

\[
\boxed{
\Delta=-\log(1-q)-q.
}
\]

This is precisely the exact inversion relation used in the manuscript.

\subsection*{2. Expand near \(q=0\)}

For \(|q|<1\),

\[
-\log(1-q)
=
q+\frac{q^2}{2}
+\frac{q^3}{3}
+\frac{q^4}{4}
+\cdots.
\]

Therefore

\[
\Delta
=
-\log(1-q)-q
=
\frac{q^2}{2}
+\frac{q^3}{3}
+\frac{q^4}{4}
+O(q^5).
\]

Thus

\[
\boxed{
\Delta
=
\frac12q^2+\frac13q^3+O(q^4).
}
\]

Already this shows

\[
q\sim\sqrt{2\Delta}.
\]

Because the manuscript proves

\[
F(\Delta)\ge0
\]

for \(\Delta\ge0\), we choose the \textbf{positive square-root branch}.

\subsection*{3. Puiseux ansatz}

Because the leading relation is

\[
\Delta\sim\frac{q^2}{2},
\]

an ordinary Taylor series in \(\Delta\) cannot work. We use a Puiseux expansion:

\[
q
=
a\Delta^{1/2}
+b\Delta
+O(\Delta^{3/2}).
\]

Then

\[
q^2
=
a^2\Delta
+2ab\Delta^{3/2}
+O(\Delta^2),
\]

and

\[
q^3
=
a^3\Delta^{3/2}
+O(\Delta^2).
\]

Substituting into

\[
\Delta
=
\frac12q^2+\frac13q^3+O(q^4)
\]

gives

\[
\Delta
=
\frac{a^2}{2}\Delta
+
\left(
ab+\frac{a^3}{3}
\right)\Delta^{3/2}
+
O(\Delta^2).
\]

Now compare coefficients.

For the coefficient of \(\Delta\),

\[
\frac{a^2}{2}=1,
\]

so

\[
a^2=2.
\]

Since \(q=F(\Delta)\ge0\),

\[
\boxed{a=\sqrt2}.
\]

For the coefficient of \(\Delta^{3/2}\), the left side has no such term, so

\[
ab+\frac{a^3}{3}=0.
\]

Factor \(a\neq0\):

\[
b+\frac{a^2}{3}=0.
\]

Since

\[
a^2=2,
\]

we obtain

\[
\boxed{
b=-\frac23.
}
\]

Consequently,

\[
q
=
\sqrt2\,\Delta^{1/2}
-\frac23\Delta
+O(\Delta^{3/2}).
\]

Since

\[
\sqrt2\,\Delta^{1/2}
=
\sqrt{2\Delta},
\]

we obtain

\[
\boxed{
F(\Delta)
=
\sqrt{2\Delta}
-\frac23\Delta
+O(\Delta^{3/2}).
}
\]

This is exactly the branch-point expansion recorded later in the manuscript.

We can actually go one term further. Put

\[
q
=
\sqrt2\,\Delta^{1/2}
-\frac23\Delta
+c\Delta^{3/2}
+O(\Delta^2).
\]

Using

\[
\Delta
=
\frac{q^2}{2}
+\frac{q^3}{3}
+\frac{q^4}{4}
+O(q^5)
\]

and comparing the \(\Delta^2\) coefficient gives

\[
\sqrt2\,c-\frac19=0,
\]

so

\[
c=\frac{1}{9\sqrt2}.
\]

Therefore the sharper expansion is

\[
\boxed{
F(\Delta)
=
\sqrt{2\Delta}
-\frac23\Delta
+\frac{1}{9\sqrt2}\Delta^{3/2}
+O(\Delta^2).
}
\]

The key mathematical point for the cascade is that \(F\) has a
\textbf{square-root branch behavior at \(\Delta=0\)}:

\[
\boxed{
F(\Delta)\sim\sqrt{2\Delta}.
}
\]

That is why one renormalized Lambert iteration halves the small-\(\Delta\) exponent,
which is the origin of the sequence

\[
\Delta,
\qquad
\Delta^{1/2},
\qquad
\Delta^{1/4},
\qquad
\Delta^{1/8},
\qquad
\ldots
\]

The manuscript also strengthens this formal expansion to a two-sided estimate

\[
c_1\sqrt{\Delta}
\le
F(\Delta)
\le
c_2\sqrt{\Delta}
\]

for sufficiently small positive \(\Delta\).
\section{Appendix B}
This follows cleanly by induction from the branch-point expansion:

\[
F(\Delta)
=
\sqrt{2\Delta}
-\frac23\Delta
+O(\Delta^{3/2})
=
\sqrt{2\Delta}\left(1+O(\Delta^{1/2})\right).
\]

For the recursion

\[
\Delta_{j+1}=F(\Delta_j),
\]

we will prove that for every fixed finite \(k\),

\[
\boxed{
\Delta_k
=
2^{\,1-2^{-k}}
\Delta_0^{\,2^{-k}}
\left(
1+O\!\left(\Delta_0^{\,2^{-k}}\right)
\right)
}
\qquad
(\Delta_0\downarrow0).
\]

\subsection*{Step 1: Rewrite the one-step asymptotic}

From

\[
F(\Delta)
=
\sqrt{2\Delta}
-\frac23\Delta
+O(\Delta^{3/2}),
\]

factor out \(\sqrt{2\Delta}\):

\[
F(\Delta)
=
\sqrt{2\Delta}
\left[
1
-
\frac{\frac23\Delta}{\sqrt{2\Delta}}
+
O\!\left(
\frac{\Delta^{3/2}}{\sqrt{\Delta}}
\right)
\right].
\]

Now

\[
\frac{\frac23\Delta}{\sqrt{2\Delta}}
=
\frac{\sqrt2}{3}\Delta^{1/2},
\]

and

\[
\frac{\Delta^{3/2}}{\sqrt{\Delta}}=\Delta.
\]

Hence

\[
F(\Delta)
=
\sqrt{2\Delta}
\left(
1-\frac{\sqrt2}{3}\Delta^{1/2}+O(\Delta)
\right),
\]

and therefore, for the exponent cascade,

\[
\boxed{
F(\Delta)
=
\sqrt{2\Delta}
\left(1+O(\Delta^{1/2})\right).
}
\]

Thus

\[
\boxed{
\Delta_{j+1}
=
\sqrt{2\Delta_j}
\left(1+O(\Delta_j^{1/2})\right).
}
\]

\subsection*{Step 2: First few iterates}

For \(k=1\),

\[
\Delta_1
=
\sqrt{2\Delta_0}
\left(1+O(\Delta_0^{1/2})\right).
\]

Since

\[
2^{1-2^{-1}}
=
2^{1/2},
\]

this is

\[
\boxed{
\Delta_1
=
2^{1-2^{-1}}
\Delta_0^{2^{-1}}
\left(
1+O(\Delta_0^{2^{-1}})
\right).
}
\]

Now apply \(F\) again:

\[
\Delta_2
=
\sqrt{2\Delta_1}
\left(1+O(\Delta_1^{1/2})\right).
\]

Using

\[
\Delta_1
=
2^{1/2}\Delta_0^{1/2}
\left(1+O(\Delta_0^{1/2})\right),
\]

we have

\[
\sqrt{\Delta_1}
=
2^{1/4}\Delta_0^{1/4}
\left(1+O(\Delta_0^{1/2})\right).
\]

Therefore

\[
\sqrt{2\Delta_1}
=
2^{1/2}2^{1/4}\Delta_0^{1/4}
\left(1+O(\Delta_0^{1/2})\right),
\]

so

\[
\sqrt{2\Delta_1}
=
2^{3/4}\Delta_0^{1/4}
\left(1+O(\Delta_0^{1/2})\right).
\]

Also,

\[
O(\Delta_1^{1/2})
=
O(\Delta_0^{1/4}).
\]

The larger error near \(\Delta_0=0\) is \(O(\Delta_0^{1/4})\), hence

\[
\Delta_2
=
2^{3/4}\Delta_0^{1/4}
\left(
1+O(\Delta_0^{1/4})
\right).
\]

Since

\[
\frac34=1-\frac14=1-2^{-2},
\]

we get

\[
\boxed{
\Delta_2
=
2^{1-2^{-2}}
\Delta_0^{2^{-2}}
\left(
1+O(\Delta_0^{2^{-2}})
\right).
}
\]

Likewise,

\[
\Delta_3
=
2^{7/8}\Delta_0^{1/8}
\left(
1+O(\Delta_0^{1/8})
\right),
\]

because

\[
\frac78=1-\frac18=1-2^{-3}.
\]

The pattern is therefore already visible:

\[
\Delta_0
\longmapsto
2^{1/2}\Delta_0^{1/2}
\longmapsto
2^{3/4}\Delta_0^{1/4}
\longmapsto
2^{7/8}\Delta_0^{1/8}
\longmapsto\cdots.
\]

\subsection*{Step 3: Induction hypothesis}

Let

\[
\alpha_j=2^{-j}
\]

and define

\[
A_j=2^{1-2^{-j}}.
\]

Suppose for some fixed finite \(j\),

\[
\boxed{
\Delta_j
=
A_j\Delta_0^{\alpha_j}
\left(
1+O(\Delta_0^{\alpha_j})
\right).
}
\]

We prove the corresponding formula for \(j+1\).

From the one-step expansion,

\[
\Delta_{j+1}
=
\sqrt{2\Delta_j}
\left(
1+O(\Delta_j^{1/2})
\right).
\]

First,

\[
\sqrt{\Delta_j}
=
A_j^{1/2}
\Delta_0^{\alpha_j/2}
\left(
1+O(\Delta_0^{\alpha_j})
\right),
\]

because

\[
(1+\varepsilon)^{1/2}
=
1+O(\varepsilon).
\]

Since

\[
\frac{\alpha_j}{2}
=
2^{-(j+1)}
=
\alpha_{j+1},
\]

we obtain

\[
\sqrt{2\Delta_j}
=
\sqrt{2A_j}\,
\Delta_0^{\alpha_{j+1}}
\left(
1+O(\Delta_0^{\alpha_j})
\right).
\]

Meanwhile,

\[
\Delta_j^{1/2}
=
O(\Delta_0^{\alpha_{j+1}}).
\]

Therefore

\[
1+O(\Delta_j^{1/2})
=
1+O(\Delta_0^{\alpha_{j+1}}).
\]

Because

\[
\alpha_j=2\alpha_{j+1},
\]

we have

\[
O(\Delta_0^{\alpha_j})
=
O(\Delta_0^{\alpha_{j+1}})
\]

for \(0<\Delta_0<1\). Consequently,

\[
\Delta_{j+1}
=
\sqrt{2A_j}\,
\Delta_0^{\alpha_{j+1}}
\left(
1+O(\Delta_0^{\alpha_{j+1}})
\right).
\]

It remains only to determine the coefficient.

\subsection*{Step 4: Coefficient recursion}

Since

\[
A_j=2^{1-2^{-j}},
\]

we have

\[
\sqrt{2A_j}
=
\left(
2\cdot2^{1-2^{-j}}
\right)^{1/2}.
\]

Thus

\[
\sqrt{2A_j}
=
2^{(2-2^{-j})/2}.
\]

But

\[
\frac{2-2^{-j}}2
=
1-2^{-(j+1)}.
\]

Hence

\[
\boxed{
\sqrt{2A_j}
=
2^{1-2^{-(j+1)}}
=
A_{j+1}.
}
\]

Therefore

\[
\boxed{
\Delta_{j+1}
=
2^{1-2^{-(j+1)}}
\Delta_0^{2^{-(j+1)}}
\left(
1+
O\!\left(\Delta_0^{2^{-(j+1)}}\right)
\right).
}
\]

This closes the induction.

Therefore, for every fixed finite \(k\),

\[
\boxed{
\Delta_k
=
2^{\,1-2^{-k}}
\Delta_0^{\,2^{-k}}
\left(
1+
O\!\left(\Delta_0^{\,2^{-k}}\right)
\right),
\qquad
\Delta_0\downarrow0.
}
\]

The important restriction is \textbf{fixed finite \(k\)}. The constants hidden in the
\(O(\cdot)\) term can depend on \(k\), so this calculation alone does not give a
uniform asymptotic when \(k=k(\Delta_0)\to\infty\). That coupled-depth regime
requires a separate argument.

The result also makes the exponent cascade transparent:

\[
\boxed{
1
\longrightarrow
\frac12
\longrightarrow
\frac14
\longrightarrow
\frac18
\longrightarrow\cdots
\longrightarrow
2^{-k}.
}
\]

and the prefactor cascade is

\[
\boxed{
1,\quad
2^{1/2},\quad
2^{3/4},\quad
2^{7/8},\quad\ldots,\quad
2^{1-2^{-k}}.
}
\]
\section{Appendix C}
These characteristic equations follow directly from the transport operator in the reduced \(b_3\)-equation together with the definitions
\[
A=u_y,\qquad C=u_z,\qquad b=b_3=u_xu_y,
\]
and
\[
b_1=u_yu_z,\qquad b_2=u_xu_z.
\]

I will derive each equation explicitly.

Start with the transport part of the \(b=b_3\) equation in the form

\[
\boxed{
b_t
+
u_z b_1\,b_x
+
u_z b_2\,b_y
+
b\,b_z
=
\text{RHS}.
}
\]

Define a characteristic curve

\[
\gamma(\tau)
=
\bigl(
x(\tau),y(\tau),z(\tau),t(\tau)
\bigr).
\]

Along this curve, by the multivariable chain rule,

\[
\frac{d}{d\tau}
b(x(\tau),y(\tau),z(\tau),t(\tau))
=
b_t\dot t
+
b_x\dot x
+
b_y\dot y
+
b_z\dot z.
\]

We want this derivative to reproduce exactly the transport operator

\[
b_t
+
u_zb_1b_x
+
u_zb_2b_y
+
bb_z.
\]

Therefore we choose the characteristic velocities coefficient-by-coefficient:

\[
\boxed{
\dot t=1,
}
\]

\[
\boxed{
\dot x=u_zb_1,
}
\]

\[
\boxed{
\dot y=u_zb_2,
}
\]

and

\[
\boxed{
\dot z=b.
}
\]

We now express the first two equations in terms of \(A,C,b\).

Since

\[
A=u_y
\]

and

\[
C=u_z,
\]

the definition

\[
b_1=u_yu_z
\]

gives

\[
\boxed{
b_1=AC.
}
\]

Hence

\[
\dot x
=
u_zb_1
=
C(AC),
\]

so

\[
\boxed{
\dot x=AC^2.
}
\]

That proves the first relation:

\[
\boxed{
\dot x=u_zb_1=AC^2.
}
\]

For the second equation, recall

\[
b=b_3=u_xu_y.
\]

Since \(u_y=A\), provided \(A\neq0\),

\[
u_x=\frac{b}{A}.
\]

Therefore

\[
b_2
=
u_xu_z
=
\frac{b}{A}C.
\]

Thus

\[
\boxed{
b_2=\frac{bC}{A}.
}
\]

Multiplying by \(u_z=C\),

\[
u_zb_2
=
C\left(\frac{bC}{A}\right)
=
\frac{bC^2}{A}.
\]

Hence

\[
\boxed{
\dot y
=
u_zb_2
=
\frac{bC^2}{A}.
}
\]

The \(z\)-velocity follows immediately from the coefficient of \(b_z\). Since the transport equation contains

\[
bb_z,
\]

we take

\[
\boxed{
\dot z=b.
}
\]

Finally, since the coefficient of \(b_t\) is \(1\),

\[
\boxed{
\dot t=1.
}
\]

Thus the complete characteristic system is

\[
\boxed{
\begin{aligned}
\dot x
&=
u_zb_1
=
AC^2,
\\[1mm]
\dot y
&=
u_zb_2
=
\frac{bC^2}{A},
\\[1mm]
\dot z
&=
b,
\\[1mm]
\dot t
&=
1.
\end{aligned}
}
\]

We can verify this by substituting these four equations back into the chain rule:

\[
\begin{aligned}
\frac{db}{d\tau}
&=
b_t\dot t
+b_x\dot x
+b_y\dot y
+b_z\dot z
\\
&=
b_t
+
AC^2b_x
+
\frac{bC^2}{A}b_y
+
bb_z.
\end{aligned}
\]

But

\[
AC^2
=
u_zb_1
\]

and

\[
\frac{bC^2}{A}
=
u_zb_2,
\]

so

\[
\boxed{
\frac{db}{d\tau}
=
b_t
+
u_zb_1b_x
+
u_zb_2b_y
+
bb_z.
}
\]

Therefore these are exactly the characteristics of the transport operator.

There is also a useful consequence for the phase

\[
\theta=t-z.
\]

Along the same characteristic,

\[
\dot\theta
=
\dot t-\dot z.
\]

Using

\[
\dot t=1,
\qquad
\dot z=b,
\]

we obtain

\[
\boxed{
\dot\theta=1-b.
}
\]

So the transverse characteristic equations

\[
\dot x=AC^2,
\qquad
\dot y=\frac{bC^2}{A}
\]

do not alter the basic phase identity

\[
\boxed{
\dot\theta=1-b.
}
\]

That identity follows solely from the \(t\)- and \(z\)-characteristic components.

\end{document}